\documentclass[centertags,leqno,11pt]{article}

\usepackage{amsmath}
\usepackage{amssymb}
\usepackage{latexsym}
\usepackage{amsxtra}
\usepackage{amscd}
\usepackage{theorem}

\usepackage{graphics}
\usepackage{epic}

\theoremstyle{change}

{\theorembodyfont{\slshape}

\newtheorem{thm}{Theorem.}[section]
\newtheorem{cor}[thm]{Corollary.}
\newtheorem{lem}[thm]{Lemma.}
\newtheorem{prop}[thm]{Proposition.}

\newtheorem{defn}[thm]{Definition.}}

{\theorembodyfont{\rmfamily}

\newtheorem{rem}[thm]{Remark.}

\newtheorem{ass}[thm]{Assumption.}
\newtheorem{notation}[thm]{Notation.}

}

\renewcommand{\em}{\sl}

\newcommand{\proof}{\noindent {\bf Proof:\ }}
\newcommand{\Endproof}{\hspace*{\fill} $\Box$ \vspace{1ex} \noindent }

\makeatletter
\renewcommand{\subsection}{\@startsection{subsection}{2}%
{\z@}{-3.25ex plus -1ex minus-.2ex}{-1em}{\bf}} \makeatother

\newcommand{\PP}{\mathbb{P}}
\newcommand{\ZZ}{\mathbb{Z}}
\newcommand{\CC}{\mathbb{C}}

\newcommand{\QQ}{\mathbb{Q}}
\newcommand{\NN}{\mathbb{N}}
\newcommand{\FF}{\mathbb{F}}
\renewcommand{\AA}{\mathbb{A}}
\newcommand{\GG}{\mathbb{G}}
\newcommand{\F}{\mathcal{F}}

\newcommand{\V}{\mathcal{V}}
\newcommand{\W}{\mathcal{W}}
\newcommand{\U}{\mathcal{U}}

\newcommand{\LL}{\mathcal{L}}
\newcommand{\OO}{\mathcal{O}}

\newcommand{\tr}{{\rm trace}}
\newcommand{\eps}{{\epsilon}}
\newcommand{\Kd}{{\rm K3}}
\newcommand{\geom}{{\rm geo}}

\newcommand{\Pic}{{\rm Pic}}

\newcommand{\g}{{\bf g}}

\newcommand{\GL}{{\rm GL}}
\newcommand{\SL}{{\rm SL}}
\newcommand{\PGL}{{\rm PGL}}
\newcommand{\PSL}{{\rm PSL}}
\newcommand{\Gal}{{\rm Gal}}
\newcommand{\Aut}{{\rm Aut}}
\newcommand{\Spec}{{\rm Spec\,}}

\newcommand{\Hom}{{\rm Hom}}
\renewcommand{\O}{{\rm O}}

\newcommand{\pib}{\bar{\pi}}

\newcommand{\MC}{{\rm MC}}
\newcommand{\SU}{{\rm SU}}
\newcommand{\GU}{{\rm GU}}
\newcommand{\bF}{{\frak F}}
\newcommand{\bG}{{\frak G}}

\newcommand{\para}{_p}

\newcommand{\et}{{\rm\acute{e}t}}

\newcommand{\an}{^{\rm an }}

\newcommand{\N}{{\rm N }}

\newcommand{\Sh}{{\rm Sh}}
\newcommand{\LS}{{\rm LS}}
\newcommand{\Rep}{{\rm Rep}}
\newcommand{\ga}{{\gamma}}
\newcommand{\NS}{{\rm NS}}

\newcommand{\inj}{\hookrightarrow}
\newcommand{\To}{\;\longrightarrow\;}

\newcommand{\liso}{\;\stackrel{\sim}{\longrightarrow}\;}
\newcommand{\Mapsto}{\;\longmapsto\;}

\newcommand{\KKK}{{\cal K}}

\newcommand{\im}{{\rm im}}
\newcommand{\sym}{{\rm sym}}
\newcommand{\al}{{\alpha}}
\newcommand{\A}{{\cal A}}
\newcommand{\D}{{\bf w}}
\newcommand{\uo}{{{\bf u}}}

\newcommand{\vo}{{{\bf v}}}

\newcommand{\q}{{{\rm d}}}

\newcommand{\pr}{{{\rm pr}}}
\newcommand{\opr}{{\overline{\rm pr}}}
\newcommand{\tim}{{\circ}}

\newcommand{\Conv}{{\rm Conv}}
\newcommand{\rk}{{\rm rk}}

\newcommand{\Var}{{\rm Var}}
\newcommand{\Constr}{{\rm Constr}^\et}

\newcommand{\Frob}{{\rm Frob}}

\newcommand{\sep}{{\rm sep}}
\newcommand{\chara}{{\rm char}}

\newcommand{\PPP}{{\rm P}}

\newcommand{\FFF}{{\cal F}}\newcommand{\GGG}{{\cal G}}
\newcommand{\ka}{{k}}
\newcommand{\bk}{{\bar{k}}}

\newcommand{\G}{{\cal  G}}

\newcommand{\nr}{{\rm  nr}}
\newcommand{\cris}{{\rm  cris}}
\newcommand{\End}{{\rm  End}}

\newcommand{\m}{{\frak m}}
\newcommand{\SO}{{\rm SO}}
\renewcommand{\char}{{\rm char}}

\renewcommand{\L}{{\cal L}}
\newcommand{\CX}{{\frak X}}

\numberwithin{equation}{subsection}
\numberwithin{thm}{subsection}
\theoremstyle{plain}

\begin{document}
\thispagestyle{empty}
\begin{center}${}^{}$
\begin{Large}

\vspace{1cm}
{\bf Habilitationsschrift}\\

\vspace{1cm}

eingereicht bei der \\

\vspace{1cm}

Fakult\"at f\"ur Mathematik und Informatik\\
\vspace{.5cm}
 der\\
\vspace{.5cm}
  Ruprecht-Karls-Universit\"at Heidelberg
 \vspace{6cm}

Vorgelegt von\\
\vspace{.5cm}
Dr. rer. nat. Michael Dettweiler\\
\vspace{.5cm}
 aus T\"ubingen\\
\vspace{1cm}
2005

\pagebreak
\thispagestyle{empty} ${}^{}$

\end{Large}
\pagebreak
 \thispagestyle{empty}

${}^{}$
 \vspace{4cm}

\begin{huge}
{ {\sc Galois realizations of classical groups\\ 
\vspace{.5cm}
 and the middle convolution}}
\end{huge}

\begin{LARGE}

\vspace{1cm}

{\sc by

\vspace{1cm}

Michael Dettweiler}

\thispagestyle{empty}
\pagebreak
 \thispagestyle{empty}

${}^{}$\end{LARGE}\end{center}
\pagebreak
\thispagestyle{empty}
${}^{}$
\vspace{2cm}

\begin{large}
\begin{center}{\bf Abstract}\end{center}
\vspace{.3cm}

Up to date, the middle convolution is the 
most powerful tool for realizing classical groups as Galois groups
over $\QQ(t).$
We study the middle convolution of local systems on the punctured
affine line in the setting of singular cohomology 
and in the setting of \'etale cohomology. We derive a
formula to compute the topological monodromy of the middle convolution
in the general case and use it to deduce  some 
irreducibility criteria. 
Then we give a geometric
interpretation of the middle convolution in the \'etale setting. 

This geometric approach to the convolution 
and the theory 
of Hecke characters yields information
on the occurring arithmetic determinants. We employ these methods to
realize  special linear groups regularly as
Galois groups over $\QQ(t).$

The geometric theory of the middle convolution can also be used 
to compute Frobenius elements for many of the known Galois realizations 
of classical groups. 
We illustrate this by investigating specializations 
of $\PGL_2(\FF_\ell)$-extensions
of $\QQ(t)$ which are  associated to a family of $\Kd$-surfaces of Picard
number $19.$
\pagebreak
 \thispagestyle{empty}

${}^{}$
\end{large}
\pagebreak
 \thispagestyle{empty}

 \thispagestyle{empty}
\tableofcontents \thispagestyle{empty}

\pagebreak


\addcontentsline{toc}{section}{Introduction}
\setcounter{page}{1}
\section*{Introduction}\label{Introduction}

The convolution  
$$ f\ast g \hspace{.05cm}(y) :=\int f(x)g(y-x) dx$$
of sufficiently integrable functions
plays an important role in many areas of mathematics and physics, see
 \cite{Bauer},  \cite{Katz96}, \cite{Schwartz}, and 
 \cite{Yosida}.\\

Suppose that 
 the functions $f$ and $g$ are solutions 
of meromorphic connections on the complex affine line $\AA^1.$
Then $f$ and $g$ can be viewed as sections of the solution sheaves of 
these connections,
see  \cite{Deligne70}.
The comparison theorem between singular cohomology
and de Rham cohomology 
shows that  the integral 
$$\int_\gamma f(x)g(y_0-x) dx, \quad {\rm where}\quad  y_0 \in \AA^1\, ,$$ can be seen 
as a cohomology class in the sheaf cohomology on $\AA^1$
($\gamma$ denoting some homology cycle),
see  \cite{Deligne70} and 
\cite{BlochEsnault}.
 This suggests the following 
generalization from the convolution of functions  to 
the convolution of sheaves  $\V_1,\,\V_2$
on  $\AA^1$:
$$\V_1\ast_{\rm aff} \V_2:=R^1 \pr_{2*}(\V_1\circ \V_2)\, .$$
Here, 
$\pr_2:\AA^2\to \AA^1$  denotes the second projection of 
$\AA^2=\AA^1_x\times \AA^1_y,$ 
$R^1\pr_{2*}$ is the first higher direct image of the functor 
$\pr_{2*}$ (see \cite{Hartshorne} and 
Section~\ref{convsheav}),
 and 
$$\V_1\circ \V_2 :=\pr_1^*(\V_1)\otimes \q^*(\V_2)\, ,$$
where $\pr_1:\AA^2\to \AA^1$ is the first projection and 
$\q:\AA^2\to \AA^1$
 is the difference map $(x,y)\mapsto y-x.$ \\

Let us assume that $\V_1$ and $\V_2$  are sheaves on $\AA^1$ which 
are pushforwards of  local systems on  open subsets of $\AA^1$
(see \cite{Deligne70} and 
\cite{dw03} for the notion of a local system). Such sheaves 
naturally arise as the sheaves of solutions of 
connections on $\AA^1,$ see \cite{Deligne70}. Then
the  ``affine'' convolution $\V_1\ast_{\rm aff} \V_2$  
contains 
a canonical subsheaf $$\im\left( R^1 \pr_{2!}(\V_1\circ \V_2)
\To R^1 \pr_{2*}(\V_1\circ \V_2)\right)\, ,$$ where 
$R^1 \pr_{2!}$
denotes the first higher direct image with compact supports. 
The latter sheaf 
 is canonically isomorphic to
$R^1\overline{\pr}_{2*}(j_*(\V_1\circ \V_2)),$ where 
$j:\AA^2\to \PP^1\times \AA^1$ denotes the natural inclusion and 
 $\overline{\pr}_{2*}:\PP^1\times \AA^1\to \AA^1$ is the second 
projection, see \cite{dw03}.
 We set 
\begin{equation}\label{eqmidconvo}
\V_1\ast \V_2:= R^1\overline{\pr}_{2*}(j_*(\V_1\circ \V_2))\,.
\end{equation}
Following Katz \cite{Katz96}, we call the sheaf 
$\V_1\ast \V_2$ the {\em middle convolution}
of $\V_1$ and $\V_2.$ (In loc.~cit., Katz gives a similar construction 
in a more general category of complexes of sheaves.)
 The term {\em middle} indicates the fact that
there is a natural interpretation of $\V_1\ast \V_2$ as 
a middle direct image, i.e.,   one obtains 
the middle convolution by taking 
the 
image of the cohomology with compact supports in the cohomology
(compare to \cite{Katz96}, Chap. 2, and  \cite{dw03}).  

One reason why one is interested in the middle 
 convolution is that 
 $\V_1\ast \V_2$ is 
often irreducible, while $\V_1\ast_{\rm aff} \V_2$ is not. 
A striking application of  the middle 
convolution is Katz' existence algorithm for irreducible 
rigid local systems, see 
\cite{Katz96}, Chap. 6. See also \cite{dr00}, \cite{Voelklein01} and 
\cite{dr03}.\\  

The formulation of the convolution in terms of sheaf 
cohomology has the advantage that this construction 
works similarly 
also in different categories, e.g., in 
the category of \'etale local systems (see Section~\ref{secetaledef}). 
The latter category 
corresponds to the category of 
 Galois representations of \'etale fundamental groups (compare to 
Prop.~\ref{propequivofetallocal}). 

In this work, we  study  the middle
convolution of (\'etale) local systems in 
view of the inverse Galois problem and related questions. 
The philosophy is that by convoluting elementary objects, 
like the \'etale local systems associated to cyclic or 
dihedral extensions of $\QQ(t),$ one obtains highly 
non-trivial Galois realizations.\\

Let  $H$ be a profinite group. One says that 
$H$ is {\em realized as a Galois group over $\QQ$} if there exists
a surjective homomorphism
$$\kappa:G_\QQ=\Gal(\bar{\QQ}/\QQ) \To H\,.$$ 
 Similarly, the group  $H$ is said to be {\em realized
 as a Galois group over $\QQ(t)$} if there 
exists a surjective homomorphism 
$$\kappa_t:G_{\QQ(t)}=\Gal(\overline{\QQ(t)}/\QQ(t))\To H\,.$$
One says that $H$ is {\em realized regularly
 as Galois group over $\QQ(t)$} if it is realized as Galois group 
over $\QQ(t)$ via a surjection $\kappa_t$ 
such that $\QQ$ is algebraically closed in the fixed field
of the kernel of $\kappa_t$. 

A homomorphism $\kappa$ or $\kappa_t$ as above 
 is called a {\em Galois realization} 
of $H$ over $\QQ,$ or over ${\QQ(t)},$ respectively.
By Hilbert's irreducibility theorem,
a Galois realization of a finite group $H$ over $\QQ(t)$ 
gives Galois realizations of $H$ over $\QQ$ by suitable specializations of
the  transcendental  $t,$ see \cite{MalleMatzat}, \cite{Voelklein}.
A similar argument often holds for profinite groups, see 
\cite{SerreMordellWeil}, 10.6, and \cite{Terasoma85}. 

It is a fundamental task in arithmetic 
to determine the 
Galois  realizations over   $\QQ$ and over 
${\QQ(t)}.$ 
One conjectures that all finite groups can be 
realized over $\QQ$ and over $\QQ(t).$ This conjecture is called  
 the {\em inverse Galois problem for $\QQ,$} and {\em for $\QQ(t),$} 
respectively
(see \cite{Voelklein}, \cite{MalleMatzat}). By Shafarevich's theorem,
any solvable group can be realized as Galois group over 
$\QQ,$ see \cite{NSW}. For 
 non-solvable groups, only partial results
are known, see \cite{MalleMatzat}, \cite{Voelklein}, \cite{dr99}  and \cite{dr00}.
 It already follows from 
class field theory that  not all 
profinite groups can be realized over $\QQ.$
See \cite{Voelklein93}, \cite{Voelklein03}, \cite{dkr}
 and \cite{dw03} for Galois 
realizations of some non-solvable profinite groups.\\

The relevance of the middle convolution for the inverse 
Galois problem was first noticed by S. Reiter and the 
author (\cite{dr99}). 
In the subsequent paper 
\cite{dr00}, many new families of classical groups were 
realized regularly as Galois groups over $\QQ(t)$ using 
the middle convolution with Kummer sheaves (compare to Remarks 
~\ref{remmclambda} and \ref{remlab}).
 Also, it 
turned out that most of the older results on Galois 
realizations  of classical groups (including the famous results 
of Belyi \cite{Belyi79}) can easily be derived 
using  the middle convolution with Kummer sheaves, 
see \cite{dr00}, Rem.~4.11. A similar approach is given 
by  V\"olklein's {\em braid companion functor}, see
\cite{Voelklein01}, \cite{VoelkleinCrelle}, and \cite{Voelklein03}. 
We remark 
that in all known applications to the inverse Galois problem, 
the braid companion functor  can be expressed in terms of the middle 
convolution.

So, up to date, the middle convolution is the 
most powerful tool for the realization of  classical groups as Galois groups 
over $\QQ(t).$
 As already mentioned, the idea is to 
convolute elementary objects (like one-dimensional 
local systems) in order to obtain 
new Galois realizations. \\

Nevertheless, the existing  methods have some limitations: 
The known results usually only encode
the topological information of the Galois realizations under consideration 
and most of the time it is
not possible to obtain more  precise information 
(like  the determination of Frobenius elements
acting on specializations).

Also, in many cases the existing methods
fail to produce simple groups like the 
projective special linear groups $\PSL_n(\FF_q)$  as Galois groups
over $\QQ(t)$  (see \cite{dr99} and 
\cite{MalleMatzat} for  realizations of some special linear groups).
The problem is that if the index of $\PSL_n(\FF_q)$
in $\PGL_n(\FF_q)$ is $>1,$ then it is usually impossible to bound the 
arithmetic part of the underlying Galois representations. \\

Based on methods
of Katz \cite{Katz96} and on previous work of S. Wewers and the author 
\cite{dw03}, we give a geometric 
approach to the middle convolution in order to overcome 
 the above mentioned limitations of the existing methods:
Using this approach, one 
obtains valuable information on the occurring determinants,
leading to new Galois realizations of special linear groups. Moreover,
 computation of Frobenius elements for many of the known Galois 
realizations of classical groups is now possible. 

Thus the geometry of the middle convolution leads not only 
to new results but also 
to a 
much finer
description of the 
existing  results on Galois realizations of classical 
groups.
It turns out that often (e.g., in the 
case of the family of $\Kd$-surfaces considered below)
 this description  of the convolution
encodes the arithmetic information of the resulting 
Galois realizations much more compactly than 
any describing polynomial would do. This shows the benefits 
of the more refined approach 
 also from a computational viewpoint.\\

Let us now discuss the contents of this work in more detail:\\

For a topological space $X$,  let
$\LS_R(X)$ denote   the category of local systems of $R$-modules on
$X$ (see Section~\ref{localsyst}).

Let $\uo$ and $\vo$ denote divisors on $\AA^1.$ We 
define the middle convolution of local systems 
$\V_1\in \LS_R(\AA^1\setminus \uo)$ and 
$\V_2\in \LS_R(\AA^1\setminus \vo)$ in analogy to Formula 
\eqref{eqmidconvo} 
with the difference that we work over an open subset of 
$\AA^2$, where the sheaf $\V_1\circ \V_2$ is a local system, see 
Section~\ref{sectconvdef}. 
The middle convolution 
$\V_1\ast \V_2$ of $\V_1$ and $\V_2$
is a local system on $\AA^1\setminus \uo\ast\vo,$
where $ \uo\ast\vo$ is the sum of the divisors $ \uo$ and $\vo,$ see 
Section~\ref{sectconvdef}.

Using  \cite{dw03}, we  derive a formula for
the rank and the  monodromy (see Section~\ref{localsyst} for the 
concept of monodromy of a local system) 
of the middle convolution
$\V_1\ast \V_2,$
see Prop.~\ref{dimensione} and Prop. 
 \ref{mondeltaij}.  Moreover, some 
irreducibility criteria and results on local monodromies
are given, see Thm.~\ref{thmirrd} and Section~\ref{secmonodrconv}. \\

  If $X$ is a variety over $\QQ,$ let
$\LS_R^\et(X)$ denote the  category of \'etale local systems of
$R$-modules on $X,$ where $R$ is a suitable coefficient ring,
 see Section~\ref{secetaledef}. (In view of the inverse 
Galois problem for the field $\QQ,$
 we assume that our ground field is $\QQ.$ Many 
of the statements in the discussion below are valid for more general 
fields.)
It is well known that any \'etale local system $\LS_R^\et(X)$
corresponds to a Galois representation of the \'etale 
fundamental group $\pi_1^\et(X,\bar{x}),$ 
 called the {\it monodromy
representation}, see
\cite{FreitagKiehl}, Prop. A.1.8.

The middle convolution of \'etale local systems (see Section
\ref{secvondefetale}) is defined in an analogous way as the middle
convolution of local systems, using the higher direct image
approach of \cite{dw03}.  
The middle convolution
$\V_1\ast\V_2$ of  $\V_1\in \LS_R^\et(\AA^1_\QQ\setminus \uo)$ and
$\V_2\in \LS_R^\et(\AA^1_\QQ\setminus \vo)$ is 
an \'etale local system on $S:=\AA^1_\QQ\setminus \uo\ast\vo.$ 
It therefore corresponds to a  continuous
representation
$$ \begin{CD} \pi_1^\et(S,{\bar{s}_0})
@>{\rho_{\V_1\ast\V_2}}>> \GL((\V_1\ast\V_2)_{\bar{s}_0})\end{CD}\,.$$ 
 
The \'etale 
fundamental group $\pi_1^\et(S,{\bar{s}_0})$
is isomorphic to the Galois group of the maximal algebraic extension
of $\QQ(t)$ which is unramified outside the divisor  $\uo\ast\vo$ and outside 
$\infty.$
After a choice of a geometric base point $\bar{s}_0$ which is 
defined over $\QQ,$ the geometric fundamental group 
 can be written canonically
as a semidirect product
$$ \pi_1^\et(S,{\bar{s}_0})=
\pi_1^\geom(S,{\bar{s}_0})\rtimes G_\QQ\,,$$
where $\pi_1^\geom(S,{\bar{s}_0}):=
\pi_1^\et(S\times \Spec(\overline{\QQ}),{\bar{s}_0})$ is the 
{\it geometric fundamental group}. (The geometric 
fundamental group is known to be isomorphic to the profinite 
closure of the topological fundamental group $\pi_1(S(\CC),s_0).$)
Thus  one 
obtains a Galois representation 
$$\begin{CD} G_{\QQ(t)}@>>> 
\pi_1^\et(S,{\bar{s}_0})=\pi_1^\geom(S,{\bar{s}_0})\rtimes G_\QQ
@>>> \GL((\V_1\ast\V_2)_{\bar{s}_0})\,.\end{CD}$$
It is the
determination of this Galois representation that we are interested in. \\

The geometric monodromy of the \'etale local system 
${\V_1\ast\V_2}$ can be determined using the convolution 
$\V_1^{\rm an}\ast \V_2^{\rm an}$ of the analytifications
 $\V_1^{\rm an},\, \V_2^{\rm an}$ which are local systems
on the punctured complex affine line, see 
Section~\ref{secetaledef}. 

Moreover,
Deligne's work on the Weil conjectures (\cite{DeligneWeil2})
yields information on the absolute values of the Frobenius
elements which occur in the arithmetic monodromy,
see Section~\ref{secfirstpropetalecase}. \\

To obtain more information
on Galois representations which are associated to 
convolutions,
we proceed as follows:

Let $\ell$ be a prime and 
let $E_\lambda$ be the completion of a number field $E$ 
at a finite prime $\lambda$ of $E$ with $\chara(\lambda)=\ell.$ 
We consider sheaves 
$$\V=(\cdots ((({\cal F}_1\ast {\cal F}_2)\otimes\GGG_1)\ast{\cal F}_3)\otimes 
\cdots \ast {\cal F}_{n})\otimes \GGG_{n-1}\,\in\, \LS_{E_\lambda}^\et(S),\quad 
S=\AA^1_\QQ\setminus {\bf w}\, ,$$ 
which are obtained from 
iteratively  convoluting and tensoring 
$E_\lambda$-valued  \'etale local systems with finite monodromy. 
(All known applications to the inverse Galois problem use these sheaves.) For 
such \'etale local systems $\V,$ we
find a geometric interpretation, see
Section~\ref{secmotiv}.
 This interpretation involves 
a smooth map $\Pi:{\frak U}\to S$ of relative dimension $n-1$ 
 and 
a certain projection $\PPP\in \End({\frak U}/S)\otimes E.$ 
Then, $\V$ is isomorphic to the weight-$(n-1)$-part of the 
subsystem $\PPP(R^{n-1}\Pi_*(E_\lambda))$ of $R^{n-1}\Pi_*(E_\lambda)$ 
which is cut out by $\PPP,$ see Thm.~\ref{thmpureweight}.

 For actual 
computations and deeper results on Frobenius elements, 
the  results of Bierstone and  Milman (see \cite{BierstoneMillman}) and 
of Encinas, Nobile and Villamayor (\cite{EV1}, \cite{ENV}) 
on equivariant desingularization in characteristic zero 
play an important role: They lead to a $\PPP$-equivariant 
normal crossings compactification
$\tilde{\CX}_{\bar{s}_0}$ of the fibre 
${\frak U}_{\bar{s}_0}$ of ${\frak U}$ over the geometric base point 
$\bar{s}_0$ of $S$ (which is assumed to be defined over 
$\QQ$).  It follows that there are isomorphisms of 
$G_\QQ$-modules
\begin{eqnarray}\nonumber 
\V_{\bar{s}_0}&  \simeq & \left(\mbox{\rm the weight-$({n-1})$-part of}\quad  \PPP(R^{n-1}\Pi_*(E_\lambda))\right)_{\bar{s}_0}\\
&\simeq &\mbox{\rm the weight-$({n-1})$-part of}\quad  \PPP(H^{n-1}({\frak U}_{\bar{s}_0},E_\lambda)) \nonumber \\
&\simeq & \im\left( \PPP(H^{n-1}(\tilde{\CX}_{\bar{s}_0},E_\lambda)) \To 
\PPP(H^{n-1}({\frak U}_{\bar{s}_0},E_\lambda))\right)\, ,\nonumber \end{eqnarray} 
see Prop.~\ref{corisogalmod}. (Here, $G_\QQ$ is viewed as a subgroup of
$\pi_1^\et(S,\bar{s}_0)=\pi_1^\geom(S,\bar{s}_0)\rtimes G_\QQ.$)\\

These results,
 a deep theorem of Henniart on the algebraicity
of one-dimensional compatible systems (see \cite{Schappacher}, Prop. 1.4,
 \cite{Henniart81}, and 
\cite{Khare03}), and the theory of Hecke characters
 lead to the following result, 
 see Thm.~\ref{thmdet}:\\

\noindent{\bf Theorem I.} {\it 
Let 
$\FFF_1,\,\FFF_2,\,\FFF_3,\,\GGG
$  be irreducible and non-trivial $E_\lambda$-valued  \'etale local systems on punctured affine lines over $\QQ$ having finite monodromy. Let 
$$ \V= ((\FFF_1\ast \FFF_2)\otimes\GGG)\ast \FFF_3\,,$$
and 
let $\rho_\V:\pi_1^{\rm geo}(S)\rtimes G_\QQ\to \GL(V)$ be the 
Galois representation
 associated to $\V.$ Assume that 
the associated monodromy 
tuple 
$T_{\V}$ (see Section~\ref{secetaledef}) 
generates an infinite and absolutely irreducible subgroup 
of $\GL(V){\it .}$
 Then the determinant of $\rho_\V$ is of the form 
\begin{equation}\nonumber \label{eqzweii} \det(\rho_\V)=
\det(\rho_\V)|_{\pi_1^\geom(S)}\otimes \chi_\ell^m\otimes
\epsilon\,.\end{equation} Here, 
$\chi_\ell:G_\QQ\to \QQ_\ell$  denotes the $\ell$-adic cyclotomic character, 
$m$ is
an integer,  and $\epsilon:G_\QQ\to E^\times$ is a finite
character. }\\

The proof of Thm.~I makes essential use of the fact that $\QQ$ is a totally
real field (which implies that any Hecke character of $\QQ$ is a 
product of the norm character and a finite character).\\

Thm.~I is crucial for the next result
 (see  Thm.~\ref{thmrealierung1} and Corollary
\ref{corsldreia}):\\

\noindent {\bf Theorem II.} {\it Let $\FF_q$ be the finite field 
of order $q=\ell^k,$ where $k\in \NN .$ 
Then the special linear group $\SL_{2n+1}(\FF_q)$  
occurs regularly as Galois
group over $\QQ(t)$ if $$q\equiv 5\mod 8\quad \mbox{and}\quad 
n> 6+2\varphi((q-1)/4)$$ ($\varphi$ denoting
Euler's $\varphi$-function){\it .}}\\

Thm.~II implies that, under the conditions  of the theorem,
the simple  group 
$\PSL_{2n+1}(\FF_q)$ occurs regularly as Galois
group over $\QQ(t).$ 
The latter result is the first result on regular Galois
realizations of the groups ${\rm PSL}_n(\FF_q)$ over $\QQ(t),$ 
 where 
$$(n,q-1)=[{\rm PGL}_n(\FF_q):{\rm PSL}_n(\FF_q)]> 2\,.$$

The idea of the proof of Thm.~II is the following:  By our assumptions,
the finite field $\FF_q$ is generated 
over its prime field $\FF_\ell$ by an element of odd order $m.$
Let $$E:=\QQ(\zeta_{m}+\zeta_{m}^{-1},i)\, ,$$
where 
$m=(q-1)/4,$  where $\zeta_m$ denotes
 a primitive $m$-th root of unity, and where $i$ is a primitive fourth root of unity. Let 
$\lambda$ be
a prime of $E$ lying over $\ell.$
One considers 
$E_\lambda$-valued  \'etale local systems $\FFF_1,\,\FFF_2,\,\FFF_3,\,\GGG$ 
associated to Galois representations with values in the dihedral group 
of order $2m$ and to Galois representations with values in  
cyclic groups of order  $2$ and $4.$ Then  
one forms their convolution
$$ \V= ((\FFF_1\ast \FFF_2)\otimes \GGG)\ast \FFF_3\,\in\, \LS^\et_{E_\lambda}(S)\,.$$
Let $\rho_\V:\pi_1^\et(S)\to \GL(V)$ be the 
Galois representation
 associated to $\V$ and let
$O_\lambda$ be the valuation 
ring of $E_\lambda.$ 
Using analytification and reduction modulo $\lambda,$ one can show that 
the image of the geometric fundamental group 
 under 
$\rho_\V$  is, up to scaling, isomorphic 
to $\SL_{2n+1}(O_\lambda).$ Here, 
$n$ depends on $m,$  enforcing the condition 
$n> 6+2\varphi(m).$
  Since $m$ is odd,
 the only roots of unity which are 
contained in $E$ are fourth roots of unity. 
It  then follows from Thm.~I that the occurring determinants which arise from 
$G_\QQ\leq \pi_1^\et(S)$ are, 
up to a twist with the cyclotomic 
character, contained in the 
group of fourth roots  of unity. Thus by a 
twist with a suitable finite character of order four, one can assume that the 
 image of the whole \'etale fundamental 
group $\pi_1^\et(S)$ 
is equal to $\SL_{2n+1}(O_\lambda).$
The result then follows from 
reduction modulo $\lambda$ and from the above interpretation 
of $\pi_1^\et(S)$ as a factor of $G_{\QQ(t)}.$

The proof of Thm.~II implies that (under 
the above restrictions) the profinite groups 
$\SL_{2n+1}(O_\lambda)$ occur regularly as Galois groups over $\QQ(t),$
see Thm.~\ref{thmrealierung1}.
This leads to towers of Galois covers with interesting arithmetic
properties. These towers 
 seem to be intimately connected to recent work of 
M. Fried on modular towers, see \cite{Fried03}.

Let us also remark
that the conditions on $\FF_q$ can be weakened in many ways.
We have chosen this formulation because of the simplicity of the 
statement. \\

In Section~\ref{secspecil} we consider specializations of Galois groups:
Using the middle convolution,
 we construct regular Galois extensions 
$L_\ell/\QQ(t)$ with 
$$\Gal(L_\ell/\QQ(t))\simeq \PGL_2(\FF_\ell),\quad {\rm where} \quad 
\ell\not\equiv \pm 1\mod 8\quad \mbox{\rm is an odd prime,} $$ and we 
determine the behaviour of the Galois groups under the specialization 
of $L_\ell$  via  $t\mapsto 1.$   (The existence of such Galois 
extensions is a classical result, see \cite{Matzat87}, Folgerung 2, p. 181.
What is new is the determination of the specializations.)

 Using the geometric description of the 
convolution, one finds a family 
$$\pi:\tilde{\frak X}\To \AA^1_\QQ\setminus \{\pm 2,0\}$$
of $\Kd$-surfaces of Picard number
$19$ such that the extensions $L_\ell$ occur in the cohomology of 
this relative surface. The action of $G_\QQ$ on the specialization 
of $L_\ell$ under  
$$t\,\mapsto\, s_0\,\in \, \QQ\setminus \{\pm 2,0\}\, ,$$  
amounts to a twist of a subrepresentation  of 
 $$G_\QQ \To \Aut\left(H^2(\tilde{\frak X}_{\bar{s}_0},\FF_\ell)\right)\,,$$ where $\tilde{\frak X}_{\bar{s}_0}$
denotes the fibre over $\bar{s}_0$ (and 
$\bar{s}_0$ denotes a geometric point 
induced by $s_0$).
 The  
Lefschetz fixed point formula for \'etale cohomology
of varieties in characteristic $p>0$ can be used to compute 
the action of the Frobenius elements. 
The crystalline 
comparison isomorphism  relates the Hodge numbers to 
the image of the tame inertia at $\ell,$ see Section 
\ref{tameinertia}. Using these ingredients, we prove that the 
specializations of 
the Galois groups remain constant for infinitely many choices of $\ell,$
see Thm.~\ref{thmspecki}:\\

\noindent {\bf Theorem III.} {\it
 Let $\ell\not\equiv \pm 1\mod 8$ be an odd prime.
Then there exists a regular Galois extension  
$L_\ell/\QQ(t)$  with Galois 
group isomorphic to $\PGL_2(\FF_\ell){\it .}$  Let $L'_\ell/\QQ$ denote the 
specialization of $L_\ell/\QQ(t)$ under $t\mapsto 1{\it .}$ 
Then the following holds:
\begin{enumerate}
\item If, in addition, $\ell\geq 11$ and
$\left(\frac{-3}{\ell}\right)=-1\,,$ then the Galois 
group $\Gal(L_\ell/\QQ(t))$ 
is preserved under specialization $t\mapsto 1,$ i.e.,
$$\Gal(L'_\ell/\QQ)\,\simeq\, \PGL_2(\FF_\ell)$$
(where $\left( \frac{\cdot}{\cdot }\right)$ denotes the Legendre symbol).
\item
For almost all $\ell\not\equiv \pm 1,2\mod 8,$ the Galois 
group $\Gal(L_\ell/\QQ(t))$ is preserved under specialization $t\mapsto 1{\it .}$ 
\end{enumerate} }

Let us close the introduction with a discussion on some 
potential applications of the methods considered here: \\

The varieties $\Pi:{\frak U}\to S,$ 
and their 
compactifications, which are constructed in  the convolution process  
appear in many areas of mathematics. For example, they appear in the 
study of Painlev\'e equations or, more generally,
in the study of 
isomonodromic deformations of differential equations, 
see \cite{dr04}. Another example is the 
 Legendre family of elliptic curves which also
occurs in this way (by convoluting a double cover of $\AA^1\setminus \{0,1\}$
with an involutory Kummer sheaf). 

Thus the geometric theory 
of the  convolution 
serves as a rich source of meaningful varieties
for which many of the conjectures in arithmetic geometry can be asked
(and solved via the convolution process?):
Modularity of 
the Galois representations on the fibre (see \cite{Langlands}), 
existence and construction 
of semi-stable models (see \cite{deJong}), to name a few.
 Of course,
the solution of 
any such question leads to a better understanding of the resulting 
Galois 
realizations.

Also, in many cases there is a natural 
map of these varieties into Shimura varieties of PEL-type
\cite{dw05}. Thus 
 the arithmetic of the Galois representations which are 
associated to convolutions translates into 
statements in the theory of Shimura varieties and vice
versa. For example, the Andr\'e-Oort conjecture (see \cite{Andre89}, Chap. 
X.4, \cite{EdixhovenYafaev} and \cite{Oort})
can be translated into 
a statement on the specializations of the occurring Galois representations. 

Another potential application 
of our methods is this:
While a full solution of the inverse Galois problem over $\QQ(t)$ still 
seems to be out of reach, the situation over the field 
$\QQ^{\rm ab}(t)$ (where  $\QQ^{\rm ab}$ denotes the 
maximal abelian Galois  extension of $\QQ$) is more hopeful: 
The author conjectures that 
any finite group can be realized in the cohomology of some 
smooth affine 
variety  over $S=\AA^1_{\QQ^{\rm ab}}\setminus \uo,$ where 
$\uo$ is some divisor.  
In this context, the 
affine varieties $\Pi:{\frak U}\to S$
which appear 
in the convolution process 
should be the major tool for realizing unipotent groups over $\QQ^{\rm ab}(t)$
and 
for solving embedding
problems of a classical group with a unipotent group.\\

{\bf Acknowledgments:} I  heartily thank  B.H. Matzat for 
his support during the last years, for his advice,
 and  for many valuable discussions 
on the subject of this work. 
I thank M. Berkenbosch, 
J. Hartmann, U. K\"uhn,  S. Reiter, A. R\"oscheisen and S. Wewers 
for valuable comments and discussions and for pointing
out several errors and misprints in earlier versions of this work. 
Part of this work was written 
during my stay at the School of Mathematics of the Tel Aviv University
in the spring of 2004. 
I thank M. Jarden and D. Haran for their hospitality during this stay
and for interesting discussions. 
Also, I acknowledge the financial support provided through 
the European Community's Human Potential Program under the contract
HPRN-CT-2000-00114, GTEM.\\

\section{Preliminaries}\label{secgrouptheoryprel}

It is the aim of this section to set up the notation
and to state some (mostly well known) 
results which are used in later sections.

\subsection{Group theoretical definitions and tensor products 
of representations.}\label{secdefacon}

In this subsection we collect useful results on tensor
representations. The content is well
known, compare to
\cite{KleidmanLiebeck90}, Chap. 4.4.\\

Let $G$ be a group. 
If  $(g_1,\ldots,g_r)\in G^r$ is a tuple of elements 
of $G,$ then we set
\begin{equation} \label{r+1} g_{r+1}:=(g_1\cdots g_r)^{-1}.
\end{equation}

 Let $R$ be a commutative ring
with a unit. If $V$ is a free 
$R$-module, then $\GL(V)$ denotes the $R$-linear 
isomorphisms of $V.$ As usual, the group of 
$R$-linear isomorphisms of $R^n$ is denoted by $\GL_n(R).$
Linear automorphisms act from the right, i.e., if $g
\in \GL(V)$ and $v\in V,$ then $vg$ denotes the image of $v$ under
$g.$ \\

A Jordan block of eigenvalue $\alpha\,\in R$
and of length $n$ is denoted by $J(\alpha,n).$ We write
$$J(\alpha_1,n_1)\oplus \cdots \oplus J(\alpha_m,n_m)$$ for a block
matrix in $\GL_{n_1+\cdots +n_m}(R)$ which is in Jordan normal
form and whose Jordan blocks are $J(\alpha_1,n_1),\ldots ,
J(\alpha_m,n_m).$  \\

Let  $V_1,\ldots,V_t$ be free $R$-modules of rank
 $n_1,\ldots,n_t,$ respectively.  For
$$g_i \in \GL(V_i),\,i=1,\ldots,t,$$
 define an element $$g_1\otimes \cdots \otimes
g_t \in \GL(V_1\otimes \cdots \otimes V_t) $$ by setting
$$ (v_1\otimes \cdots \otimes v_t)(g_1\otimes \cdots \otimes g_t):=
v_1g_1\otimes \cdots \otimes v_tg_t\,. $$ This tensor product of
matrices is also called the {\em Kronecker product}. Let $\rho_1:
H_1\to \GL(V_1)$ and $\rho_2: H_2\to \GL(V_2)$ be 
representations. Then the tensor product defines a representation
$$ \rho_1\otimes \rho_2: H_1\times H_2 \To \GL(V_1\otimes V_2),\;
(h_1,h_2)\Mapsto \rho_1(h_1)\otimes \rho_2(h_2)\,.$$

Let us assume for the rest of this subsection that $R=K$ is a field.
If $f_1,\ldots,f_t$ are bilinear forms on $V_1,\ldots,V_t,$
respectively, then the formula 
\begin{equation}\label{form}
f(v_1\otimes \cdots \otimes v_t,w_1\otimes \cdots \otimes w_t):=
\prod_{i=1}^t f_i(v_i,w_i)\end{equation}
defines by bilinear extension a bilinear form $f=f_1\otimes
\cdots \otimes f_t$ on $V.$ 

Let the characteristic of $K$ be $\not=2.$ If $W$ is a $K$-vector
space and $f$ is a bilinear form on $W,$ then we set $\sym(f):=-1$
if $f$ is symplectic, $\sym(f):=1$ if $f$ is orthogonal, and
$\sym(f):=0$ otherwise. If the characteristic of $K$ is $2,$ then we
set $\sym(f):=1$ if $f$ is orthogonal and $\sym(f):=0$ otherwise. Then
\begin{equation}\label{sym}
\sym(f_1\otimes \cdots \otimes f_t)=\prod_{i=1}^t
\sym(f_i)\,.\end{equation}

It is often important to compute the Jordan normal form of the 
tensor product $A\otimes B$ of two
matrices $A \in \GL_n(K)$ and $B\in \GL_m(K).$ This can be done
using the following lemma (see \cite{MartsinkVlassov} for a
proof):

\begin{lem} \label{lemkroneckerjordan} Let $K$ be an algebraically
closed field of characteristic zero, let $\alpha,\beta \in K,$
and  let 
$n_1\leq n_2.$ Let $J(\alpha,n_1)\in \GL_{n_1}(K)$ and
$J(\beta,n_2)\in \GL_{n_2}(K)$ be two Jordan blocks. Then the
Jordan normal form of $J(\alpha,n_1) \otimes J(\beta,n_2)$ is
given by
$$ \bigoplus_{i=0}^{n_1-1} J(\alpha\beta,n_1+n_2-1-2i)\,.$$
\end{lem}

 \subsection{Braid groups and affine fibrations.}\label{braidgroups}
We
will write $\AA^1,\, \PP^1 , \ldots$ instead of $\AA^1(\CC),\,
\PP^1(\CC), \ldots$ and  view these objects equipped with their
associated topological and complex analytic structures.\\

 Let $X$ be a connected topological manifold and let 
$\pi_1(X,x)$ denote the fundamental group of $X$ with base point $x.$
 The multiplication in $\pi_1(X,x)$
is induced by the path product which is defined 
using the following convention: Let  $\gamma, \gamma'$ be two 
closed paths at $x\in X.$ Then their product $\gamma \gamma'$ is given
by first walking along $\gamma$ and then walking along $\gamma'$. \\

 Let $r\in \NN_{>0}$ and let $U_0= \AA^1\setminus \uo,$
where $\uo:= \{ u_1,\ldots,u_r\}$ is a finite subset 
of $ \AA^1.$  We will identify
$U_0$ with $\PP^1 \setminus (\uo \cup \{\infty\})$ in the
obvious way. Using a suitable homeomorphism $\kappa:
\PP^1\to\PP^1,$ (called a {\em marking} in \cite{dw03},
Section 1.2) and the conventions of loc.~cit., Section 2.3, one
obtains generators $\al_1,\ldots,\al_{r+1}$ of $\pi_1(U_0,u_0)$
which satisfy the product relation $\al_1\cdots\al_{r+1}=1.$ \\

Consider the configuration spaces 
$$ \OO_{r}:=\{ \uo \subseteq \AA^1\, \mid\, |\uo|=r\}$$
and
$$\OO_{r,1}:=\{ (\uo,x) \in \OO_r\times \AA^1 \, \mid\, x \notin \uo\}\,.$$
The sets $\OO_r$ and $\OO_{r,1}$ are connected topological manifolds in a
natural way. In fact, there is a standard way to identify  
the space $\OO_r$ with $\AA^r\setminus \Delta_r,$ where 
$\Delta_r$ denotes the discriminant locus 
 (see \cite{Voelklein}). We set $ \A_r:=\pi_1(\OO_r,\uo)$ and
$\A_{r,1}:=\pi_1(\OO_{r,1}, (\uo,u_0)).$ The marking on $U_0$
also defines standard generators $\beta_1,\ldots,\beta_{r-1}$ of
$\A_r$ which satisfy the usual relations of the standard
generators of the Artin braid groups (\cite{dw03}). One obtains a split exact
sequence
\begin{equation}\label{eqbraids1} 1\To\pi_1(U_0,u_0)\To\A_{r,1}
\To\A_r \To 1\end{equation} such that the following equations
hold (with respect to the splitting):
\begin{equation} \label{locals5eq2}
  \beta_i^{-1}\,\alpha_j\beta_i \;=\quad
  \begin{cases}
    \quad \alpha_i\alpha_{i+1}\alpha_i^{-1}, &
                              \quad\text{\rm for $j=i$,} \\
    \quad \alpha_i, & \quad \text{\rm for $j=i+1$,} \\
    \quad \alpha_j, & \quad \text{\rm otherwise.}
  \end{cases}
\end{equation}

As usual, one sees that $\A_r$ acts (product preserving)  on $G^r,$ where $G$ is any group, as follows:
\begin{multline}\label{braidaction}
(g_1,\ldots,g_r)^{\beta_i}=\\
(g_1,\ldots,g_{i-1},g_{i+1},g_{i+1}^{-1}g_ig_{i+1},g_{i+2},\ldots,g_r),
\; \forall\; (g_1,\ldots,g_r) \in G^r. \end{multline}

Let $$\OO^r:=\{(v_1,\ldots,v_r)\in \AA^r\mid v_i\not= v_j \quad
{\rm for} \quad i\not= j\}\,.$$ Let
$\A^r:=\pi_1(\OO^r,(u_1,\ldots,u_r)).$ The map $$ \OO^r \To
\OO_r,\; (v_1,\ldots,v_r)  \Mapsto \{v_1,\ldots,v_r\}$$ is a
unramified covering map. Thus (via the lifting of paths) $\A^r$
can be seen as a subgroup of $\A_r.$ It is well known that 
$\A^r$ is then generated by the following braids:
\begin{equation}\label{purebraidgens}
\beta_{i,j}:=(\beta_i^2)^{\beta_{i+1}^{-1}\cdots
\beta_{j-1}^{-1}}= (\beta_{j-1}^2)^{\beta_{j-2}\cdots \beta_i},
\end{equation}
where $1\leq i <j \leq r.$\\

 Let $S$ be a  connected complex manifold, let $X:=\PP^1_S=\PP^1\times S,$
and let $\D\subseteq X$ be a smooth relative divisor of degree $r+1$ over $S$ which
contains the section $\{\infty\} \times S.$ Let $U:=X \setminus
\D,$ let $j: U \to X$ be 
the natural inclusion, let $\bar{\pi}: X \to S$ be the
projection onto S, and let $\pi: U \to S$ be the restriction of
$\bar{\pi}$ to $U.$ Let further $s_0 \in S$ and suppose that
$U_0=\pi^{-1}(s_0).$ There is a continuous map
$$  S \To \OO_r,\,\, s \Mapsto \pi'(\bar{\pi}^{-1}(s)\cap \D)\setminus
\infty\, ,$$ where $\pi'$ denotes the projection of $X$ onto
$\PP^1.$ This map induces a homomorphism of fundamental groups
$\phi: \pi_1(S,s_0) \to \A_r.$ Similarly, the map
$$  U \To \OO_{r,1},\,\, (u,s) \Mapsto (\pi'(\bar{\pi}^{-1}(s)\cap \D)\setminus
\infty, \pi'(u))\, ,$$ gives rise to a homomorphism
$\tilde{\phi}:\pi_1(U,(u_0,s_0)) \to \OO_{r,1}.$

In \cite{dw03} it is shown how to obtain a commutative diagram whose
rows are split exact  sequences:
\begin{equation}\label{splitdiag}\begin{CD} 1\To\pi_1(U_0,u_0)
@>>> \pi_1(U,(u_0,s_0))
    @>>>  \pi_1(S,s_0) \To 1 \\
   @VVV   @V\tilde{\phi}VV @V{\phi}VV       \\
  1\To \pi_1(U_0,u_0)@>>>\A_{r,1} @>>>\quad\quad\quad\A_r\To 1\, .    \\
\end{CD}\end{equation}
Let $\iota_1:\pi_1(S,s_0)\to\pi_1(U,(u_0,s_0)) $ denote the
splitting of the upper row and let $\iota_2:\A_r \to\A_{r,1}$ be the
splitting of the lower row. Then one may assume that 
\begin{equation}\label{rowsplitting}
\tilde{\phi}\circ \iota_1=\iota_2\circ\phi \end{equation}
 (see loc.~cit., Section 2.3 and Rem. 2.6).

\subsection{Conventions on sheaves.}\label{convsheav} For the
definition and basic properties of sheaves we refer to \cite{Hartshorne}.
Let $X$ be a topological space. The category of sheaves of abelian
groups on $X$ is denoted by $\Sh(X).$ If $R$ is a ring, then the
category of sheaves of $R$-modules on $X$ is denoted by
$\Sh_R(X).$
The stalk of $\V \in \Sh(X)$ at $x \in X$ is denoted by $\V_x.$ 
 Let $f:X\to Y$ be a continuous map of
topological spaces and let $\V\in \Sh(X),$ resp. $\W\in \Sh(Y).$ Then
$f_*\V \in \Sh(Y)$ denotes the sheaf
$$ U \Mapsto \V(f^{-1}(U)) \quad (U \subseteq Y \; {\rm open})$$
and  $f^*\W \in \Sh(X)$ denotes the sheaf associated to the
presheaf
$$U \Mapsto \lim_{f(U)\subseteq V }\W(V) \quad (U \subseteq X\; {\rm and}\;
V \subseteq Y
 \; {\rm open})\,.$$
Let $f$ be as above and $\V\in \Sh(X).$ The $i$-th higher direct
image of the functor $f_*$ is  denoted by $R^if_*,$
see e.g. \cite{Hartshorne}. Thus $R^if_*(\V)$ 
 coincides with the sheaf associated
to the presheaf
\begin{equation}\label{eqhigh} U \mapsto H^i(f^{-1}(U), \V|_U) \quad (U \subseteq Y \; {\rm open}),\end{equation}
see loc.~cit. Prop. 8.1.

 The constant sheaf which is associated
to an abelian group $A$ is again denoted by ${A}.$
The tensor product of sheaves of $R$-modules is defined in the
usual way. 

\subsection{Local systems and representations of fundamental
groups.}\label{localsyst}

Let $R$ be a commutative ring with a unit and let 
$X$ be a connected topological manifold. A {\em local system of $R$-modules}
is a sheaf $\V \in \Sh_R(X)$ for which there exists an $n \in \NN$
such that $\V$ is locally isomorphic to ${R}^n.$ The
number $n$ is called the {\em rank} of $\V$ and is denoted by
$\rk(\V).$ Let $\LS_R(X)$ denote the category of local systems of
$R$-modules on $X.$ For any path $\gamma: [0,1] \to X$ on has a
sequence of isomorphisms
$$ \V_{\gamma(0)} \To (\gamma^*\V)_{0} \To \gamma^*\V([0,1])
\To (\gamma^*\V)_{1} \To \V_{\gamma(1)}\,.$$ The composition of these
isomorphisms is an isomorphism of $ \V_{\gamma(0)}$ to
$\V_{\gamma(1)}$ which is denoted by $\rho_\V(\gamma).$ It is easy
to see that $\rho_\V(\gamma)$  only depends on the homotopy class
of $\gamma.$ Thus any local system $\V \in \LS_R(X)$ gives rise
to its {\em monodromy representation}
$$ \rho_\V:\pi_1(X,x) \To \GL(V), \, \gamma \Mapsto \rho_\V(\gamma)$$
(here we have used the convention $V=\V_{x}$). We  always let
$\pi_1(X,x)$ act from the right on $V.$

Let $\Rep_R(\pi_1(X,x))$ denote the category of representations
$\pi_1(X,x) \to \GL(V),$ where $V\simeq R^n$ for some $n\in
\NN.$ The following proposition is well known, see \cite{Deligne70}:

\begin{prop}\label{p11} There is an equivalence of categories
$$ \LS_R(X) \cong \Rep_R(\pi_1(X,x))$$ under which 
$ \V$ corresponds to $\rho_\V\,.$
\end{prop}
The following notation will be useful later:

\begin{notation} \label{notfortupeloflocsys} Let 
$$U_0:=\AA^1\setminus \uo,\, \uo=\{u_1,\ldots,u_r\}
\in \OO_r\,,$$ and fix generators $\al_1,\ldots,\al_{r+1}$ of
$\pi_1(U_0,u_0)$ which satisfy 
the product relation
$\al_1\cdots \al_{r+1}=1$ 
as in Section~\ref{braidgroups}. 
If $\V$ is a given local system on $U_0,$ then
$$ \begin{array}{rcl} \V \in\LS_R(U_0)&\longleftrightarrow &
\rho_\V \in \Rep_R(\pi_1(U_0,u_0))\\
&\longleftrightarrow&
T_\V:=(T_1:=\rho_\V(\al_1),\ldots,T_{r+1}:=\rho_\V(\al_{r+1}))\\
&&\quad \in
\GL(V)^{r+1},\,\, T_1\cdots T_{r+1}=1\end{array}$$ 
(where ``$\longleftrightarrow$'' stands for ``corresponds to'').
We call 
$T_\V$ the {\em associated tuple} of $\V.$\end{notation}

Let $\V \in \LS_R(X)$ be equipped with a {\em bilinear pairing}
$\kappa :\V \otimes \V \to {R},$ corresponding 
to an isomorphism $\V\hookrightarrow \V^*,$
where $\V^*$ denotes the dual local system. It follows from the
definitions that any such pairing induces  a bilinear form
$f_\kappa: V\otimes_R V \to R$ which is compatible with the action of
$\pi_1(X,x).$ We set $\sym(\kappa)=\sym(f_{\kappa}),$
where $\sym(f_{\kappa})$ is as in Section~\ref{secgrouptheoryprel}.\\

\subsection{Field theoretic notation.}\label{secgeneralnotation}

Let $k$ be a field. An algebraic closure of $k$ is denoted by
$\bar{k}$ and a separable closure is denoted by $k^\sep.$ We set
$G_k:=\Gal(k^\sep/k).$ The characteristic of $k$ is denoted by
$\chara (k).$ The set of non-archimedian primes of $k$ is denoted
by $\PP^f(k).$ The characteristic of $\nu \in \PP^f(k)$ is 
denoted by $\chara(\nu).$\\

For $n\in \NN,$ we set
$ \zeta_n:=e^{\frac{2\pi i}{n}},$ where $i$ denotes the primitive 
fourth 
root of unity which induces a positive orientation on $\CC.$\\

A {\it Galois representation} is a continuous
homomorphism $$\rho \;:\; G_k \To \GL(V)\, ,$$ where $k$ 
is a field and $V$ is a free module of finite rank over 
some topological 
 ring $R.$
In this section and in the sections which are concerned 
with \'etale sheaves, the topological ring $R$  will be a 
complete subfield of
$\bar{\QQ}_\ell$ equipped with the induced $\ell$-adic 
topology, the valuation ring 
of such a field, 
or a finite field equipped with the 
discrete topology.\\

The cyclotomic character
is an important example of a  Galois representation: 
 If $\ell$ is a prime different from $\chara (k),$ let $\mu_{\ell^m}\in
k^\sep$ denote the group of $\ell^m$-th roots of unity and let 
$$
\mu_\ell^\infty=\lim_{\leftarrow}\mu_{\ell^m}\simeq \ZZ_\ell\,.$$
 The {\em $\ell$-adic cyclotomic character}
 $\chi_\ell:G_k\to
\ZZ_\ell^\times $ is defined as follows:
$$ g(z)=z^{\chi_\ell(g)}\quad {\rm if}\quad g \in G_k\quad {\rm and}\quad 
 z\in \mu_{\ell^m} \,.$$

Let $\rho:G_k\to \GL(V)$ be a Galois representation,
where $V$ is a vector space of finite dimension over 
 a complete subfield $K$ of
$\bar{\QQ}_\ell.$  
Let $O$ denote the valuation ring of $K$ and let $\m$ be its maximal ideal.  
By \cite{Se}, Remark  on page I.1,
there exists a $G_k$-invariant $O$-lattice $W$ in $V.$
Let $\overline{V}:=W/\m W$ and let 
$ \hat{\rho}: G_k \to \GL(\overline{V})$
be the  representation which arises as the composition 
of $\rho$ with the residual homomorphism $\GL(V)\to \GL(\overline{V}).$
Then the {\it residual representation}
 $$\overline{\rho}:=\hat{\rho}^{\rm s}:G_k \to \GL(\overline{V})$$
does not depend on the chosen lattice (see Lemma 2 in \cite{Wor}),
where the superscript ${}^{\rm s}$ stands for  {\em semisimplification}
 in the 
following sense:

Let $\rho: H \to \GL(W)$ be a  representation, where
$W$ is a finite dimensional vector space over some 
field. Then $W$ has a
composition series
$$ W=W_0 \supset W_1 \supset \ldots \supset W_q =0$$
of $\rho$-invariant submodules such that
$W_i/W_{i+1}, \,
i=0,\ldots,q-1,$ is irreducible.  The {\em semisimplification of
$\rho$} is then the representation on the $H$-module
$$W^{\rm s}:=\bigoplus_{i=0}^{q-1}W_i/W_{i+1}\,.$$

Let $k$ be a number field. 
For $\nu \in \PP^f(k),\; k_\nu$ denotes its completion, $O_\nu$
denotes the valuation ring and $\FF_\nu=O_\nu/\m_\nu$ is
the residue field.
For any $\nu \in \PP^f(k)$ there is the inclusion
$G_{k_\nu} \hookrightarrow G_k$ by viewing $G_{ k_\nu } $ as the decomposition
subgroup of an extension $\overline{\nu}$ of $\nu$ to $\overline{k_\nu}.$
Let $I_\nu:= \Gal(\overline{k_\nu}/k_\nu^{\nr})$ be the {\it inertia subgroup}
at $\nu,$ where $k_\nu^{nr}$
denotes the maximal unramified algebraic 
extension of $k_\nu,$ and let 
$R_\nu:=\Gal(\overline{k_\nu}/k_\nu^{t}),$ where $k_\nu^{t}$ is the maximal
tamely ramified algebraic extension of $k_\nu.$ The {\em tame inertia subgroup}
 at $\nu$ 
is defined to be 
$I_\nu^t:=I_\nu/R_\nu.$ 
A Galois representation $\rho: G_k\to \GL(V)$
 is called {\em unramified at} $\nu \in \PP^f(k)$ if
$\rho(I_\nu)=1.$

There is an exact sequence
$$ 1\To I_\nu\To G_{k_\nu}\To G_{\FF_\nu} \To 1\, ,$$
where $G_{\FF_\nu}$ is topologically generated by the geometric
Frobenius element $F_\nu$ (the inverse to the arithmetic
Frobenius $\alpha \mapsto \alpha ^{|\FF_\nu|}$). (See 
\cite{Katz932} for a discussion 
on the action of  the geometric Frobenius on the \'etale cohomology.)
A {\em Frobenius
element} is any element $\Frob_\nu \in G_k$ conjugate to an
element of $G_{k_\nu}$ which is mapped to $F_\nu$ under the
above homomorphism $G_{k_\nu}\to G_{\FF_\nu}.$
If a Galois representation $\rho:G_k\to \GL(V)$ is unramified at 
$\nu \in \PP^f(k),$  then the element $\rho(\Frob_\nu)$
is determined by $\nu$ up to conjugacy in the image of $\rho.$

\begin{defn}\label{defeinss}
{\rm  Let $k$ be a number field and let 
$\rho: G_k \To \GL(V)$ be a Galois representation,
where $V$ is a vector space over a complete subfield 
$E_\lambda$ of 
$\bar{\QQ}_\ell.$  
Then $\rho$ (resp. $V$) is called {\em pure of weight $w\in \QQ,$}
if for any embedding $\iota : \bar{\QQ}_\ell \to \CC$ and any $\nu
\in \PP^f(k)$ for which $\rho$ is unramified, the following holds:
If $\alpha$ is an eigenvalue of $\rho(\Frob_\nu),$ then
$$ |\iota(\alpha)|_\CC =|\FF_\nu |^{\frac{w}{2}}\, ,$$
where $|\cdot |_\CC$ denotes the complex absolute value.}
\end{defn}

\subsection{Geometric preliminaries.}\label{secgeomprell}

Let $k$ be a field and let $\Var_k$ denote the category of
varieties over $k,$ i.e., the category of reduced schemes of finite type
over $k.$ For any
field  $k'$ which contains $k,$ there is an extension
functor
$$ \Var_k \To \Var_{k'},\; X=X_k \Mapsto X_k \times \Spec(k') =X_{k'}\,.$$

Let $S$ be a variety over 
$k$ and let 
$s:\Spec(k')\to S$ be a $k'$-rational point of $S,$ where 
$k\subseteq k'.$ Then $\bar{s}$ denotes a geometric point which 
arises from $s$ via a composition 
$$\Spec(\overline{k'})\to \Spec(k')\to S\,.$$
Let
$\pi:X\to S$ be a morphism. 
 Then $\pi_s$ denotes the pullback of $\pi$ along $s$ and 
$X_s$ denotes the fibre over $s.$ \\

The following varieties appear throughout the paper:
For $r \in \NN_{>0},$ define 
$\OO_{r,k}:=\AA^r_k\setminus \Delta_r,$
where $\Delta_r$ denotes the discriminant locus. 
Let 
$$\OO_r(k):=\{ \{u_1,\ldots,u_r\}\mid u_i\in \AA^1_k({\bar{k}}),\,\,
i\not= j \Rightarrow u_i\not= u_j,\,\,(x-u_1)\cdots (x-u_r)\in k[x]\}\,.$$
There is a natural identification $\OO_r(k)\simeq \OO_{r,k}(k)$
see \cite{Voelklein}, Lemma 10.17.   For 
$\uo:=\{u_1,\ldots,u_r\} \in \OO_{r}({k}),$ let   
$$\AA^1_k\setminus \uo:=\Spec\left(k[x,\frac{1}{(x-u_1)\cdots (x-u_r)}]\right)\,.$$

\subsection{Etale sheaves and local systems.}\label{secetaledef}

In this section we recall the basic properties of \'etale sheaves.
The standard references are the books of Freitag and Kiehl
\cite{FreitagKiehl} and of Milne \cite{MilneEC}.\\ 

Let   $X\in \Var_k$ be a smooth and geometrically irreducible 
variety over a field $k$ and let $\ell$ be a prime $\not=\chara(k).$
  Our coefficient ring $R$  will be a topological 
ring as in Section~\ref{secgeneralnotation}, i.e.,
$R$ will be a complete subfield of
$\bar{\QQ}_\ell$ with the induced topology, the valuation ring 
of such a field, 
or a finite field of characteristic $\ell$
 equipped with the 
discrete topology.\\

 Let $\Sh_R^\et(X)$ denote the category of
\'etale sheaves of finitely generated 
$R$-modules on $X.$  The constant \'etale sheaf associated to $R^n$ is again denoted by 
${R}^n.$ The stalk of $\V \in \Sh_R^\et(X)$
at a geometric point $\bar{x}$ of $X$ is denoted by $\V_{\bar{x}}$ (again, we
will often set $V=\V_{\bar{x}}$). For \'etale sheaves in 
$\Sh_R^\et(X)$ one
has the notions of { direct image, inverse image, extension by
zero, tensor product, higher direct image} etc. (in short:
Grothendieck's six operations), which are parallel to the case of
sheaves on topological spaces 
and which are denoted by the same symbols. \\

\begin{defn}{\rm  \begin{enumerate}

\item  An {\em \'etale local system} $\V$ of $R$-modules on $X$ is
 a { locally constant sheaf of $R$-modules} on
$X$ (in the \'etale topology) 
whose stalks are
free $R$-modules of finite rank, see \cite{FreitagKiehl}. 
An \'etale local system 
is also
called a {\em lisse $\ell$-adic sheaf} \cite{SGA4}. The category
of \'etale local systems on $X$ is denoted by $\LS_R^\et(X).$
\item An \'etale sheaf $\V\in \Sh_R^\et(X)$ is called {\em constructible} if
for any nonempty closed subscheme $Y\subseteq X$ there exists a nonempty
Zariski open subset $U\subseteq Y$ for which the restriction
$\V|_U$ is locally constant. The category of constructible sheaves
of $R$-modules on $X$ is denoted by $\Constr_R(X).$ 
\end{enumerate}}
\end{defn}

Let $\bar{x}:\Spec(\bar{k})\to X$ be a geometric point. The 
\'etale fundamental group of $X$ with
base point $\bar{x}$ is denoted by 
$\pi_1^\et(X,\bar{x})$ (see \cite{FreitagKiehl}, App. A, or \cite{MilneEC}).
There is a short exact
sequence \begin{equation}\label{eqsplitpione} 1 \To
\pi_1^\geom(X,\bar{x})\To \pi_1^\et(X,\bar{x})\To G_{k}\To 1\,,\end{equation}
where $\pi_1^\geom(X,\bar{x})$ denotes 
the {\em geometric fundamental group} 
$\pi_1^\et(X_{\bar{k}},\bar{x}).$
Moreover, if $\bar{x}$ is defined over $k,$  
then the sequence \eqref{eqsplitpione}
splits.\\

 Let
$\Rep_R(\pi_1^\et(X,\bar{x}))$ denote the category of continuous
representations $$\pi_1^\et(X,\bar{x}) \To \GL(V)\, ,$$ 
where $V\simeq R^n$ for
some $n\in \NN.$ If $\rho\in \Rep_R(\pi_1^\et(X,\bar{x})),$ then 
$\rho^\geom$ denotes the restriction 
of $\rho$ to the { geometric fundamental group} 
$\pi_1^\geom(X,\bar{x}).$\\

The following proposition is well known, see
\cite{FreitagKiehl}, Prop. A.1.8:

\begin{prop}\label{propequivofetallocal} 
The map $ \V\mapsto \rho_\V$ induces an an equivalence of categories
$$ \LS_R^\et(X)\quad  \cong \quad \Rep_R(\pi_1^\et(X,\bar{x}))\,.$$
\end{prop}

Now suppose that $k\subset\CC$ is a subfield of the complex
numbers.  The set of $\CC$-rational points of $X$ has a canonical
structure of a complex manifold which is denoted by $X^{\rm an}$.
Moreover, there is a functor $\F\mapsto\F\an$ from \'etale 
sheaves (of $R$-modules) on $X$ to sheaves on $X\an$, called {\em
analytification} (see e.g.\ \cite{FreitagKiehl}, Chap. I.11). 

If $\V$
is an \'etale local system on $X$ corresponding to a
representation $$\rho_\V:\pi_1^\et(X,\bar{x})\To\GL(V)\,,$$ then the
analytification $\V\an$ of $\V$ is the local system corresponding
to the composition of $\rho_\V$ with the natural homomorphism
$\pi_1(X^{\rm an},\bar{x})\to\pi_1^\et(X,\bar{x})$.\\
 
 Let $k \subseteq \CC$ and let 
 $$U_0=\AA^1_k\setminus \uo,\quad  \uo=\{u_1,\ldots,u_r\} \in \OO_r(k)\,,$$ be as above.
Fix
generators  $\al_1,\ldots,\al_{r+1}$ of $\pi_1(U_0\an,\bar{u}_0)$ as
in Section~\ref{braidgroups}. The natural map
$$\iota :\pi_1(U_0\an,\bar{u}_0)\To \pi_1^\et(U_0,\bar{u}_0)$$ is an injection under which
$\alpha_i$ maps to a generator of the inertia group of the missing
point $u_i$ (see e.g. \cite{MalleMatzat}). Moreover,
$$\pi_1^\geom(U_{0},\bar{u}_0)=\widehat{\langle
\iota(\al_1),\ldots,\iota(\al_{r+1}) \rangle}\, ,$$ where
$\,\widehat{\;}\,$ denotes the profinite closure.

\begin{defn}{\rm For $\V\in \LS_R^\et(U_0),$  let 
\begin{eqnarray}\nonumber
T_\V&:=&\left(\rho_\V(\iota(\al_1)),\ldots,\rho_\V(\iota(\al_{r+1}))\right)\in
\GL(V)^{r+1} \end{eqnarray} be the {\em associated tuple} of $\V.$}
\end{defn}

The following remark is immediate:

\begin{rem}\label{remanaletal} 
$$ T_\V =T_{\V^{\rm an}}\,.$$
\end{rem}

\subsection{Galois covers and  fundamental groups.}\label{seccovv}

The proofs of the statements in this section can be found in
\cite{FreitagKiehl}, App.~I.\\

 Let $k$ be a subfield of $\CC,$  let $X$ be 
 a smooth and geometrically irreducible 
variety over $k,$ and let $x$ be a geometric point of $X.$ Any
finite  \'etale Galois cover $f:Y\to X$ with Galois group
$G=G(Y/X)$ corresponds to  a surjective homomorphism   $\Pi_f:
\pi_1^\et(X,x)\to G.$

Let $U_0:=\AA^1_k\setminus \uo$ and let 
$$\pi_1^\geom(U_{0},\bar{u}_0)=\widehat{\langle
\iota(\al_1),\ldots,\iota(\al_{r+1}) \rangle}$$ be as in Section
\ref{secetaledef}. 
 Let  $f : Y_\bk \to U_{0,\bar{k}}$ be a finite \'etale Galois cover.
Then
$$ \begin{array}{rcl} f &\longleftrightarrow &
\Pi_f : \pi_1^\geom(U_{0},\bar{u}_0)\to G\\
&\longleftrightarrow&
\g_f:=\left(g_1:=\Pi_f(\iota(\al_1)),\ldots,g_{r+1}:=\Pi_f(\iota(\al_{r+1}))
\right)\in
G^{r+1}\,.\end{array}$$
If $f : Y \to U_{0,{k}}$ is a finite \'etale Galois cover, then we
define  $\g_f$ to be the tuple corresponding to $f_{\bar{k}}:Y_\bk
\to U_{0,\bar{k}}.$\\

Let $R$ be as in Section~\ref{secetaledef}, let $V\simeq R^n,$ and let
$\chi:G \hookrightarrow \GL(V)$ be a representation. If
$f:Y\to U_{0,k}$ is an \'etale  Galois cover with Galois group $G$, then
$\L_{(f,\chi)}\in \LS_R^\et(U_{0,k})$ 
denotes the local system associated to the composition
$\chi\circ \Pi_f.$

\begin{rem}\label{remaha}
 Any \'etale local system $\V$ with finite monodromy (i.e.,
$\im(\rho_\V)$ is finite) arises as a local system $\L_{(f,\chi)},$
where $f$ and $\chi$ are as above.
\end{rem}

\subsection{Tate twists.}\label{sectatetwi}

Using the concept of Tate twists,  one can control the 
eigenvalues of Frobenius elements. 

\begin{defn}\label{deftwist} {\rm 
Let $k$ be a field and let $\ell$ be a prime which is different from 
the characteristic of $k.$
\begin{enumerate}
\item Consider the $G_k$-module
$$ \ZZ_\ell(1)=\lim_{\leftarrow}  \mu_{\ell^n} \, ,$$
where $G_k$ acts via the cyclotomic character (see 
Section~\ref{secgeneralnotation}) and let 
$R$ be as in Section~\ref{secetaledef}. We set
$$R(1):=\ZZ_\ell(1)\otimes_{\ZZ_\ell} R\,,\quad 
R(-1):=\Hom(R(1),R)\,,\quad  R(0)=R(1)\otimes
R(-1)\,,$$ and for $n\in \NN_{>0},$ we set
$$ R(n):=R(1)\otimes \cdots \otimes R(1)\quad \text{\rm
($n$ times)} $$ and
$$ R(-n):=R(-1)\otimes \cdots \otimes R(-1)\quad \text{\rm
($n$ times)}\,.$$
\item  For $X\in
\Var_k$ and $\V\in \Sh_R^\et(X)$ we set
$$\V(n):=\V \otimes_R {R(n)}\,.$$
The sheaf $\V(n)$ is called the {\em $n$-th Tate twist} of $\V.$

\item Let $V$ be a free $R$-module of finite rank 
and let 
$\rho:G_k\to \GL(V)$ be a Galois representation. Then the $n$-th {\em 
Tate twist} of the $G_k$-module $V$  is defined to be 
$V(n):=V\otimes_R R(n)$
(where $R(n)$ is viewed as a $G_k$-module via the construction
 above).
The $n$-th {\em 
Tate twist}  of $\rho$ is 
 the corresponding representation 
$\rho(n): G_k\to \GL(V(n)).$
\end{enumerate}}
\end{defn}

\begin{rem}
\begin{enumerate}
\item If $R\subseteq \bar{\QQ}_\ell,$ taking the 
$n$-th Tate twists
of $\rho:G_k\to \GL(V)$ amounts to tensoring $\rho$ 
by the $n$-th power of the cyclotomic character $\chi_\ell:$
\begin{equation}\nonumber \rho(n)=\rho\otimes \chi_\ell^n\,.\end{equation}
\item If $R\subseteq \bar{\FF}_\ell,$ this amounts to tensoring $\rho$ 
by the $n$-th power of the mod-$\ell$-cyclotomic character $\bar{\chi}_\ell.$
\end{enumerate}
\end{rem}

\section{Variation of parabolic cohomology and 
the  middle convolution}\label{sectconvdef1}

In this section we start by briefly recalling the results of \cite{dw03} on
the parabolic cohomology of local systems. 
 This leads to the definition of the middle convolution
in Section~\ref{sectconvdef}. \\

We will freely use the notation introduced in the
last sections.
As in Section~\ref{braidgroups},  we
will write $\AA^1,\, \PP^1 , \ldots$ instead of $\AA^1(\CC),\,
\PP^1(\CC), \ldots$ and  view these objects equipped with their
associated topological and complex analytic structures.

\subsection{Cohomology of local systems on $U_0$.}\label{secparabcohom11}
Let $U_0=\AA^1\setminus \uo$ be as in Section~\ref{braidgroups} and 
let $R$ be a commutative Ring with a unit.
Let $\V_0\in \LS_R(U_0)$ 
and let 
$$T:=T_{\V_0}=(T_1,\ldots,T_{r+1})\in \GL(V)^{r+1}$$
denote the associated tuple as in Notation 
\ref{notfortupeloflocsys}. It is shown in \cite{dw03} that the group
$H^1(U_0,\V_0)$ is  isomorphic to $H_T/E_T,$ where
\begin{equation}\nonumber H_T:=\{ (v_1,\ldots,v_{r+1}) \in V^{r+1}\mid v_1(T_2\cdots
T_{r+1})+ v_2(T_3\cdots T_{r+1})+\cdots +v_{r+1}=0 \}\end{equation}
and
\begin{equation} \nonumber E_T:=\{ \left( v(T_1-1),\ldots,v(T_{r+1}-1)\right)\mid v \in V\}\,.\end{equation}
(The isomorphism is given by the composition of the natural isomorphism
$$H^1(U_0,\V_0) \to H^1(\pi_1(U_0,u_0),V)$$ with the evaluation map,
which associates to the equivalence class of a crossed
homomorphism $[\delta]\in H^1(\pi_1(U_0,u_0),V)$ the corresponding
equivalence class of $[(\delta(\al_1),\ldots,\delta(\al_{r+1}))]$ in
$V^{r+1}/E_T.$)

Let $j:\U_0 \to \PP^1$ be the natural inclusion. It is shown
in {loc.~cit.} that the parabolic cohomology group
$H^1_p(U_0,\V_0):=H^1(\PP^1,j_*(\V_0))$ is isomorphic to
$U_T/E_T,$ where
\begin{equation} \nonumber U_T:=\{ (v_1,\ldots,v_{r+1}) \in
H_T \mid  v_i \in \im(T_i-1),\, i=1,\ldots,r+1 \}\,.\end{equation}
Here, 
the additional relations arise from the natural isomorphism
$$H^1_p(U_0,\V_0)\simeq \im\left(H^1_c(U_0,\V_0)\to H^1(U_0,\V_0)\right)\,.$$

\subsection{Variation of parabolic cohomology.}\label{secvarpar}
 Let $S$ be a  connected complex analytic
 manifold, let $X:=\PP^1_S=\PP^1\times S,$
and let $\D\subseteq X$ be a smooth relative divisor of degree $r+1$ over $S$ which
contains the section $\{\infty\} \times S.$ Let $U:=X \setminus
\D,$ let $j: U \to X$ be the natural inclusion, let $\bar{\pi}: X \to S$ be 
the
projection onto $S,$ and let $\pi: U \to S$ be the restriction of
$\bar{\pi}$ to $U.$ Let further $s_0 \in S$ and let 
$U_0:=\pi^{-1}(s_0).$\\

 A local system $\V\in
\LS_R(U)$ is called a {\em variation of} $\V_0\in \LS_R(U_0)$
over $S$ if
$\V_0 = \V|_{U_0}.$ The {\em parabolic cohomology} of this
variation is by definition the first higher direct image
$\W:=R^1\bar{\pi}_*(j_*\V).$ It is a local system on $S$ whose
stalk $\W_{s_0}$ is canonically isomorphic to the parabolic
cohomology group $H^1_p(U_0,\V_0)$ (see loc. cit). Thus $\W$ corresponds to its
monodromy representation
$$\rho_\W:
\pi_1(S,s_0)\To \GL(H^1_p(U_0,\V_0))\cong \GL(U_T/E_T) \, ,$$ where
$T:=T_{\V_0}$ is the associated tuple  of
$\V_0$ and $U_T$ and $E_T$ are as in the last section. 

We want to 
determine the representation $\rho_\W.$ To this end, let 
 $\beta_1,\ldots,\beta_{r-1}$ denote the standard
generators of $\A_r$ (see Section \ref{braidgroups}). 
Consider linear automorphisms $\Phi(T,\beta_i)$ of $V^{r+1}$ which
are defined as follows:
\begin{multline} \label{locals5eq6}
 \qquad (v_1,\ldots,v_{r+1})^{\Phi(T,\beta_i)} \\
    \;=\;  (v_1,\ldots,v_{i-1},\,v_{i+1},\,
    \underbrace{v_{i+1}(1-T_{i+1}^{-1}T_iT_{i+1})+v_iT_{i+1}}_{
        \text{\rm $(i+1)$th entry}},\,v_{i+2},\,\ldots,v_{r+1})\,.
\end{multline}
These automorphisms multiply by the following rule:
\begin{equation} \label{locals5eq7}
    \Phi(T,\beta)\cdot\Phi(T^{\beta},\beta') \;=\;
       \Phi(T,\beta\beta')\,.
\end{equation}

It is easy to see that the spaces $U_T$ and $E_T$ are mapped
under $$\Phi(T,\phi(\gamma)),\quad  \gamma \in  \pi_1(S,s_0)\, ,$$ 
isomorphically to the spaces $U_{T^{\phi(\gamma)}},$ and
$E_{T^{\phi(\gamma)}},$ respectively (where $\phi(\gamma)$ is as in Section 
\ref{braidgroups} and acts  as in
\eqref{braidaction} on $\GL(V^{r+1})^r$). Let
$$\bar{\Phi}(T,\phi(\gamma)):U_T/E_T  \To
U_{T^{\phi(\gamma)}}/E_{T^{\phi(\gamma)}}$$ be the isomorphism
induced by $\Phi(T,\phi(\gamma)).$ The next result immediately follows
from \cite{dw03}, Thm. 2.5 and Rem. 2.6 (using the
above diagrams \eqref{splitdiag} and \eqref{rowsplitting}):

\begin{prop}\label{parmono} Assume that 
 $\rho_\V(\pi_1(S,s_0))=\{1\},$ where 
 $\pi_1(S,s_0)$ is viewed as a subgroup 
of $\pi_1(U,(u_0,s_0))$ as in Section~\ref{braidgroups}.
 If the setup is chosen in such a way that
\eqref{rowsplitting} holds, then
$$ \rho_\W(\gamma)=\bar{\Phi}(T,\phi(\gamma)), \quad \forall
\gamma \in\pi_1(S,s_0)\,.$$
\end{prop}

\begin{prop}  \label{prop1}\begin{enumerate}
 \item {\rm (Ogg-Shafarevich)}
 Suppose that $R=K$ is a field and that the stabilizer $V^{\pi_1(U_0)}$ is
  trivial. Then
  \[
    \rk(\W)= \dim_K H^1\para(U_0,\V_0) \;=\;
        (r-2)\dim_K V - \sum_{i=1}^r \dim_K{\rm Ker}(T_i-1)\,.
  \]

\item {\rm (Poincar\'e Duality)}
Let $\V\otimes\V\to{R}$ be a non-degenerate symmetric
(resp.\
  alternating) bilinear pairing of sheaves corresponding to an injective
  homomorphism $\kappa:\V\inj\V^*$ with $\kappa^*=\kappa$ (resp.\
  $\kappa^*=-\kappa$) and  let $ \W \;:=\; R^1\pib_*(j_*\V).$
Then the cup product
  defines a non-degenerate alternating (resp.\ symmetric) bilinear
pairing of sheaves $\W \otimes \W \to  {R}.$
\end{enumerate}
\end{prop}

\proof Claim (i) is \cite{dw03}, Rem. 1.3. See 
\cite{dw04} for (ii).
 \Endproof


\subsection{The definition of the middle convolution.}\label{sectconvdef}

Let $\OO_r,\,r\in \NN_{>0},$ denote the configuration space 
of subsets of $\AA^1$ having cardinality $r$ (see 
Section~\ref{braidgroups}).\\

For $\uo :=\{ x_1,\ldots,x_p\} \in \OO_p$ and $\vo
:=\{y_1,\ldots,y_q\} \in \OO_q$ set
$$ \uo \ast \vo := \{ x_i + y_j \mid i=1,\ldots,p, \,\, j = 1,\ldots
,q \}\,.$$ Let $U_1:=\AA^1\setminus \uo,$ $U_2:=\AA^1\setminus \vo$
and $S:=\AA^1\setminus \uo \ast \vo.$ Set
$$\tilde{f}(x,y):= \prod_{i=1}^{p} (x-x_i) \prod_{j=1}^q(y-x-y_j)
\prod_{i,j} (y-(x_i+y_j)) $$ and let $f\in \CC[x,y]$ be the
associated reduced polynomial. One has
 $\tilde{f}= f$ if and only if $|\uo \ast \vo| = p\cdot q,$ in which
 case we call $\uo \ast \vo$ {\em generic}.
Let
$$ \tilde{\bf w}:=\{ (x,y)\in \AA^2 \mid f(x,y) =0\}$$
and let $U:=\AA^2 \setminus \tilde{\bf w}.$

\vspace{.5cm}
\begin{center}
\includegraphics{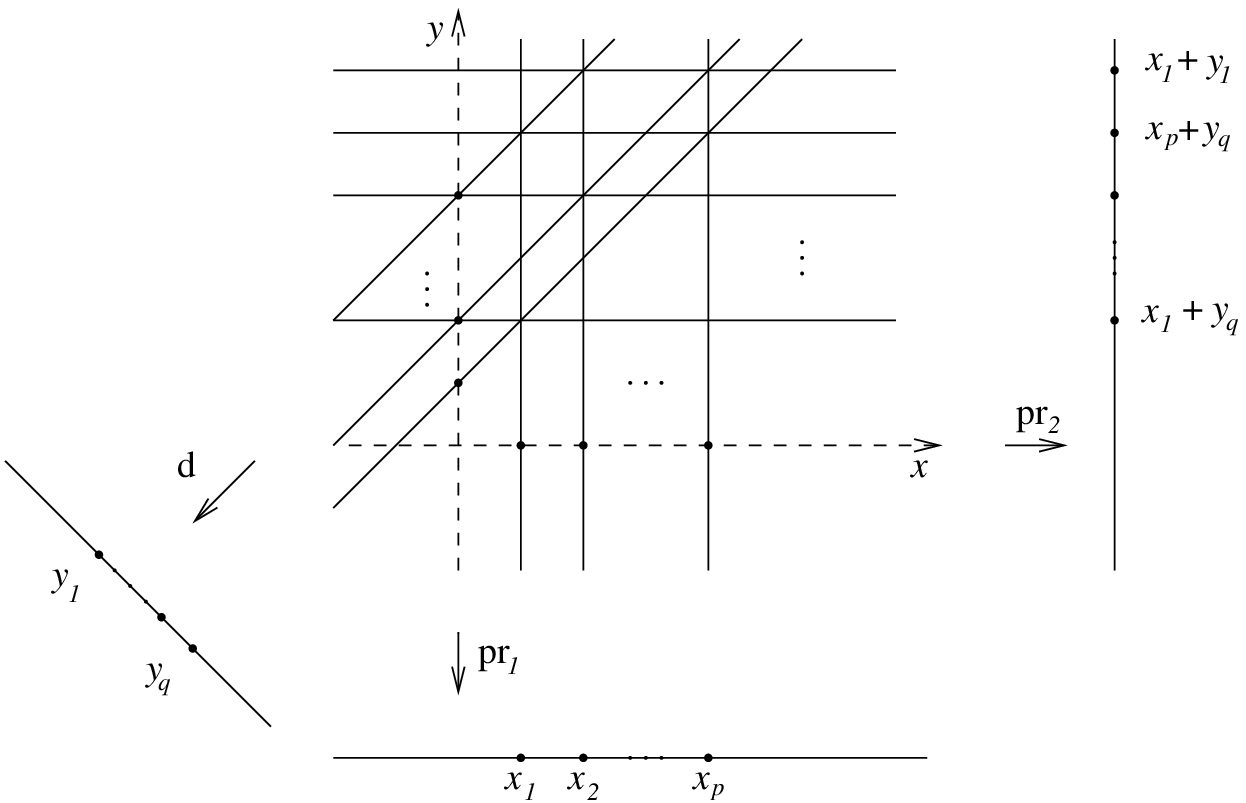}

\vspace{.5cm}
{\bf Figure 1:} The complement $U=\AA^2\setminus \tilde{\bf w}$ and its projections.
\vspace{.5cm}
\end{center}

The set $U$ is equipped with three
maps which play an important role in all that follows: One the one
hand, there are the projections
$$\pr_1: U \To U_1,\quad  (x,y) \Mapsto x$$
and
$$ \pr_2: U \To S,\quad (x,y) \Mapsto y\,.$$
On the other hand, there are the { subtraction map}
$$ \q: U \To U_2,\quad  (x,y) \Mapsto y-x\,.$$
Let 
 $$j:U \To X:=\PP^1_{S}, \quad (x,y)\Mapsto ([x,1],y)$$ and 
let 
$ {\bf w}:=X\setminus U.$ 
 Since ${\bf w}$
 is a smooth relative divisor of degree 
$p+q+1$ over $S,$  we are in
the situation of Section~\ref{secvarpar} with 
$r=p+q$ and 
$\pi=\pr_2.$  The
 second projection $X=\PP^1_{S} \to S$ is
denoted by $\opr_2.$ The fibre $\pr_2^{-1}(y_0)$ is denoted by
$U_0.$ The first projection $X=\PP^1_{S} \to \PP^1$
yields an identification of $U_0$ with
$\AA^1 \setminus (\uo \cup (y_0 - \vo)),$ where
$$\uo \cup (y_0 - \vo):=\uo \cup \{ y_0-y_1,\ldots,y_0 - y_q\}
\in \OO_{p+q}\,.$$

Let  $\V_1 \in \LS_R(U_1)$ and let 
$\V_2 \in \LS_R(U_2).$ The local system
$$\V_1 \tim \V_2:= \pr_1^*\V_1\otimes \q^*\V_2$$
is a local system on $U$ which is a variation of $\V_1 \tim
\V_2|_{U_0}$ over $S.$ The middle convolution of $\V_1$ and $\V_2$ 
is now defined to be the parabolic cohomology of this 
variation:

\begin{defn}{\rm  The {\em middle convolution} of $\V_1 \in \LS_R(U_1)$ and
$\V_2 \in \LS_R(U_2)$ is the local system
$$ \V_1 \ast \V_2 := R^1 (\opr_2)_*(j_*(\V_1\tim \V_2)) \in \LS_R(S)\,.$$

The {\em $\ast$-convolution} of $\V_1 \in \LS_R(U_1)$ and $\V_2
\in \LS_R(U_2)$ is the local system
$$ \V_1 \ast_\ast \V_2:=R^1 (\pr_2)_*(\V_1\tim \V_2) \in \LS_R(S)\,.$$}
\end{defn}

\begin{rem}\label{remmclambda} \begin{enumerate}

\item   In \cite{Katz96},
N. Katz gives a similar construction in a more general 
category of complexes of  sheaves, see \cite{Katz96} and 
\cite{KiehlWeissauer}. 
The approach of Katz has many advantages, since it measures the part
of the middle convolution  which 
appears at the points of infinity of $S.$ 
But in most applications  it  
suffices to deal with the middle convolution of local systems. 
Since the
category of local systems is 
 accessible to explicit computation via the concept
of monodromy,  it is often actually necessary to deal with local 
systems (this is the reason why we consider only the convolution 
of local systems here).
Also, 
 there is a dictionary which translates Katz' approach (over an open subset) 
to the  middle
convolution considered here (see \cite{Katz96}, Chap. 2.8).

\item An important case of the middle convolution is Katz' middle
convolution functor $\MC_\chi,$ see \cite{Katz96}: Let $\chi$ be a
character of $\pi_1({\mathbb G}_m),\,{\mathbb G}_m=\AA^1\setminus\{0\},$ 
and let $\V_\chi \in
\LS_R({\mathbb G}_m)$ be the associated local system. We call
$\V_\chi$ the {\em Kummer sheaf} associated to $\chi.$ 
 Then one obtains a functor
$$ \LS_R(U_1) \To \LS_R(U_1) ,\, \V \Mapsto \V \ast \V_\chi\,.$$
In \cite{dr03} it is shown that this functor coincides with Katz'
middle convolution functor $\MC_\chi$ if one restricts it to the
category of convolution sheaves (see Def.~\ref{defconvosheaf}
below).

\item Let $\ga$ be a counterclockwise generator of $\pi_1({\mathbb
G}_m)$ and let $\lambda=\chi(\ga),$ where $\chi$ is as in (ii). Then one
obtains a transformation of tuples $T_\V\mapsto T_{\V\ast \V_\chi}$
which corresponds to the tuple transformation $\MC_\lambda$
considered in \cite{dr00}, see \cite{dr03}. In Section 
\ref{secirrkumm} we give an alternative proof of this statement,
 see Rem.~\ref{remlab}.

\end{enumerate}

\end{rem}

\section{Properties of the 
convolution}\label{secfirstpropsing}

\subsection{First properties.}\label{secfirstpropsing1}

In this section we collect some facts about the middle convolution
of local systems which will be useful in later applications.\\

Let $\V_1 \in \LS_R(U_1)$ and $\V_2 \in \LS_R(U_2),$
where $U_1=\AA^1\setminus \uo$ and $U_2=\AA^1\setminus \vo$ are 
as in the last section. Let us fix a base point
$(x_0,y_0)$ in $U.$ This induces base points
$x_0=\pr_1(x_0,y_0),\,y_0-x_0=\q(x_0,y_0), \, y_0=\pr_2(x_0,y_0) $ of $U_1,\,U_2,$
and $S=\AA^1\setminus \uo\ast \vo,$ respectively. Let
 $V_1$ denote  the stalk of $\V_1$
at $x_0$ and   let $V_2$ denote 
the stalk of $\V_2$ at $y_0-x_0.$ Let also 
  $$U_0=\pr_2^{-1}(y_0)=\AA^1 \setminus \uo \cup (y_0 - \vo)$$ be as in the last section. \\

 The representation
$\rho_{\V_1\circ \V_2|_{U_0}}: \pi_1(U_0,(x_0,y_0)) \to \GL(V_1\otimes V_2)$
factors as
\begin{equation}\label{rhocirc}\rho_{\V_1\circ \V_2|_{U_0}}=(\rho_{\V_1}\otimes
\rho_{\V_2})\,\circ\,(\pr_1 \times \q)_* \, ,\end{equation} where
$$(\pr_{1} \times \q)_*:\pi_1(U_0,(x_0,y_0))\To \pi_1(U_1,x_0)\times
\pi_1(U_2,y_0-x_0)$$ is the map which is induced 
by $\pr_{1}|_{U_0} \times \q|_{U_0}.$
 Let $\al_1,\ldots,\al_{p+q}$ be generators of
$\pi_1(U_0,(x_0,y_0))$ which are chosen as in Figure 2 in Section \ref{secbasicsetup}
below. Let
$$ \ga_1:=\pr_{1\ast}(\al_1)\,,\,\ldots\,,\,\ga_{p}:=\pr_{1\ast}(
\al_p)$$ be the induced generators of $\pi_1(U_1,x_0)$ and
let 
$$\eta_1:=\q_*(\al_{p+1})\,,\,\ldots\,,\,\eta_{q}:=\q_*(\al_{p+q})$$ be 
those
of $\pi_1(U_2,y_0-x_0).$ Let $\ga_{p+1}:=(\ga_1\cdots \ga_p)^{-1}$
and let $\eta_{q+1}:=(\eta_1\cdots \eta_q)^{-1}.$
With respect to the generators $\ga_1,\ldots,\ga_{p+1},$
and $\eta_1,\ldots,\eta_{q+1},$ let
$$T_{\V_1}=(A_1,\ldots,A_{p+1})\in \GL(V_1)^{p+1}$$ and
$$T_{\V_2}=(B_1,\ldots,B_{q+1})\in \GL(V_2)^{q+1}$$ be the associated 
tuples (respectively).

 It follows from our choice of homotopy generators and
 \eqref{rhocirc}
 that
\begin{multline} \label{Tcirc}
T_{\V_1\circ \V_2|_{U_0}} = (C_1= A_1\otimes 1_{V_2}\,,\,\ldots,\,
C_p= A_p\otimes 1_{V_2}\,,\\
 C_{p+1}=1_{V_1}\otimes B_1\,,\,\ldots\,,\, C_{p+q}=1_{V_1}\otimes 
B_q\,,\,C_{p+q+1}=A_{p+1}\otimes B_{q+1})\,.
\end{multline}

\begin{prop}\label{dimensione} Suppose that $R=K$ is a field and
that one of the stabilizers
$$V_1^{\pi_1(U_1)}\quad\quad {  and}\quad\quad
V_2^{\pi_1(U_2)}$$ is trivial. Let $\dim_K V_i=n_i.$
Then
\begin{multline}\label{dimensionconv} \rk(\V_1\ast \V_2)= (p+q-1)n_1 n_2 -
\sum_{i=1}^p n_2\dim_K
\ker(A_i - 1_{V_1}) \\
-\sum_{j=1}^q n_1\dim_K \ker(B_j - 1_{V_2})
-\dim_K\ker(A_{p+1}\otimes B_{q+1}-1_{V_1\otimes
V_2})\,.\end{multline}
\end{prop}

\proof It follows from \eqref{Tcirc} and the properties of the
tensor product that
$$\dim_K \ker(C_i-1_{V_1\otimes V_2})= n_2\dim_K
\ker(A_i - 1_{V_1}),\;i=1,\ldots, p\,,$$ and
$$\dim_K \ker(C_i-1_{V_1\otimes V_2})= n_1\dim_K
\ker(B_i - 1_{V_2}),\; i=p+1,\ldots, p+q\,.$$ The claim now follows
from Prop.~\ref{prop1} (i). \Endproof

\begin{rem} The dimension $\dim_K\ker(A_{p+1}\otimes B_{q+1}-1_{V_1\otimes
V_2})$ can be easily computed using Lemma~\ref{lemkroneckerjordan}.
\end{rem}

\begin{rem}\label{remalternative} There is an alternative description of the middle convolution which
is sometimes useful: The linear automorphism
\begin{equation*} \AA^2 \To \AA^2,\; x\Mapsto x,\; y \Mapsto y+x
\end{equation*} leaves the first projection
$\pr_1$ unchanged but it transforms $\q$ into $\pr_2$ and $\pr_2$
into the addition map
$$ {\rm a}: \AA^2 \To \AA^1,\, (x,y) \Mapsto x+y\,.$$
After adapting the notation in the obvious way one obtains a
natural isomorphism
$$ \V_1\ast \V_2 \cong R^1 {\bar{\rm a}}_*(\tilde{j}_*(\V_1 \tim
\V_2))\, ,$$ where $\tilde{j}$ is the compactification in the
direction $x+y.$\end{rem}

Using the last remark and the  coordinate switch $y\mapsto x,\,
x \mapsto y,$ one obtains the following result:

\begin{prop}\label{propconviscommutative}
\begin{equation}\label{convcommn}
\V_1\ast \V_2 \cong \V_2\ast \V_1 \,.\end{equation}\Endproof
\end{prop}
 Let
$\kappa_1: \V_1 \otimes \V_1 \to {R}$ and $\kappa_2:\V_2
\otimes \V_2 \to {R}$ be bilinear pairings. By Prop.~\ref{prop1}~(ii) 
and by the definition of the middle convolution as
the parabolic cohomology of a variation, one obtains a bilinear
pairing \begin{equation}\label{bilinearpairofconv} \kappa_1 \ast
\kappa_2: (\V_1\ast \V_2)\otimes(\V_1\ast \V_2) \to
{R}\,.\end{equation} Let $\sym(\kappa)$ be as in 
Section~\ref{localsyst}.

\begin{prop} \label{fform} $$ \sym({\kappa_1 \ast
\kappa_2}) = -\, \sym({\kappa_1})\cdot \sym({\kappa_2}). $$
\end{prop}

\proof It is easy to see (using \eqref{form}) that the sheaf
$\V_1\circ \V_2\in \LS_R(U)$ carries a bilinear pairing
$\kappa_1\circ \kappa_2$ such that
$$\sym({\kappa_1\circ \kappa_2})=\sym({\kappa_1})\cdot
\sym({\kappa_2})\,.$$ The result then follows from Prop.
\ref{prop1} (ii). \Endproof

\begin{lem}\label{naive} 
Let $R=K$ be a field and  let $\V_2 \in \LS_K(U_2)$ 
such that the stabilizer of the fundamental group
on the stalk is zero, i.e., $V_2^{\pi_1(U_2)}=0.$  Let
further $S^o\subseteq S$ be an open
subset. Then the following holds:
\begin{enumerate}
\item The sheaves $R^i(\pr_2)_*(\V_1\tim \V_2)$
vanish for any local system $\V_1\in\LS_K(U_1)$ and 
 $i\not=1.$
\item The functor
$$ \LS_K(U_1) \To \LS_K(S^o),\; \V_1 \Mapsto (\V_1 \ast_\ast \V_2)|_{S^o}\,,$$
is exact\,.\end{enumerate}
\end{lem}

\proof  Let $\FFF$ be any
local system on $U_0,$ where $U_0=\AA^1\setminus \uo\cup (y_0-\vo)$ is as above.
It is well known that $H^2(U_0,\F)=0$ for 
any locally constant sheaf $\F$ (this follows from the fact
that $U_0$ is affine). Thus, by  
\eqref{eqhigh}, 
$$R^2({\pr}_2)_*(\V_1\circ \V_2)=0\,.$$
By the 
properties of the tensor product (see Equation \eqref{rhocirc}), 
the condition 
$V_2^{\pi_1(U_1)}=0$ implies that 
$$H^0(U_0,\V_1\circ \V_2|_{U_0})=(V_1\otimes V_2)^{\pi_1(U_0)}=0\,.$$
Thus, again by 
\eqref{eqhigh}, $R^0({\pr}_2)_*(\V_1\circ \V_2|_{U_0})=0.$ This proves 
claim (i).

For (ii): It suffices to prove the claim for the stalks. 
Let
$$ 0 \To \V_1' \To \V_1 \To \V_1'' \To 0$$
be an exact sequence of local systems on $U_1.$ One obtains an
exact sequence
$$ 0 \To (\V_1'\circ \V_2)|_{U_0} \To (\V_1 \circ \V_2)|_{U_0}\To
 (\V_1''\circ \V_2)|_{U_0} \To 0 $$
of local systems on $U_0.$ 
 The long exact cohomology
sequence yields an exact sequence \begin{multline}
 H^0(U_0,(\V_1''\circ \V_2)|_{U_0})\To H^1(U_0,(\V_1'\circ
 \V_2)|_{U_0})\To
 H^1(U_0,(\V_1\circ \V_2)|_{U_0})\\
\To H^1(U_0,(\V_1''\circ \V_2)|_{U_0})\To
 H^2(U_0,(\V_1'\circ \V_2)|_{U_0})\,.\end{multline}
From (i), one obtains
$$H^0(U_0,(\V_1''\circ \V_2)|_{U_0}) =0=H^2(U_0,(\V_1'\circ
\V_2)|_{U_0})\,,$$ and the claim follows. \Endproof

\begin{prop}\label{propexactconv} 
We use the notation of Lemma~\ref{naive}. Assume that for 
$\V_2 \in \LS_K(U_2),$ the entries of the associated tuple 
$T_{\V_2}\in \GL(V_2)^{q+1}$ generate an absolutely irreducible and non-trivial 
subgroup of $\GL(V_2).$ 
 Then the following holds:
\begin{enumerate}
\item The sheaves $R^i(\overline{\pr}_2)_*(j_*(\V_1\tim \V_2))$
vanish for any local system $\V_1\in\LS_K(U_1)$ and 
 $i\not=1.$
\item The functor
$$ \LS_K(U_1) \To \LS_K(S^o),\; \V_1 \Mapsto \V_1 \ast \V_2|_{S^o}\,,$$
is exact.\end{enumerate}
\end{prop}

\proof Let $\FFF$ be any
local system on $U_0.$
 It is well known that  $H^0(\PP^1,j_*\FFF)$ is isomorphic  to
 the module of invariants $(\FFF_x)^{\pi_1(U_0,x)}$ and that
$H^2(\PP^1,j_*\FFF)$ is isomorphic to
 the module of coinvariants 
$$({\FFF_x})_{\pi_1(U_0,x)}={\FFF_x}/\langle fg-f\mid f \in {\FFF_x},\,
g \in \pi_1(U_0,x)\rangle\, ,$$ see e.g.
\cite{Looijenga92}, Lemma 5.3. 

With
the above identification of $H^0$ and $H^2,$ the 
irreducibility and non-triviality assumption, 
 and the
properties of the tensor product, one obtains
$$H^0(\PP^1_{y_0},j_*(\V_1\circ \V_2)|_{\PP^1_{y_0}})) =0=H^2(\PP^1_{y_0},j_*(\V_1\circ
\V_2)|_{\PP^1_{y_0}}))$$ for any 
$\V_1\in \LS_K(U_1).$
Thus, by \eqref{eqhigh},
the sheaves $R^i(\overline{\pr}_2)_*(j_*(\V_1\tim \V_2))$
vanish for $i\not=1.$ This proves (i). 

For (ii): It suffices to prove the claim for the stalks. 
Let
$$ 0 \To \V_1' \To \V_1 \To \V_1'' \To 0$$
be an exact sequence of local systems on $U_1.$ One obtains an
exact sequence
$$ 0 \To \V_1'\circ \V_2 \To \V_1 \circ \V_2\To
 \V_1''\circ \V_2 \To 0 $$
of local systems on $U$ and thus an exact sequence
$$ 0 \To j_*(\V_1'\circ \V_2)|_{\PP^1_{y_0}} \To j_*(\V_1 \circ \V_2)|_{\PP^1_{y_0}}\To
 j_*(\V_1''\circ \V_2)|_{\PP^1_{y_0}} \To 0 \,.$$
 The long exact cohomology
sequence yields an exact sequence \begin{multline}
 H^0(\PP^1_{y_0},j_*(\V_1''\circ \V_2)|_{\PP^1_{y_0}})\To H^1(\PP^1_{y_0},j_*(\V_1'\circ \V_2)|_{\PP^1_{y_0}})\To
 H^1(\PP^1_{y_0},j_*(\V_1 \circ \V_2)|_{\PP^1_{y_0}})\\
\To H^1(\PP^1_{y_0},j_*(\V_1''\circ \V_2)|_{\PP^1_{y_0}})\To
 H^2(\PP^1_{y_0},j_*(\V_1'\circ \V_2)|_{\PP^1_{y_0}})\,.\end{multline}
From (i), one obtains
$$H^0(\PP^1_{y_0},j_*(\V_1''\circ \V_2)|_{\PP^1_{y_0}}) =0=H^2(\PP^1_{y_0},j_*(\V_1'\circ \V_2)|_{\PP^1_{y_0}})\,,$$ and the claim follows.
\Endproof

\subsection{The basic setup in the generic case.}\label{secbasicsetup}

It is the aim of this subsection to provide the 
setup for the next two subsections, where we derive some statements 
on the 
irreducibility and the local monodromy of the middle convolution.\\

In the situation and notation of Subsections~\ref{sectconvdef} and 
\ref{secfirstpropsing1}: 
Let $\V_1 \in \LS_R(U_1),$  $\V_2 \in \LS_R(U_2)$ and $\V_1\circ \V_2|_{U_0} \in \LS_R(U_0)$
with 
$T_{\V_1}=(A_1,\ldots,A_{p+1})\in \GL(V_1)^{p+1},$
$T_{\V_2}=(B_1,\ldots,B_{q+1})\in \GL(V_2)^{q+1}$ and
\begin{multline}
T_{\V_1\circ \V_2|_{U_0}} = (C_1= A_1\otimes 1_{V_2}\,,\,\ldots,\,
C_p= A_p\otimes 1_{V_2}\,,\\
 C_{p+1}=1_{V_1}\otimes B_1\,,\,\ldots\,,\, C_{p+q}=1_{V_1}\otimes 
B_q\,,\,A_{p+1}\otimes B_{q+1})\,.
\end{multline}
Throughout this and the following two subsections, we
assume   that $R=K$ is a field and that $\uo \ast \vo$
is  generic, i.e., 
$\uo \ast \vo \in \OO_{p+q}.$\\

 We want to describe
the monodromy of $\V_1\ast \V_2.$ We can assume (using a suitable marking
as in \cite{dw03}) that we are in the following situation: 
The sets
$\uo=\{x_1,\ldots,x_p\},\,\vo=\{y_1,\ldots,y_q\},\, \{y_0\}$ are
elementwise real and 
$$x_1<x_2<\ldots < x_p <
y_0-y_1<y_0-y_2<\ldots< y_0-y_q\,.$$ Moreover, we can assume that
\begin{equation}\label{differenceofvu}
|x_p-x_1|< |y_{i+1}-y_i|\quad {\rm for}\quad  i=1,\ldots,q-1. \end{equation}

Let us fix a base point $(x_0,y_0)$ of $U_0$ and of $U.$ 
We assume that the imaginary part of $x_0$ is large enough, i.e.,
larger  than  the maximal 
imaginary part of  $\delta_{i,j}(t),$ where $\delta_{i,j}$ is 
as in Figure 4 below. One obtains base points
$$x_0=\pr_1(x_0,y_0),\quad y_0-x_0=\q(y_0,x_0), \quad  y_0=\pr_2(y_0,x_0) $$
on $U_1,\,U_2$ and $S$ (respectively). We choose
generators $\al_1,\ldots,\al_{p+q}$ of $\pi_1(U_0,(x_0,y_0))$ as follows:

\vspace{.6cm}
\begin{center}
\includegraphics{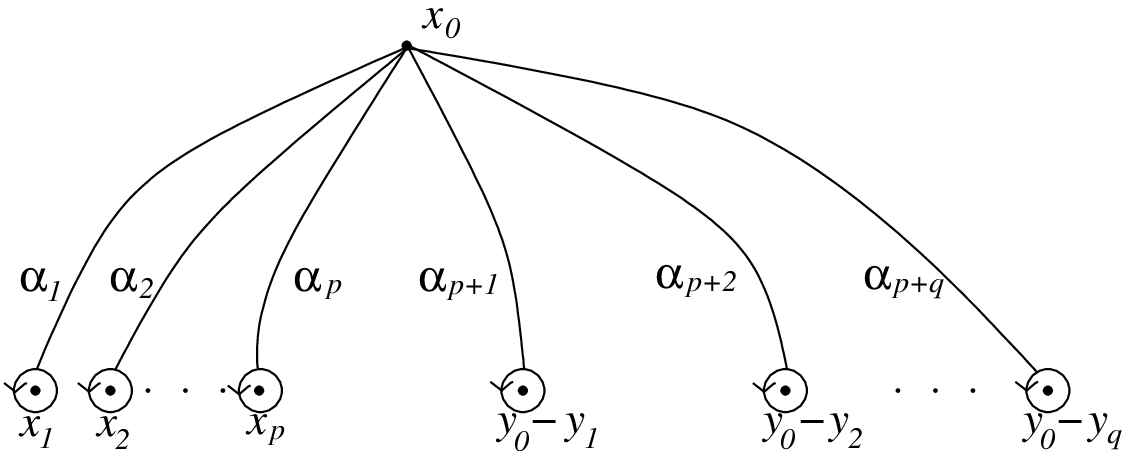}

\vspace{.5cm}
{\bf Figure 2:} The generators $\al_1,\ldots,\al_{p+q}$
\vspace{.5cm}

\end{center}

Next we choose generators 
$\beta_1,\ldots,\beta_{p+q-1}$ of 
$$\A_{p+q}=\pi_1(\OO_{p+q},\uo\cup\{y_0-y_1,\ldots, y_0-y_q\}) $$ as follows:

\vspace{.6cm}
\begin{center}
\includegraphics{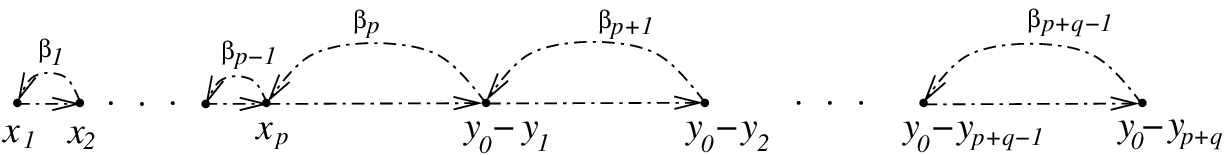}

\vspace{.6cm}
{\bf Figure 3:} The generators $\beta_1,\ldots,\beta_{p+q-1}$
of $\A_{p+q}=\pi_1(\OO_{p+q},\uo \cup
(y_0-\vo)) $
\vspace{.6cm}

\end{center}

Here, the generator $\beta_i$ is the path in $\OO_r$ which  fixes 
the points $$\{x_1,\ldots,x_{p+q}\}\setminus \{x_i,x_{i+1}\}$$ and which moves 
the point $x_i$ along the real axis to $x_{i+1}$ and $x_{i+1}$ to 
$x_i$ as indicated in  Figure 3.\\

Then we choose generators $\delta_{i,j},\,
i=1,\ldots,p,\,j=1,\ldots ,q$ of $\pi_1(S,y_0)$ as follows:

\vspace{.6cm}
\begin{center}
\includegraphics{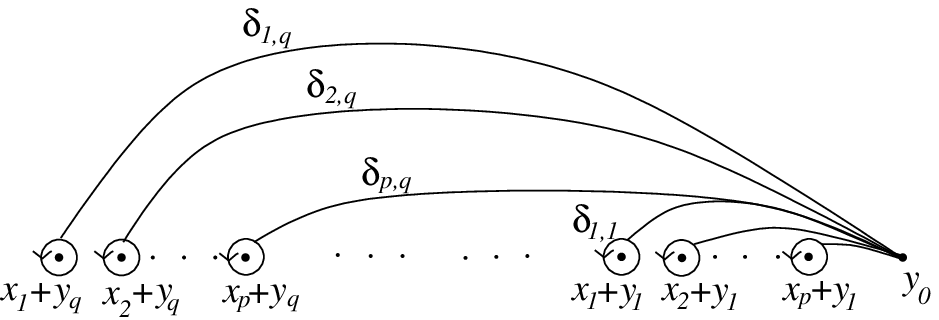}

\vspace{.6cm}
{\bf Figure 4:} The generators $\delta_{i,j}$
\vspace{.6cm}

\end{center}

Note that the product
\begin{multline}\label{eqproductdeltai}
\delta_{1,q}\,\delta_{2,q}\,\cdots\,
\delta_{p,q}\,\delta_{1,q-1}\,\delta_{2,q-1}\,\cdots\,
\delta_{p,q-1}\,\cdots\,\delta_{1,1}\delta_{2,1}\,\cdots\,
\delta_{p,1}
\end{multline}
is homotopic to a simple loop around $\infty.$\\

\begin{prop}\label{propdeltt}
Let $$\phi:\pi_1(S,y_0) \to \A_{p+q}=\pi_1(\OO_{p+q},\uo \cup
(y_0-\vo)) $$ be as in Section~\ref{braidgroups}. Then
\begin{equation}\label{eqimageofdeltai} 
\phi(\delta_{i,1})=\beta_{i,p+1},\; i=1,\ldots,p\,,
\end{equation}
and \begin{equation}\label{eqimageofdeltai2} 
\phi(\delta_{i,j})=\beta_{i,p+1}^{\beta_{p+1}\cdots \beta_{p+j-1}},\; i=1,\ldots,p,\; j=2,\ldots,q \,,
\end{equation} where the generators $\beta_i,\, i=1,\ldots, p+q-1,$
are as in Figure $3$ and  
 $$\beta_{i,j}=(\beta_i^2)^{\beta_{i+1}^{-1}\cdots \beta_{j-1}^{-1}}\,.$$
\end{prop}
\proof
 Using
\eqref{differenceofvu} and the methods of 
\cite{de00}, it is easy to see that
$$\phi(\delta_{i,1}) =\beta_{i,p+1}$$
and 
$$\phi(\delta_{i,j}) =(\beta_{i+j-1,p+j})^{(\beta_{j-1}^{-1}\cdots\, \beta_{p+j-2}^{-1})
\,\cdots\, (\beta_{2}^{-1}\cdots \beta_{p+1}^{-1})\,(\beta_1^{-1}\cdots \beta_p^{-1})}\,,\quad j=2,\ldots,q\,.$$
Using  a suitable homotopy argument in $\OO_{p+q}$ 
(deform the paths with initial points $y_0-y_1,\ldots,y_0-y_{p+j-1}$
to paths with constant real part and large enough imaginary part), 
one can see that for $j\geq 2$ these braids coincide 
with $\beta_{i,p+1}^{\beta_{p+1}\cdots \beta_{p+j-1}}.$
\Endproof

Using the choice of our setup, one obtains a diagram
\begin{equation}\label{splitdiag2}\begin{CD}
1\To\pi_1(U_0,(x_0,y_0)) @>>> \pi_1(U,(x_0,y_0))
    @>>>  \pi_1(S,y_0) \To 1 \\
   @VVV   @VVV @V{\phi}VV       \\
  1\To \pi_1(U_0,(x_0,y_0))@>>>\A_{p+q,1} @>>>\quad\quad\quad\A_{p+q}\To 1\, ,    \\
\end{CD}\end{equation}
such that the rows are split exact sequences 
and such that the vertical arrows are compatible with the
splittings of the rows (see \eqref{rowsplitting}). 

\begin{prop}\label{mondeltaij} The monodromy of $\V_1\ast \V_2$ is
given by
$$ \rho_{\V_1\ast \V_2}(\gamma)=\bar{\Phi}(T_{\V_1
\circ\V_2|_{U_0}},\phi(\gamma))\quad \forall \gamma \in \pi_1(S,y_0)\, ,$$
where $T_{\V_1 \circ\V_2|_{U_0}}=(C_1,\ldots,C_{p+q+1})$ is as in
\eqref{Tcirc} and $\bar{\Phi}$ is as in Section 
\ref{secvarpar}.
\end{prop}

\proof By the properties of the tensor 
product, the elements
$$C_i=A_i\otimes 1_{V_2} \quad {\rm and}\quad 
C_{p+j}=1_{V_1}\otimes B_j,\quad  \;
i=1,\ldots,p,\,j=1,\ldots ,q\, ,$$ commute. 
Using this, and using the braiding action of 
the elements $\phi(\delta_{i,1})$ 
(see \cite{dr03}, page 10), one can verify that  
 $$\rho_{\V_1\circ
\V_2}(\pi_1(S,y_0))=1\,.$$   By the above
discussion, Equation  \eqref{rowsplitting} can be assumed to hold for
\eqref{splitdiag2}. Thus Prop.~\ref{parmono} gives the claim.
\Endproof
\begin{rem}{\rm \begin{enumerate}
\item 
The author has used the last proposition 
in order to write  a computer program
in the computer algebra language GAP which computes the associated tuple 
$T_{\V_1\ast \V_2}$ of the middle convolution 
$\V_1\ast \V_2$ in the generic case from the associated tuples 
 $T_{\V_1}$ and $T_{\V_2},$ see \cite{DettweilerConv}.
\item In the non-generic case, a deformation argument shows that 
one obtains the associated tuple  $T_{\V_1\ast \V_2}$ from the generic case,
by multiplying the monodromy generators of the generic case 
suitably.
\end{enumerate}}
\end{rem}

\subsection{Irreducibility of the middle convolution and 
Kummer sheaves.}\label{secirrkumm}

The middle convolution of two irreducible local systems is in general 
not an irreducible local system.  In this section we give some 
irreducibility criteria for the middle convolution. 

\begin{defn}\label{defconvosheaf}{\rm  A local system $\V \in \LS_K(U_0)$ is a {\em
convolution sheaf} if it has no factors or quotients which are
isomorphic to the restriction $\V_\chi|_{U_0},$ where $\chi:
\pi_1(\AA^1\setminus\{x_1\})\to K^\times$ is a character of
$\pi_1(\AA^1\setminus\{x_1\})$ for some $x_1\in \AA^1$. The
category of convolution sheaves on $U_0$ is denoted by
$\Conv_K(U_0).$}
\end{defn}

\begin{rem} A local system $\V\in \LS_K(U_0)$ is a
convolution sheaf if and only if $T_\V=(T_1,\ldots,T_{r+1})$
satisfies the following two conditions, see 
\cite{dr00}, Section 3:\begin{enumerate}

\item[($\ast$)]  Let  $\V \in \LS_K(U_0)$ and $T_\V=(T_1,\ldots,T_{r+1})\in
\GL(V)^{r+1}.$  Then
$$ \bigcap_{j\not= i }{\rm ker}(T_j-1) \cap {\rm ker} (\tau T_i -1) =0,\; i=1,
\ldots , r , \;\forall \tau \in K^\times \,.$$

\item[($\ast\ast$)] Let
 ${ W}_i(\tau):=
 \sum_{j \neq i} {\rm im}(T_j-1) +{\rm im} (\tau T_i-1),
\; i=1,\ldots,r, \;\tau \in K^\times. $ Then 
$$ {\rm dim}({W}_i(\tau))=\dim (V) ,\;
i=1,\ldots,r, \;\forall \tau \in K^\times\,. $$ \end{enumerate}
\end{rem}

Throughout this section we assume that $\V_1$ is an
irreducible convolution sheaf with 
$T_{\V_1}=(A_1,\ldots, A_{p+1})$
and that $\V_2$ is a non-trivial rank one system 
 with $T_{\V_2}=(\lambda_1,\ldots ,\lambda_{q+1})$
and $\lambda_i\not=1,\, i=1,\ldots,q.$ \\

Set $\rho:=\rho_{\V_1\circ \V_2|_{U_0}}.$
Recall from Section~\ref{secparabcohom11} that 
$$H^1(U_0,\V_1\circ \V_2|_{U_0})\simeq H^1(\pi_1(U_0),V_1\otimes V_2)\,.$$
Here, 
$$H^1(\pi_1(U_0),V_1\otimes V_2)=C^1(\pi_1(U_0),V_1\otimes V_2)/B^1(\pi_1(U_0),V_1\otimes V_2)\,,$$ where 
\begin{eqnarray}\label{eqcoc}
C^1(\pi_1(U_0),V_1\otimes V_2)&:=&
\{(\delta:\pi_1(U_0)\to V_1\otimes V_2)\mid\nonumber \\
&&\quad  \delta(\alpha \beta)=
\delta(\alpha) \rho(\beta)+\delta(\beta)
\quad \forall\; \alpha,\beta\in \pi_1(U_0)\}\end{eqnarray}
is the group of {\em $1$-cocycles} and 
$$B^1(\pi_1(U_0),V_1\otimes V_2):=\{ \delta_v\mid v \in V_1\otimes V_2,\quad 
\delta_v(\gamma)=v(1-\rho(\gamma))\quad \forall\; \gamma\in \pi_1(U_0)\}$$
is the group of {\em $1$-coboundaries}.\\

In order to deal with monodromy questions, it is often convenient
 to use a base consisting of ``Pochhammer cycles'' (see 
\cite{dr03}):

\begin{defn}{\rm   Let $G$ be a group.
Then we define the {\em commutator of} $\alpha,\beta \in G$ as  
 $$[\alpha,\beta]:=\alpha^{-1}\beta^{-1}\alpha\beta\,.$$

 The linear map
$$
 \tau :C^1(\pi_1(U_0),V_1\otimes V_2)\To (V_1\otimes V_2)^{pq}\, ,$$
$$\delta \Mapsto 
\left(\,\delta([\alpha_1,\alpha_{p+1}]),\cdots,\,
\delta([\alpha_p,\alpha_{p+1}]),\ldots,\,
\delta([\alpha_1,\alpha_{p+q}]),\cdots,
\,\delta([\alpha_p,\alpha_{p+q}])\,\right)$$
is called the {\em twisted evaluation map.}}\end{defn}

\begin{lem}\label{3.4}  The kernel of the twisted evaluation map 
coincides with the coboundaries 
 $B^1(\pi_1(U_0),V_1\otimes V_2)$ and thus induces a map
$$ \tau :H^1(\pi_1(U_0),V_1\otimes V_2) \To (V_1\otimes V_2)^{pq}\,.$$
\end{lem}

 \proof Assume that $\delta\in \ker(\tau).$
The cocycle relation \eqref{eqcoc} implies
$$\delta([\alpha_i,\alpha_{p+j}])= \delta(\alpha_i)(\lambda_j-1)+
\delta(\alpha_{p+j})(1-A_i),\quad i=1,\ldots,p,\,j=1,\ldots,q \,.$$ 
Since $\lambda_j$ is assumed to be $\not=1$ and since 
 $\delta\in \ker(\tau),$
\begin{equation}\label{eqcotz} \delta(\alpha_i)=\frac{1}{1-\lambda_j}\delta(\alpha_{p+j})(1-A_i),
\quad i=1,\ldots,p,\,j=1,\ldots,q\,.\end{equation}
Set $v_j:=\frac{1}{1-\lambda_j}\delta(\alpha_{p+j}).$ By \eqref{eqcotz},
$$ v_j(1-A_i)=v_{j'}(1-A_i)\quad {\rm for} \quad 
i=1,\ldots,p,\quad {\rm and}\quad j,j'=1,\ldots,q \,.$$
If $v_j\not=v_{j'}$ then the above equality shows that 
the vector $v_j-v_{j'}$ spans a trivial 
$\langle A_1,\cdots, A_p\rangle$-submodule of $V_1.$ Since $\V_1$ was assumed
to be a convolution sheaf, this is impossible. Thus 
$$\frac{1}{1-\lambda_j}\delta(\alpha_{p+j})=
\frac{1}{1-\lambda_{j'}}\delta(\alpha_{p+j'})$$
so, by  \eqref{eqcotz}, 
$$\delta(\alpha_i)=\frac{1}{1-\lambda_q}\delta(\alpha_{p+q})(1-\rho(\alpha_i))\,.$$ Thus 
$\delta$ is a coboundary. Now the claim follows from dimension
reasons together with the fact that every vector $v\in V_1\simeq V_1\otimes 
V_2$ ($V_2$ was supposed to be one-dimensional) 
appears as some $\delta_v(\alpha_{p+q}).$
\Endproof 

Consider the projection 
$$p_1: (V^p)^q\to V^p,\, (v_1,\ldots,v_q)\mapsto v_1\,\,.$$
Via $p_1,$ the space 
$V^p$ turns into a $\langle \delta_{1,1},\ldots,\delta_{p,1}\rangle$-module:
The action is  induced
by the action of $\pi_1(S)=\langle \delta_{i,j}\rangle$ on 
$\im(\tau)$ (which in turn is induced by sending 
$\delta(\alpha)$ to $\delta(\alpha^{ \delta_{i,j}^{-1}}),$ see 
\cite{dw03}, Lemma 2.2). 

Let 
$$\tilde{D}_{i,1},\quad i=1,\ldots,p\, ,$$ be the linear 
transformation on $V^p$ induced by $\delta_{i,1}.$ These matrices 
 coincide with the Pochhammer matrices considered in \cite{dr00} and 
\cite{dr03}:

\begin{lem}\label{lempochmat} The linear transformation 
$\tilde{D}_{i,1}$ is of the following form:
\[ \left( \begin{array}{ccccccccc}
                  1 & 0 &  & \ldots& & 0\\
                   & \ddots &  & & &\\
                    & & 1 &&&\\
               \lambda_1 (A_1-1) & \ldots& \lambda_1 (A_{i-1}-1)  & \lambda_1 A_{i} & (A_{i+1}-1) & \ldots 
&   (A_p-1) \\
     &&&&1&&\\
               &   &  & && \ddots  &   \\             
                   0 &  &  & \ldots& &0 & 1
          \end{array} \right)
,\] where $\tilde{D}_{i,1}$ is the identity block matrix outside the 
$i$-th block row.
\end{lem}

\proof The cocycle relation \eqref{eqcoc} implies that
\begin{equation}\label{eqcoomraus}
 \delta(\alpha[\beta,\gamma]\epsilon)=\delta(\alpha\epsilon)+
\delta([\beta,\gamma])\rho(\epsilon)\end{equation}
and $\delta(\alpha^{-1})=-\delta(\alpha)\rho(\alpha)^{-1}.$
If $m<i,$ then these formulas together 
with Prop.~\ref{propdeltt} yield
\begin{eqnarray}\delta[\alpha_m,\alpha_{p+1}]\tilde{D}_{i,1}&=&
\delta[\alpha_m^{(\delta{i,1})^{-1}},
\alpha_{p+1}^{(\delta{i,1})^{-1}}]\nonumber \\
&=& \delta[\alpha_m,\alpha_{p+1}^{\alpha_i\alpha_{p+1}}]\nonumber \\
&=&\delta(\alpha_m^{-1}\alpha_{p+1}^{-1}[\alpha_i,\alpha_{p+1}]\alpha_m
[\alpha_{p+1},\alpha_i] \alpha_{p+1})\nonumber \\
&=&\delta[\alpha_i,\alpha_{p+1}]\lambda_1(A_m-1)+\delta[\alpha_m,\alpha_{p+1}]\,.
\nonumber \end{eqnarray}
If $m=i$ then 
\begin{eqnarray}\delta[\alpha_i,\alpha_{p+1}]\tilde{D}_{i,1}&=&
 \delta[\alpha_i^{\alpha_{p+1}},
\alpha_{p+1}^{\alpha_i\alpha_{p+1}}]\nonumber \\
&=&\delta(\alpha_{p+1}^{-1}\alpha_i^{-1}[\alpha_i,\alpha_{p+1}]\alpha_i\alpha_{p+1})\nonumber \\
&=&\delta[\alpha_i,\alpha_{p+1}]A_i\lambda_1\nonumber\,, 
\end{eqnarray}
where in the last equation we have used \eqref{eqcoomraus} and  
$\delta(1)=0.$
If $m>i,$ then 
\begin{eqnarray}\delta[\alpha_m,\alpha_{p+1}]\tilde{D}_{i,1}&=&
 \delta[\alpha_m^{[\alpha_i,\alpha_{p+1}]},\alpha_{p+1}^{\alpha_i\alpha_{p+1}}]\nonumber \\
&=&\delta[\alpha_i,\alpha_{p+1}](A_m-1)+\delta[\alpha_m,\alpha_{p+1}]\,
\nonumber. \end{eqnarray}
This proves the claim.\Endproof

Remember that the parabolic cohomology 
$H^1_p(U_0,\V_1\circ \V_2|_{U_0})=(\V_1\ast\V_2)_{y_0}$ was considered 
as the subspace of 
$H^1(\pi_1(U_0),V_1\otimes V_2)$ consisting of the elements $$[\delta]\in 
H^1(\pi_1(U_0),V_1\otimes V_2)$$ with 
$$\delta(\gamma)\in \im(\rho(\gamma)-1),\quad \forall \gamma \in \pi_1(U_0)\,.$$ 

\begin{prop}\label{propimparab} Let 
 $W$ be the image of the parabolic cohomology group 
$H^1_p(U_0,\V_1\circ \V_2|_{U_0})$ in $V^{p}$ 
under the composition of the twisted evaluation map and the projection
onto the first coordinate
$$ (V^{p})^{q}\To V^p,\, (v_1,\ldots,v_q)\mapsto v_1\,.$$ Then 
the following statements hold:
\begin{enumerate}
\item
The space $W$ is equal to ${\cal K}\cap {\cal L},$ where 
$$ {\cal K}=\{\;(w_1.\ldots,w_p)\mid w_i\in \im(A_i-1)\;\}$$ and 
$$ {\cal L}=\{\;(w_1A_2\cdots A_p,w_2A_3\cdots A_p,\ldots,w_p)\mid w_i\in \im(A_1\cdots A_p\lambda_1-1)\;\}\,.$$
\item The space $W$ is an irreducible 
$\langle \delta_{1,1},\ldots,\delta_{p,1}\rangle$-module.
\end{enumerate}
\end{prop}

\proof By the interpretation of $H^1_p(U_0,\V_1\circ \V_2|_{U_0})$
in terms of the cohomology with compact supports
(compare to Section~\ref{secparabcohom11} and 
\cite{dw03}), $W$ is 
contained in ${\cal K}.$ By the same reasoning, 
 $$\delta(\alpha_1\cdots\alpha_{p+1})\in 
\im(A_1\cdots A_p \cdot \lambda_1-1)\,,$$
thus 
one also has 
$W\subseteq {\cal L}.$ By looking at the definitions 
of ${\cal K}$ and ${\cal L}$ and at the structure of 
the matrices $\tilde{D}_{i,1},$ 
one immediately sees that ${\cal K}\cap {\cal L}$ is isomorphic to the dual
of the $\langle \delta_{1,1},\ldots,\delta_{p,1}\rangle$-module $MC_{\lambda_1}(V_1),$ where 
$MC_{\lambda_1}(V_1)$ is as in  \cite{dr03}. Thus 
${\cal K}\cap {\cal L}$ is an irreducible 
$\langle \delta_{1,1},\ldots,\delta_{p,1}\rangle$-module
by  \cite{dr03}, Thm. 2.4. It follows that $W$ coincides with 
${\cal K}\cap {\cal L}.$
\Endproof

\begin{rem}\label{remlab} By taking 
$q=1,$ the proof of the last corollary yields a new proof of the 
fact that 
$$MC_\lambda(T_{\V})=T_{\V\ast \V_\chi}\, \quad (\lambda \in K^\times)\,,$$
where $MC_\lambda$ is the tuple transformation of  
\cite{dr03} and $\V_\chi \in
\LS_K({\mathbb G}_m)$ is the  Kummer sheaf associated to 
$$\chi:\pi_1(\GG_m)\to K^\times,\, \gamma \mapsto \lambda$$
(where $\gamma$ denotes a generator of $\pi_1(\GG_m)$),
see Rem.~\ref{remmclambda} and \cite{dr03}.
\end{rem}

\begin{thm}\label{thmirrd} 
Let $K$ be a field. Let 
$\V_1\in \LS_K(U_1)$ be an irreducible convolution sheaf with 
$T_{\V_1}=(A_1,\ldots,A_{p+1})\in\GL_{n}(K)^{p+1}$ such that 
$$A_{p+1}=(A_1\cdots A_p)^{-1}=1\,.$$
Let 
 $\V_2\in \LS_K(U_2)$ be a rank one system with 
$T_{\V_2}=(\lambda_1,\ldots,\lambda_{q+1})$ such that  $\lambda_i\not=1$
for $i=1,\ldots,q.$
Assume that $\uo\ast \vo$ is generic.
Then the local system $\V_1\ast\V_2\in \LS_K(S)$ is irreducible if
$$ (p-2)n-
\sum_{i=1}^p\dim_K(\ker(A_i-1))>0\,.$$
\end{thm}

\proof This follows from  induction on $q.$ For $q=1$ this is 
Rem.~\ref{remlab} and 
\cite{dr03}, Thm. 2.4 (iii). If $q>1,$ then we can assume that 
${\V_1}\ast \tilde{\V}_2$ is irreducible, where 
$$ \tilde{\V}_2\in \LS_K(\AA^1\setminus \{y_2,\ldots,y_q\})\quad {\rm with}\quad T_{\tilde{\V}_2}=(\lambda_2,\ldots,\lambda_q,\lambda_1\cdot \lambda_{q+1})\,.$$
Let 
$$ p_1:(V^p)^q\To V^p ,\, (v_1,\ldots,v_q)\Mapsto v_1\,,$$
and 
$$p_2:(V^p)^q\To (V^p)^{q-1} ,\, (v_1,\ldots,v_q)\Mapsto (v_2,\ldots,v_q)\,.$$
Let $$G_1:=\langle {\delta_{i,1}},\, i=1,\ldots,p\rangle\,\leq \,\pi_1(S,y_0)$$ and 
$$G_2:=\langle {\delta_{i,j}},\, i=1,\ldots,p,\, j=2,\ldots,q\rangle\,\leq \,\pi_1(S,y_0)\,.$$

By Lemma \ref{3.4}, the twisted evaluation map 
$$\tau:C^1(\pi_1(U_0),V_1\otimes V_2)\To (V_1\otimes V_2)^{pq}$$
induces a well defined map of 
$$(\V_1\ast \V_2)_{y_0}=H^1_p(U_0,V_1\otimes V_2|_{U_0})$$
to $(V_1\otimes V_2)^{pq},$ which is also denoted by $\tau.$

 By Rem.~\ref{remlab}, the $G_1$-module
$\im(p_1\circ \tau)$ is isomorphic to the 
$G_1$-module 
$MC_{\lambda_1}(V_1)$ ($V_1$ denoting the stalk 
of $\V_1$) and is irreducible
by \cite{dr03}, Thm 2.4.  Since  $\pi_1(S,y_0)$ is the free product 
of $G_1$ and $G_2,$ it follows (by taking the $\pi_1(S,y_0)$-closure 
of the inverse image of $p_1\circ \tau$)
that 
$(\V_1\ast \V_2)_{y_0}$ contains an irreducible 
$\pi_1(S,y_0)$-submodule 
$W_1$ of
rank greater than or equal to 
$$\rk(MC_{\lambda_1}(V_1))=n_1:=pn-\sum_{i=1}^p \dim(\ker(A_i-1))\,.$$ 

By the induction hypothesis, the $G_2$-module
$$ (\V_1\ast \tilde{\V}_2)_{y_0}=
H^1_p(\AA^1\setminus (\uo\cup \{y_0-y_2,\ldots,y_0-y_q\}),
\V_1\circ\tilde{\V}_2)$$ is irreducible. 
Also, the map $$p_2\circ \tau: (\V_1\ast \V_2)_{y_0} \To 
V^{p(q-1)}$$ can be easily seen to be $G_2$-equivariant ($G_2$ acting 
via $\tau:(\V_1\ast \tilde{\V}_2)_{y_0} \to V^{p(q-1)}$) and 
non-trivial. Thus the $\pi_1(S,y_0)$-module
$(\V_1\ast \V_2)_{y_0}$ contains an irreducible 
submodule 
 $W_2$ of rank greater than or equal to
$$\rk(\V_1\ast \tilde{\V}_2)\geq n_2:=(p+q-3)n-
\sum_{i=1}^p\dim_K(\ker(A_i-1))\,.$$ 
The rank of $\V_1\ast \V_2$ is  smaller than or equal to
$$ n_3:=(p+q-1)n-\sum_{i=1}^{p}\dim_K(\ker(A_i-1))\,.$$
One has 
$$n_1+n_2=2pn +qn-3n-2\cdot \sum_{i=1}^p\dim_K(\ker(A_i-1))$$
Thus 
$$ n_1+n_2-n_3\geq (p-2)n-
\sum_{i=1}^p\dim_K(\ker(A_i-1))\,.$$ 
By assumption,  $(p-2)n-\sum_{i=1}^p\dim_K(\ker(A_i-1))>0,$ thus
$ n_1+n_2>n_3.$ It follows that the intersection of the irreducible submodules 
$W_1$ and $W_2$ is 
non-trivial and hence $\V_1\ast \V_2$ is irreducible.
 \Endproof

\subsection{The local monodromy of the middle convolution.}\label{secmonodrconv}

The local monodromy at the elements of $\uo\ast \vo$ in the
generic and semisimple case is given as follows:

\begin{lem} \label{lemmonodromy1} Let $\V_1\in 
\Conv_K(U_1)$
and $\V_2 \in
\LS_K(U_2)$ be irreducible local systems. Suppose that $\uo\ast \vo$ is
generic. Let $T_{\V_1}=(A_1,\ldots,A_{p+1})\in \GL(V_1)^{p+1}$ and
$T_{\V_2}=(B_1,\ldots,B_{q+1})\in \GL(V_2)^{q+1}$ and suppose that
the elements $B_1,\ldots,B_{q}$ are semisimple. Let
\begin{multline}\nonumber 
T_{\V_1\ast\V_2}=(D_{i,j})\in \GL(V_1\ast V_2)^{pq+1},\\
  D_{i,j}:=\rho_{\V_1\ast \V_2}(\delta_{i,j}),\;
i=1,\ldots,p,\;j=1,\ldots,q\end{multline} 
be the associated tuple 
of $\V_1\ast \V_2,$ where the $\delta_{i,j}$ are the generators of
$\pi_1(S)$ as above and the tuple is
ordered according to \eqref{eqproductdeltai}.  Then the following
holds:

Every non-trivial eigenvalue $\beta$ of $B_j$ and every Jordan
block $J(\alpha,l)\not= J(1,1)$ occurring in the Jordan
decomposition of $A_i$ contribute a Jordan block $J(\alpha
\beta,l')$ to the Jordan decomposition of $D_{i,j},$ where
$$ l':\;=\quad
  \begin{cases}
    \quad l &
                              \quad\text{\rm if $\alpha \not= 1,\beta^{-1}$,} \\
    \quad  l-1& \quad \text{\rm if $\alpha =1$,} \\
    \quad l+1 & \quad \text{\rm if $\alpha =\beta^{-1}$.}
  \end{cases}
  $$
  The only other Jordan blocks which occur in the Jordan
  decomposition of $D_{i,j}$ are blocks of the form $J(1,1).$
\end{lem}

\proof Using  a deformation argument,
one can assume that $i=p$ and $j=1,$ as well as $x_p=y_1=0.$ Let
$$S^o:=\{y\in S\mid |y|<\epsilon\}\, ,$$ where
$\epsilon$ is chosen small enough so that no other element of
$\uo\ast\vo$ lies in $S^o.$ We can assume that $y_0\in S^o$ and
that the support of $\delta_{p,1}$ lies also in $S^o.$ 
For any eigenvalue $\beta_m$ of $B_1,$ let 
$\V_{\beta_m}$ be the Kummer sheaf associated 
to 
$$ \pi_1(\GG_m(\CC),y_0)\To K^\times,\, \gamma \Mapsto \beta_m\,.$$ Let 
$$\tilde{\V}\in
\LS_K(U_0\cup \{y_0-y_1\})$$ be the local
system associated to the tuple
$$(A_1,\ldots,A_p,B_2,\ldots,B_{q},(A_1\cdots A_p\cdot B_2\cdots B_q )^{-1})\,.$$
Using the above cocycle-calculus and Prop.~\ref{propexactconv}, 
one can see
 that
there is an isomorphism of $\langle \delta_{p,1}\rangle$-modules
$$ 
\V_1\ast \V_2|_{S^o}\simeq \bigoplus_{\beta_m\not=
1}(\tilde{\V}\ast \V_{\beta_m})|_{S^o}\, ,$$ 
where the sum is over the non-trivial
$\beta_m$ which are counted with multiplicity. 
Thus the claim is
\cite{dettweiler04}, Lemma 3
(see also  \cite{dr00}, Lemma 4.1, Rem.
\ref{remmclambda}  and  Rem.~\ref{remlab}).
\Endproof

The following lemma gives the monodromy at infinity in the  case
of the convolution with Kummer sheaves:

\begin{lem}\label{lemmonodromy2} Let 
$\V\in \Conv_K(U_1)$ be a convolution sheaf and let 
$\V_\chi$ be a Kummer sheaf associated to a non-trivial 
character 
$$\chi: \pi_1(\GG_m(\CC))\to K^\times,\, \gamma \mapsto 
\lambda\,.$$  Let $T_{\V}=(A_1,\ldots,A_{p+1})\in \GL(V)^{p+1.}$
 Let
$$T_{\V\ast\V_\chi}=(D_1,\ldots,D_{p+1})\in \GL(V\ast V_\chi)^{p+1},\quad
D_{i}:=\rho_{\V\ast \V_\chi}(\delta_{i,1}),
$$ be the associated tuple 
of $\V\ast\V_\chi.$  Then the
following holds:

Every Jordan
block $J(\alpha,l)$ occurring in the Jordan
decomposition of $A_{p+1}$ contributes a Jordan block $J(\alpha
\lambda^{-1},l')$ to the Jordan decomposition of $D_{p+1},$ where
$$ l':\;=\quad
  \begin{cases}
    \quad l, &
                              \quad\text{\rm if $\alpha \not= 1,\lambda$,} \\
    \quad  l-1& \quad \text{\rm if $\alpha =\lambda$,} \\
    \quad l+1 & \quad \text{\rm if $\alpha =1$.}
  \end{cases}
  $$
  The only other Jordan blocks which occur in the Jordan
  decomposition of $D_{p+1}$ are blocks of the form $J(\lambda^{-1},1).$
\end{lem}

\proof The claim is
\cite{dettweiler04}, Lemma 3.\Endproof

\section{Convolution of \'etale local systems} \label{secetaleconN}

\subsection{Variation of parabolic cohomology in the \'etale
case.}\label{secvarparcohometale}

In this section we recall the main results of \cite{dw03} on the
variation of \'etale parabolic cohomology. Throughout Section 
\ref{secetaleconN}, our coefficient ring $R$ is 
 a topological ring as in the Sections~\ref{secgeneralnotation}
 and~\ref{secetaledef}.\\

 Let $S$ be a smooth, affine and
geometrically irreducible variety over  a field $k\subseteq \CC,$  
 let $\D\subseteq \PP^1_S$ be a smooth relative divisor of degree $r+1$ over
$S$ which contains the section $\{\infty\}\times S.$ Let $U:=\PP^1_S
\setminus \D,$ let $j: U \to \PP^1_S$ be the natural inclusion,
and let $\bar{\pi}: \PP^1_S
\to S$ be the projection onto $S.$ Let 
 $\pi: U \to S$ be the restriction of
$\bar{\pi}$ to $U,$  let $\bar{s}_0$ denote a geometric point of
$S,$ and let $U_0=U_{0,\overline{k}}$ denote 
the fibre over $\bar{s}_0.$

\begin{defn}{\rm  \label{varconvdef}
  A {\em variation} of $\V_0\in \LS_R^\et(U_0)$ over $S$ is an \'etale
  local system $\V$ on $U$ whose restriction to $U_0$ is equal
  to $\V_0$. The {\em parabolic cohomology} of the variation $\V$ is
  the sheaf of $R$-modules on $S$
  \[
        \W \;:=\; R^1\pib_*(j_*\V).
  \]}
\end{defn}

 See loc.~cit., Thm 3.2, for the next result:

\begin{prop}\label{propdw0333} \begin{enumerate}
\item The parabolic cohomology $\W$ is an \'etale local system of
$R$-modules on $S.$
\item There is a natural isomorphism of local systems on $S$
$$ \W\an \;\liso\; R^1\pib_*(j_*\V\an)\,.$$ In particular, the
fibre of $\,\W$ at  $\bar{s}_0$ may be identified with
$H^1_p(U_0\an,\V_0\an).$\end{enumerate}
\end{prop}

Poincar\'e duality implies the following result (see 
\cite{dw03}, \cite{dw04} and \cite{MilneEC}):

\begin{prop}  \label{proppoincareetale}{\rm (Poincar\'e Duality)}
Let $U$ be as above, let
  $\V\in \LS^\et_R(U),$ and let
$$\V\otimes\V(n)\to{R}$$ be a non-degenerate 
bilinear pairing of sheaves ($\V(n)$ denoting 
the $n$-th Tate twist as in Section~\ref{sectatetwi}), corresponding to an isomorphism
   $\V(n)\inj \V^*.$ 
Then the cup product
  defines a non-degenerate  bilinear
pairing of sheaves $$\W \otimes \W(n+1) \to  {R}\, ,$$
where $\W=R^1\pib_*(j_*\V)$ is as above.
\end{prop}

\subsection{The middle convolution for \'etale local
systems.}\label{secvondefetale} We use the notation of 
Section~\ref{secgeomprell}.\\

Let $k$ be a subfield of $\CC.$
 Let $\uo
:=\{ x_1,\ldots,x_p\} \in \OO_p(k),$ let $\vo :=\{y_1,\ldots,y_q\} \in
\OO_q(k),$ and let 
$$ \uo \ast \vo := \{ x_i + y_j \mid i=1,\ldots,p, \,\, j = 1,\ldots
,q \}$$ (compare to Section~\ref{sectconvdef}). 
Let $U_1:=\AA^1_k\setminus \uo,$
$U_2:=\AA^1_k\setminus \vo$ and $S:=\AA^1_k\setminus \uo \ast \vo$
(as affine varieties over $k$).  Define the polynomial $f(x,y)$ as
in Section~\ref{sectconvdef} and let $$U=\AA^2_k\setminus
\{f(x,y)=0\}:={\rm Spec}(k[x,y,\frac{1}{f(x,y)}]) \in \Var_k\,.$$ Define
$$\pr_1: U \To U_1,\, (x,y) \Mapsto x\, ,$$
$$ \pr_2: U \To S,\; (x,y) \Mapsto y$$ and
$$ \q: U \To U_2,\, (x,y) \Mapsto y-x$$
(as morphisms of varieties).
Let $j:U\to \PP^1_S$ be the natural 
inclusion.  The second projection $\PP^1_S=\PP^1_k\times S \to S$ is
denoted by $\opr_2.$ Let $y_0\in S(k)$ be a $k$-rational point of 
$S$ and let $\bar{y}_0$ be a geometric point of $S$ extending
$y_0.$ 
The fibre over $\bar{y}_0$ is denoted by
$U_0=U_{0,\bar{k}}.$ \\

Let $\V_1\in \LS_R^\et(U_1),$ let 
 $\V_2 \in \LS_R^\et(U_2),$ and let 
$$\V_1\tim \V_2:=\pr_1^*(\V_1)\otimes \q^*(\V_2)\,.$$
As in the context of local systems, the 
middle convolution of \'etale local systems is defined as the 
\'etale parabolic cohomology of 
the variation $\V_1\tim \V_2:$

\begin{defn} \begin{enumerate}{\rm 
\item  The {\em middle convolution} of $\V_1 \in \LS_R^\et(U_1)$
and $\V_2 \in \LS_R^\et(U_2)$ is the \'etale local system
$$ \V_1 \ast \V_2 := R^1 (\opr_2)_*(j_*(\V_1\tim \V_2)) \in\LS_R^\et(S)\,.$$
\item The {\em $\ast$-convolution} of $\V_1 \in \LS_R^\et(U_1)$ and
$\V_2 \in \LS_R^\et(U_2)$ is the \'etale local system
$$ \V_1 \ast_\ast \V_2:=R^1 (\pr_2)_*(\V_1\tim \V_2) \in\LS_R^\et(S)\,.$$
\item The {\em $!$-convolution} of $\V_1 \in \LS_R^\et(U_1)$ and
$\V_2 \in \LS_R^\et(U_2)$ is the \'etale local system
$$ \V_1 \ast_! \V_2:=R^1 (\pr_2)_!(\V_1\tim \V_2) \in\LS_R^\et(S)\, ,$$
where $R^1 (\pr_2)_!$ denotes higher direct image with compact
support, see \cite{FreitagKiehl}, Chap. I.8.}
\end{enumerate}
\end{defn}

\subsection{First properties of the middle convolution in the \'etale
case.}\label{secfirstpropetalecase}

Choose a $k$-rational point $(x_0,y_0)\in U(k)$ and choose 
an extension  $(\bar{x}_0,\bar{y}_0)\in U(\bar{k}).$ \\

The following properties are the \'etale analogues of 
Prop.~\ref{dimensione} (using 
Prop.~\ref{propdw0333} and 
Remark~\ref{remanaletal}) and  Prop.~\ref{propconviscommutative}:

\begin{prop}\label{dimension} Let  $\V_1 \in \LS_R^\et(U_1)$ and $\V_2
\in \LS_R^\et(U_2).$ Let $V_1:=(\V_1)_{\bar{x}_0}$ and 
let $V_2:=(\V_2)_{\bar{y}_0-\bar{x}_0}.$
Let $T_{\V_1}=(A_1,\ldots,A_{p+1})\in \GL(V_1)^{p+1}$ and
$T_{\V_2}=(B_1,\ldots,B_{q+1})\in \GL(V_2)^{q+1}$ be the associated 
monodromy tuples
as in Section~\ref{secetaledef}. 
Then
\begin{enumerate}
\item
Suppose that $R=K$ is a field and that one of the stabilizers
$$V_1^{\pi_1^\et(U_1,\bar{x}_0)}\quad \quad  and \quad \quad 
V_2^{\pi_1^\et(U_2,\,\bar{y}_0-\bar{x}_0)}$$ is trivial. 
Let $\dim_K V_i=n_i.$
Then
\begin{multline}\label{dimeeee} \rk(\V_1\ast \V_2)= (p+q-1)n_1 n_2 -
\sum_{i=1}^p n_2\dim_K
\ker(A_i - 1_{V_1}) \\
-\sum_{j=1}^q n_1\dim_K \ker(B_j - 1_{V_2})
-\dim_K\ker(A_{p+1}\otimes B_{q+1}-1_{V_1\otimes
V_2})\,.\end{multline}

\item \begin{equation*}\label{convcommv} \V_1\ast \V_2 \cong \V_2\ast
\V_1 \,.\end{equation*}
\end{enumerate}
\end{prop}

The following result is the \'etale analogue of 
 Lemma~\ref{naive}: 

\begin{lem}\label{naive2} Let $\V_1 \in \LS_R^\et(U_1)$
and  $\V_2 \in \LS_R^\et(U_2).$ Then the following holds:
\begin{enumerate}
\item
The local systems $R^i{\pr}_{2*}(\V_1\circ \V_2)$ vanish for
$i\not= 1$ if the stabilizer $V_2^{\pi_1^\et(U_2)}$ is
trivial. 

\item {\rm (Exactness of the $\ast$-convolution)} Let 
$V_2^{\pi_1^\et(U_2)}= 0.$
Then the functor
 $$\LS_R(U_1) \To \LS_R(S),\quad 
\V_1 \Mapsto \V_1\ast_\ast \V_2\,,$$
is exact.
\end{enumerate}
\end{lem}

\proof This follows analogously to Lemma 
\ref{naive}, using \cite{FreitagKiehl}, Thm. I.9.1. \Endproof

\begin{rem} \label{remconnectiontokatzconv} 
 Let $\V_2 \in \LS_R^\et(U_2)$  such that its associated 
tuple 
$T_{\V_2}$ generates an absolutely irreducible 
and non-trivial subgroup of $\GL(V_2).$
As in Prop.~\ref{propexactconv} one can show that
$$R^0(\overline{\pr}_2)_*(j_*(\V_1\circ \V_2))=0=R^2(\overline{\pr}_2)_*(j_*(\V_1\circ
\V_2))\quad \forall\, \V_1\in \LS_R^\et(U_1)\, ,$$ so that the functor
$$ \LS_R(U_1) \To \LS^\et_R(S),\; \V_1 \Mapsto \V_1 \ast \V_2\,,$$
is exact.
\end{rem}

\begin{defn}\label{defpuree}{\rm Let $k$ be a number field, 
let $X$ be a smooth and
geometrically irreducible variety over $k,$ and let $\F \in
\Constr_R(X)$ be a constructible sheaf 
(see Section~\ref{secetaledef}).  
Consider  a geometric point $\bar{x}\in X(\overline{k})$
which is defined over a finite extension $k'$ of $k.$
 One obtains a
Galois representation $ \rho_{\bar{x}} : G_{k'} \to \GL(\V_{\bar{x}}).$ 
Let $E_\lambda$ be a complete subfield 
of $\bar{\QQ}_\ell.$

 \begin{enumerate}
\item A constructible 
sheaf  $\V\in \Constr_{E_\lambda}(X)$ 
is called {\em punctually pure of weight $w\in \QQ$} if for any
finite extension $k'$ of $k$ and any
geometric point  $\bar{x}\in X(\bar{k})$ which is defined over 
$k',$ the
$G_{k'}$-module  $\V_{\bar{x}}$ is  pure of weight $w$ (see Def.~\ref{defeinss}).

\item Consider a
 local system $\V \in \LS^\et_{{E_\lambda}}(X)$ and $w_1,w_2\in
\ZZ$ such that $w_1\leq w_2.$ Then $\V$  is called {\em mixed of
weights $[w_1,w_2]$} if there exists a filtration
$$ \V=\V^{w_2}\supset \V^{w_2-1} \supset \ldots \supset
\V^{w_1}\supset \V^{w_1-1}:=0$$ such that $\V^i/\V^{i-1}$ is punctually pure
of weight $i$ for $i=w_1,\ldots,w_2$ and if $\V^{w_1}\not=0.$ The
quotient $\V^i/\V^{i-1}$ is called the {\em weight-$i$-part} of
$\V$ and is denoted by $W^i(\V).$ \end{enumerate}}
\end{defn}

\begin{prop}\label{propweilzwei} Suppose that 
$U_1$ and $U_2$ are defined over a number field $k.$ 
Let $\V_1\in \LS_{E_\lambda}^\et(U_1)$ 
be a sheaf which is punctually pure of  weight $w_1$ and let 
$\V_2\in \LS_{E_\lambda}^\et(U_2)$ be punctually pure 
of weight
 $w_2.$ Then the following statements hold:
\begin{enumerate}
\item
The local system $\V_1\ast \V_2\in \LS_{E_\lambda}^\et(S)$ is
punctually pure of weight $w_1+w_2+1.$
\item If one of the sheaves is 
irreducible and non-trivial, then there is
 an exact sequence of sheaves on $S:$
$$ 0 \To {\cal K}\To \V_1\ast_!\V_2 \To \V_1\ast\V_2 \To 0\, ,$$
where the weight of ${\cal K}$ is $\leq w_1+w_2.$
\end{enumerate}
\end{prop}

\proof Clearly, the local system $\V_1\circ \V_2$ is punctually pure of
weight $w_1+w_2$ on $U.$  Let $k'$ be a finite  extension of $k$ 
and let $z_0$ be a geometric point of $S$ which is defined over $k'.$ 
The stalk $(\V_1\ast
\V_2)_{\bar{z}_0}$ is, as a $G_{k'}$-module, isomorphic to
$H^1(\PP^1_{\bar{k}},j_*(\V_1\circ \V_2)|_{\PP^1_{\bar{k}}}).$ The claim (i) follows
thus from the base change theorem (see \cite{FreitagKiehl} or \cite{MilneEC})
and 
Deligne's work on the Weil conjectures
 (Weil II) \cite{DeligneWeil2}.

The claim (ii) follows verbatim as in \cite{Katz96}, proof of Lemma 8.3.2:
There is an exact sequence of sheaves on $\PP^1_S$
(where the subscript ${}_!$ stands for extension by zero in the 
sense of \cite{MilneEC}):
$$ 0\To j_!\left(\pr_1^*(\V_1)\otimes \pr_2^*(\V_2)\right)\To 
j_*\left(\pr_1^*(\V_1)\otimes \pr_2^*(\V_2)\right)\To Q\To 0\, ,$$
and the sheaf $Q$ is the direct sum of sheaves $Q_s,\, s=1,\ldots,p+q+1,$
 concentrated along the irreducible components $d_i,\, i=1,\ldots,p+q+1,$ of the relative divisor 
$$(\PP^1\times S)\setminus U\,.$$
Without loss, we can assume that $\V_2$ is irreducible and 
non-trivial. 
Then the long exact cohomology sequence of the above short 
exact sequence is
$$ 0\To R^0(\overline{\pr}_{2})_*(Q)\To \V_1\ast_! \V_2 \To \V_1\ast \V_2 \To 
0\, ,$$ since (by the irreducibility and non-triviality of 
$\V_2$) $$R^0(\overline{\pr}_{2})_*
(j_*(\pr_1^*(\V_1)\otimes \pr_2^*(\V_2)))=0$$ and since 
$R^2(\overline{\pr}_{2})_*
(Q)=0$ ($Q$ is concentrated at the divisors $d_i$). The sheaf 
$R^0(\overline{\pr}_{2})_*(Q)=(\overline{\pr}_{2})_*(Q)$
 coincides with the direct sum of the 
sheaves $Q_i$ viewed as a sheaf on the base. Thus the weight of 
$(\overline{\pr}_{2})_*(Q)=R^0(\overline{\pr}_{2})_*(Q)$ is 
$\leq w_1+w_2.$ 
 \Endproof

By Poincar\'e duality (see Prop.~\ref{proppoincareetale}) and since 
the comparison isomorphism between singular cohomology 
and \'etale cohomology is compatible with the Poincar\'e pairing 
(see \cite{DMOS}),  
one obtains the following 
result:

\begin{prop}\label{proppoin} Let $\V_1 \in \LS_R^\et(U_1)$ and $\V_2 \in \LS_R^\et(U_2).$ Let
$\kappa_1: \V_1 \otimes \V_1(n_1) \to {R}$ and
$\kappa_2:\V_2 \otimes \V_2(n_2) \to {R}$ be non-degenerate
bilinear
pairings of \'etale local systems.
Then there  is a non-degenerate bilinear pairing
\begin{equation}\label{bili} \kappa_1 \ast \kappa_2:
(\V_1\ast \V_2)\otimes(\V_1\ast \V_2)(n_1+n_2+1) \to
{R}\,. \end{equation} Moreover, if $\kappa_i\an,\,i=1,2,$ denotes 
the induced pairing on $\V_i\an$ and if $(\kappa_1\ast \kappa_2)\an$
denotes the induced pairing on  $(\V_1\ast \V_2)\an$ (see Prop. 
\ref{fform}), then 
$$ \kappa_1\an\ast \kappa_2\an=(\kappa_1\ast \kappa_2)\an=\kappa_1\ast \kappa_2\,.$$
\end{prop}

\begin{rem} \label{rempoinabs} 
Here are some  remarks on Poincar\'e duality and 
weights which are useful in later applications:
Let ${\bar{s}}$ be a geometric point of $S$
which is defined over $k.$ Then the following statements hold:
\begin{enumerate}
\item 
If in the pairing of Prop.~\ref{proppoin}, the number 
$n_1+n_2+1$ is an even number equal to $2m,$ then the action of 
$\pi_1^\et(S)$ on  the stalk 
of the $m$-fold Tate twist
$((\V_1\ast \V_2)(m))_{\bar{s}}$ respects the bilinear form 
given by Poincar\'e
duality (see Prop.~\ref{fform}), since \eqref{bili} can be written as
\begin{equation} 
(\V_1\ast \V_2)(m)\otimes(\V_1\ast \V_2)(m) \to
{R}\,.\nonumber \end{equation}
\item  
Assume that the $G_k$-module 
 $ (\V_1\ast \V_2)(m)_{\bar{s}}$ is unramified at $\pi\in \PP^f(k).$ 
If $\V_1\ast\V_2$ is punctually pure of weight $2m,$ then 
 the eigenvalues 
of the Frobenius elements $\Frob_\pi$ 
acting on $(\V_1\ast \V_2)(m)_{\bar{s}}$ have absolute value equal to $1.$
\end{enumerate}
\end{rem}

\section{Geometric interpretation of the middle
convolution}\label{secmotiv}

In  Sections~\ref{secoo}--\ref{secgeommi} we provide a geometric interpretation 
of sequences of middle convolutions of \'etale sheaves with finite 
monodromy. This is similar to \cite{Katz96}, 
Chap. 8. The main difference is that, on the one hand, our base 
rings are less general than in loc.~cit., but on the other hand,
the convolutions considered here 
are more general (in loc.~cit., Katz considers only the case 
of convolutions with Kummer sheaves). As a corollary, in Section 
\ref{secconse}, one obtains 
a bound on the occurring determinants, see Thm.~\ref{thmdet}.

\subsection{The underlying fibre spaces.}\label{secoo}

Let $\ka$ be a subfield of $\CC$ and 
 fix sets $$\uo_i \in \OO_{p_i}(\ka), \quad i=1,\ldots
,n\,.$$  For $1\leq j\leq l\leq n,$ define the space
\begin{multline*} \OO_{[j,l]}(\CC):=\{ (x_j,\ldots,x_l)\in
\CC^{l-j+1}\quad \mid\quad  x_i \notin \uo_1\ast\uo_{2}\ast \cdots \ast
\uo_i\,,\;  i=j,\ldots ,l \\
{\rm and} \quad x_{i+1}-x_i\notin \uo_{i+1}\,,\;
i=j,\ldots, l-1,\quad {\rm if}\quad l>j  \}. 
\end{multline*} Let $\OO_{[j,l]}$
be the underlying variety over $k.$ 
For $1\leq j_1\leq j_2\leq l_2 \leq l_1 \leq n,$ there is the
projection map
$$ \pr^{[j_1,l_1]}_{[j_2,l_2]}:\OO_{[j_1,l_1]}\To
\OO_{[j_2,l_2]},\;(x_{j_1},\ldots,x_{{l_1}})\Mapsto
(x_{j_2},\ldots,x_{{l_2}})\,.$$ 
A special case of such a projection is the map 
 $$ \pr^{[j_1,l_1]}_{[l_2,l_2]}:\OO_{[j_1,l_1]}\To
\OO_{[l_2,l_2]}=\AA^1_k\setminus \uo_1\ast \cdots \ast 
\uo_{l_2},\;(x_{j_1},\ldots,x_{{l_1}})\Mapsto x_{l_2}\,.$$ 
Moreover, for $j \leq i <i+1
\leq  l$ there are  the difference maps
$$ \q^{[j,l]}_{[i,i+1]}: \OO_{[j,l]} \To \AA^1_k\setminus
\uo_i,\, (x_j,\ldots,x_\ell)\Mapsto x_{i+1}-x_i\,.$$

\begin{prop}\label{propcart1} 
The following  commutative diagram is cartesian in $\Var_k:$
\begin{equation}\label{eqschneider1}\begin{CD} \OO_{[1,n]}
@>\pr^{[1,n]}_{[1,n-1]}>>\OO_{[1,n-1]} \\
   @V\pr^{[1,n]}_{[n-1,n]}VV   @VV\pr^{[1,n-1]}_{[n-1,n-1]}V \\
  \OO_{[n-1,n]} @>>\pr^{[n-1,n]}_{[n-1,n-1]}>\OO_{[n-1,n-1]}
\end{CD}\end{equation}
\end{prop}
\proof This follows from the fact that the 
equations 
of $\OO_{[1,n-1]}$ and $\OO_{[n-1,n]}$ 
add up to  the equations
of $\OO_{[1,n]}.$ \Endproof

\subsection{Geometric interpretation of the $\ast$-convolution.}\label{secast}

It is the aim of this section to provide a geometric interpretation 
of iterative sequences of tensor products and the $\ast$-convolutions. 
This is used in the next section to describe the corresponding 
sequences of tensor products of middle convolutions.\\

We proceed using the definitions and assumptions of the last section. 
Throughout this and the following section, we use the following 
 assumptions:

\begin{ass}\label{ass1}
{\rm Let $$f_i: \bF_i \to \AA^1_\ka  \setminus \uo_i,\, i=1,\ldots,n\, ,$$
be  finite \'etale Galois covers with Galois groups $F_i\not=\{1\}$
and let 
$$g_j:\bG_j
\to \AA^1_\ka  \setminus \uo_1\ast \cdots \ast \uo_{j+1},\,
j=1,\ldots,n-1\, ,$$ 
be finite \'etale Galois covers with Galois groups $G_j\not=\{1\}.$
Let $E$ be a number field and let
 $$\chi_i:F_i\to \GL_{m_i}(E),\, i=1,\ldots,n\, ,$$ and 
 $$\xi_j:G_j\to \GL_{n_j}(E),\, j=1,\ldots,n-1\, ,$$
be faithful 
 and 
absolutely irreducible representations, respectively. }
\end{ass}

There exist idempotent  elements $\PPP_{\chi_i}\in E[F_i]$
such that $\PPP_{\chi_i}\in E[F_i]$ applied to the 
regular representation of $E[F_i]$ is isomorphic to  
$\chi_i.$ Similarly, there exist  idempotent  elements
$\PPP_{\xi_j}\in E[G_j]$ such that $\PPP_{\xi_j}$ applied to the 
regular representation of $E[G_j]$ is isomorphic to  $\xi_j.$  In the sequel, fix a set of such 
projectors
(the choice $\PPP_{\chi_i}$ and $\PPP_{\xi_j}$ 
 is not canonical in general).

Let $\lambda$ be a finite prime of $E$ with $\char(\lambda)=\ell$ and 
let $E_\lambda$ denote the completion of $E$ with respect to $\lambda.$
Using the notation of Section~\ref{seccovv}, let 
$${\cal F}_i:=\LL_{(f_i,\chi_i)} \in
\LS_{E_\lambda}^\et(\AA^1_\ka  \setminus \uo_i)\,,\, i=1,\ldots,n\, ,$$ 
resp.
$${\cal G}_j:=\LL_{(g_j,\xi_j)} \in
\LS_{E_\lambda}^\et(\AA^1_k\setminus \uo_1\ast \cdots \ast
\uo_{j+1})\,,\,
j=1,\ldots,n-1$$
(where $\GL_n(E),\,n\in \NN,$ is viewed as a subgroup of $\GL_n(E_\lambda)$ in the 
obvious way).\\

Let 
$$ \pi_n^\times: \bF_1 \boxtimes 
\cdots \boxtimes
\bF_{n} \boxtimes \bG_1 \boxtimes \cdots \boxtimes \bG_{n-1} \To \OO[1,n]$$ denote the direct product over $\OO[1,n]$
of the  pullback
of the cover $f_1$ along $\pr^{[1,n]}_{[1,1]},$
 of the pullbacks 
 of the covers 
$f_i,\, i=2,\ldots,n$ along  the difference maps $\q^{[1,n]}_{[i-1,i]}$
and of the 
pullbacks of the covers $g_j,\,j=1,\ldots,n-1,$ 
along the projections 
$\pr^{[1,n]}_{[j+1,j+1]}.$\\

By construction, the map $\pi_n^\times$ is an  \'etale Galois cover
  with Galois
group $G^\times$ isomorphic to
$F_1\times \cdots \times F_{n}\times G_1\times\cdots \times G_{n-1}$
(that $G^\times$ is the direct product can easily be seen by embedding
a general line into $\OO[1,n]$).
Since 
\begin{eqnarray}
{G}&:=&F_1\otimes \cdots \otimes F_{n}\otimes G_1\otimes\cdots \otimes 
G_{n-1} \nonumber \\
&:=& \im(\chi_1\otimes \cdots \otimes \chi_{n}\otimes 
\xi_1\otimes \cdots \otimes \xi_{n-1}
)\nonumber \end{eqnarray}
is a factor of $G^\times$ (it is a central product of 
the groups $F_1,\ldots ,F_n$ and $G_1,\ldots,G_{n-1},$ compare to 
Section~\ref{secgrouptheoryprel}), one obtains 
an \'etale Galois cover (as a quotient of $\pi_n^\times$)
$$ \pi_n: {\frak U}:=\bF_1 \otimes \cdots \otimes
\bF_{n}\otimes \bG_1 \otimes \cdots \otimes \bG_{n-1}  \To \OO[1,n]$$   with Galois
group isomorphic to $G.$

 Let $$\Delta_j:=\chi_1\otimes \cdots \otimes \chi_j\otimes 
\xi_1 \otimes \cdots \otimes \xi_{j-1}\,, \quad j=2,\ldots,n\,
$$ and let
 $\PPP_{\Delta_j} \in E[\im(\Delta_j)]$  be the    idempotent
which cuts out $\Delta_j$ from the regular 
representation of $\im(\Delta_j).$
Moreover, for $j=2,\ldots,n,$ let 
\begin{eqnarray}\nonumber 
 \KKK_j&:=&{\cal F}_1\circ {\cal F}_2\circ\GGG_1\circ  {\cal F}_3 \circ \GGG_2
\circ \cdots \circ {\cal F}_j \circ \GGG_{j-1}\\
&:=&(\pr^{[1,j]}_{[1,1]})^*(
{\cal F}_1)\otimes (\q^{[1,j]}_{[1,2]})^*({\cal F}_2)\otimes
 (\pr^{[1,j]}_{[2,2]})^*(\GGG_1)\otimes  (\q^{[1,j]}_{[2,3]})^*({\cal F}_3) \nonumber\\
&&\otimes  (\pr^{[1,j]}_{[3,3]})^*(\GGG_2)
\otimes \cdots \otimes (\q^{[1,j]}_{[j-1,j]})^*({\cal F}_j) \otimes
(\pr^{[1,j]}_{[j,j]})^*(\GGG_{j-1})\,.\nonumber 
\end{eqnarray}

The next result gives the geometric interpretation of 
 sequences of 
$*$-convolutions and tensor products of the above 
local systems. 
This is the key result 
for later applications. 

\begin{thm} \label{thmmot1} 
Let $n\geq 2,$ let
 $$S:=\OO_{[n,n]}=\AA^1_k\setminus \uo_1\ast \cdots \ast \uo_n\,,$$ and let 
$\Pi:= \pr^{[1,n]}_{[n,n]}\circ
\pi_n:{\frak U} \to S.$
  Let 
$ \V_1:={\cal F}_1\in 
\LS_{E_\lambda}^\et(\OO_{[1,1]})$
and for 
$j=2,\ldots,n,$ define sheaves $ \V_j$
by setting 
$$ \V_j:=(\V_{j-1}\ast_\ast {\cal F}_j)\otimes \GGG_{j-1}\in \LS_{E_\lambda}^\et(\OO_{[j,j]})\,.$$ 
Then
$$ \V_n=\PPP_{\Delta_n} \left(R^{n-1}
\Pi_*({E_\lambda})\right)\,.$$
Here,   $\PPP_{\Delta_n}$ acts as a projector on $R^{n-1}
\Pi_*({E_\lambda})$ via the embedding 
$$E[G]\leq \End({\frak U}/S)\otimes E$$
(and hence cuts out a sub-local system of $R^{n-1}
\Pi_*({E_\lambda})$). 
\end{thm}

\proof We use induction on $n.$ Let $n=2:$ By definition, 
$R^0{\pi_2}_*({E_\lambda})={\pi_2}_*({E_\lambda}).$
 Since $\pi_2$ is a Galois covering,
the local system  $R^0{\pi_2}_*({E_\lambda})$ corresponds to 
the regular representation ${\operatorname{reg}}$ of $G_1\otimes F_1\otimes F_2,$ i.e.,
$$R^0{\pi_2}_*({E_\lambda})=\LL_{(\pi_2,{\operatorname{reg}})}\,,$$
in the notation of Section~\ref{seccovv}.
Thus
\begin{equation}\label{33}{\cal F}_1\circ {\cal F}_2\circ \GGG_1=
\PPP_{  \chi_1 \otimes
\chi_2\otimes\xi_1}R^0{\pi_2}_*({E_\lambda}),\end{equation} 
where $\PPP_{
\chi_1 \otimes
\chi_2\otimes\xi_1 }$ acts as an element of $\End( \bF_1\otimes \bF_2\otimes \bG_1)\otimes E$ 
on $R^0{\pi_2}_*({E_\lambda}).$ 
It follows that
\begin{eqnarray}
\V_2=({\cal F}_1\ast_\ast {\cal F}_2)\otimes \GGG_1&=&R^1{\pr^{[1,2]}_{[2,2]}}_*({\cal F}_1\circ{\cal F}_2)\otimes \GGG_1 \label{eq34}\\
&=& R^1{\pr^{[1,2]}_{[2,2]}}_*({\cal F}_1\circ{\cal F}_2\circ \GGG_1)
\label{35}\\
&=&R^1{\pr^{[1,2]}_{[2,2]}}_*(\PPP_{\chi_1 \otimes
\chi_2 \otimes\xi_1 }(R^0{\pi_2}_*({E_\lambda})))\label{36}\\
&=& \PPP_{ \chi_1 \otimes
\chi_2\otimes \xi_1 }\left(R^1{\pr^{[1,2]}_{[2,2]}}_*(R^0{\pi_2}_*({E_\lambda}))\right)\label{37}\\
&=& \PPP_{ \chi_1 \otimes
\chi_2\otimes \xi_1}\left(R^1(\pr^{[1,2]}_{[2,2]}\circ
\pi_2)_*({E_\lambda})\right),\label{38}
\end{eqnarray} where Equation 
\eqref{eq34} holds by definition, \eqref{35} follows from  the projection
formula, \eqref{36} follows from
\eqref{33}, \eqref{37} follows from the fact that taking higher
direct images commutes with automorphisms (smooth base change, see 
\cite{FreitagKiehl}, Thm. 7.3) and
\eqref{38} is an application of the Leray spectral sequence. This proves the claim for $n=2.$
 
Assume that the claim is true for $n-1\geq 2.$ 
By a similar argument as above, 
\begin{equation}\label{00schneid} \KKK_j=\PPP_{\Delta_j}(\pi_{j*}({E_\lambda})),\,\, j=2,\ldots,n\,.\end{equation}
  By the induction hypothesis, 
\begin{eqnarray}\label{39} \V_{n-1}&=&
R^{n-2}(\pr^{[1,n-1]}_{[n-1,n-1]})_*( \KKK_{n-1})\label{40000} \\
&=& \PPP_{\Delta_{n-1}}\left(
R^{n-2}(\pr^{[1,n-1]}_{[n-1,n-1]}\circ\pi_{n-1})_*({E_\lambda})\right).
\label{4000}
\end{eqnarray}
 Then
\begin{eqnarray}
 \V_n&=& 
R^1{\pr^{[n-1,n]}_{[n,n]\ast}}(  \pr^{[n-1,n]*}_{[n-1,n-1]}(\V_{n-1})\otimes\q^{[n-1,n]*}_{[n-1,n]} ({\cal F}_n))\otimes\GGG_{n-1}\label{400} \\
 &=&
R^1{\pr^{[n-1,n]}_{[n,n]_*}}( \pr^{[n-1,n]*}_{[n-1,n-1]}
R^{n-2}\pr^{[1,n-1]}_{[n-1,n-1]*}(\KKK_{n-1})\otimes
\label{4004} \\
&&\quad \quad 
 \q^{[n-1,n]*}_{[n-1,n]} ({\cal F}_n)\otimes\pr^{[n-1,n]*}_{[n,n]}(\GGG_{n-1})) \nonumber \\
&=&
R^1{\pr^{[n-1,n]}_{[n,n]_*}}( R^{n-2}\pr^{[1,n]}_{[n-1,n]*}(
\pr^{[1,n]*}_{[1,n-1]}(\KKK_{n-1}))\otimes
\label{404}\quad \quad  \\
&&\quad \quad 
\q^{[n-1,n]*}_{[n-1,n]} ({\cal F}_n) \otimes \pr^{[n-1,n]*}_{[n,n]}(\GGG_{n-1})) \nonumber
\end{eqnarray}
\begin{eqnarray}
\quad \quad 
&=&
R^1\pr^{[n-1,n]}_{[n,n]*}(
R^{n-2}\pr^{[1,n]}_{[n-1,n]*}(\KKK_{n})) \label{44} \\
&=&
R^{n-1}\pr^{[1,n]}_{[n,n]*}(\KKK_{n}) \label{41} \\
&=& R^{n-1}\pr^{[1,n]}_{[n,n]*}(\PPP_{\Delta_n}(R^0\pi_{n*}({E_\lambda}))) \label{45} \\
&=& \PPP_{\Delta_n} \left(R^{n-1}\pr^{[1,n]}_{[n,n]*}(R^0\pi_{n*}({E_\lambda})\right)\label{46}\\
&=& \PPP_{\Delta_n} \left(R^{n-1}(\pr^{[1,n]}_{[n,n]}\circ\pi_{n})_*({E_\lambda})\label{47}\right),
\end{eqnarray}
where the above formulas \eqref{400} -- \eqref{47} are derived 
using  the following arguments: 

Equation \eqref{400} holds by definition and
\eqref{4004} holds by the projection formula and \eqref{40000}. 
 The diagram \eqref{eqschneider1}
is cartesian. Thus, by smooth base change,  
$$\pr^{[n-1,n]*}_{[n-1,n-1]}\left(R^{n-2}\pr^{[1,n-1]}_{[n-1,n-1]*}
(\KKK_{n-1})\right)=R^{n-2}\pr^{[1,n]}_{[n-1,n]*}\left(\pr^{[1,n]*}_{[1,n-1]}(\KKK_{n-1})\right)\,.$$ This yields 
\eqref{404}. Equation \eqref{44} follows from the K\"unneth-formula,
\eqref{41} follows from the Leray spectral sequence,
\eqref{45} follows from \eqref{00schneid}, 
\eqref{46}  follows from smooth base change and \eqref{47} is 
again an application of the 
Leray spectral sequence. 
\Endproof

\subsection{Geometric interpretation of the middle convolution.}\label{secgeommi}

Let \begin{multline}\nonumber 
 \V_n=(\cdots ((({\cal F}_1\ast_\ast {\cal F}_2)\otimes\GGG_1)\ast_\ast{\cal F}_3)\otimes 
\cdots \ast_\ast {\cal F}_n)\otimes \GGG_{n-1}\in \\
\LS_{E_\lambda}^\et(\AA^1_k\setminus \uo_1\ast \cdots \ast
\uo_n)\end{multline} be as in the last section. It follows from
Deligne's work on the Weil conjectures
(\cite{DeligneWeil2}), that $\V_n$ is mixed of weights $[n-1,2(n-1)].$
The following result gives an interpretation 
of the corresponding sequence of {\em middle} convolutions 
$\V^n$ in terms of the weight filtration:

\begin{thm}\label{thmpureweight} Let $n\geq 2.$ Let $\F_i,\,i=1,\ldots,n,$ 
and $\G_i,\, i=1,\ldots,n-1,$  be irreducible and nontrivial 
local systems with finite monodromy as in the last section, and let
\begin{multline}\nonumber 
 \V_j=(\cdots ((({\cal F}_1\ast_\ast {\cal F}_2)\otimes\GGG_1)\ast_\ast{\cal F}_3)\otimes 
\cdots \ast_\ast {\cal F}_j)\otimes \GGG_{j-1}\in \\\LS_{E_\lambda}^\et(\AA^1_k\setminus \uo_1\ast \cdots \ast
\uo_j),\,\, j=2,\ldots,n\, ,\end{multline}
be
as in Thm.~\ref{thmmot1}.  Let
$$ \V^j:=(\cdots ((({\cal F}_1\ast {\cal F}_2)\otimes\GGG_1)\ast{\cal F}_3)\otimes 
\cdots \ast {\cal F}_j)\otimes \GGG_{j-1}\in  \LS_{E_\lambda}^\et(\AA^1_k\setminus \uo_1\ast \cdots \ast
\uo_j)\;$$ be the corresponding sequence of middle convolutions. 
 Then 
 $$ \V^{n}= W^{n-1}(\V_n)\, ,$$ where $W^{n-1}$ denotes the weight-$(n-1)$-part of $\V_n$ in the sense of Def.~\ref{defpuree}.
\end{thm}

\proof  Since the weight of $\GGG_{n-1}$ is zero, 
we can assume that $\GGG_{n-1}$ is trivial. 
By Poincar\'e duality, the statement of the theorem 
is  equivalent to saying that 
$${\V}^{n}= W^{n-1}\left({\V}_n^!\right)=W^{n-1}\left({(\cdots (({\cal F}_1\ast_! {\cal F}_2)\otimes\GGG_1)\ast_!{\cal F}_3)\otimes 
\cdots )\ast_! {\cal F}_n} \right)\,.$$
For $n=2,$ the statement follows from Prop.~\ref{propweilzwei}. Thus we can assume that 
$n\geq 3.$
Let us assume that 
$$ {\V}^{n-1}=(\cdots ((({\cal F}_1\ast {\cal F}_2)\otimes\GGG_1)\ast{\cal F}_3)\otimes 
\cdots \ast {\cal F}_{n-1})\otimes \GGG_{n-2}\;$$
coincides with 
$W^{n-1}({\V}_{n-1}^{!}),$ where 
$${\V}_{n-1}^{!}=(\cdots ((({\cal F}_1\ast_! {\cal F}_2)\otimes\GGG_1)\ast_!{\cal F}_3)\otimes 
\cdots \ast_! {\cal F}_{n-1})\otimes \GGG_{n-2}\,.$$
Thus
one obtains a short 
 exact sequence 
$$ 0\to K \to {\V}_{n-1}^! \to {\V}^{n-1} \to 0\, ,$$ 
where  $K$ is mixed of weights 
$\leq n-3$ and ${\V}^{n-1}$ is punctually pure of 
weight $n-2.$ Since 
the $!$-convolution with $\F_n$ is an exact functor 
(the dual of Lemma~\ref{naive2}) and by the description 
of the middle convolution in terms of the cohomology with compact supports
(see \cite{dw03}), the composite  $\pi$ of the 
following maps
$$  {\V}_n^!={\V}_{n-1}^!\ast_!\F_n \to{\V}^{n-1}\ast_! \F_n\to {\V}^{n}={\V}^{n-1}\ast\F_n$$
is surjective. By Prop.~\ref{propweilzwei} (i),
 the sheaf  $\V^{n-1}\ast \F_n$ is punctually pure of weight 
$n-1.$ Moreover, the  kernel 
of $\pi$ is mixed of weights $\leq n-2,$ which 
 follows again from the exactness of the $!$-convolution and from 
~\ref{propweilzwei} (ii).  
\Endproof

Let us recall the results of Bierstone and Milman, see
\cite{BierstoneMillman}, Thm. 13.2 as well of the results 
of Encinas, Nobile, and Villamayor \cite{EV1}, \cite{ENV}:

\begin{thm}\label{thmbiermil}
Let $k$ be a field of characteristic zero and let $X\in \Var_k.$ Then there
exists a morphism $\sigma_X: \tilde{X}\to X$ such that:\begin{enumerate}
\item $\tilde{X}$ is smooth over $k.$

\item Let ${\rm sing}(X)$ denote singular locus of
$X$ and let $E:=\sigma_X^{-1}({\rm sing}(  X)).$  Then $E$ is a
normal crossings divisor, i.e., $E$ is the union of smooth 
hypersurfaces $E_1,\ldots,E_l$ which have only normal crossings. 

\item Let $X,Y\in \Var_k$ and let $\phi:X|_U\to Y|_V$ be an
isomorphism over (nonempty) open subsets $U\subseteq X$ and
$V\subseteq Y.$ Then there exists an isomorphism $\phi':
\tilde{X}|_{\sigma_X^{-1}(U)}\to \tilde{Y}|_{\sigma_Y^{-1}(V)}$ such that the
following diagram commutes:
\begin{equation*}\label{cartesian6}\begin{CD}
\tilde{X}|_{\sigma_X^{-1}(U)} @>{\phi'}>>  \tilde{Y}|_{\sigma_Y^{-1}(V)} \\
   @VVV   @VVV     \\
U  @>{\phi}>> V\,.    \\
\end{CD}\end{equation*}
\end{enumerate}
\end{thm}

\begin{rem} Desingularization of varieties over fields of 
characteristic zero is quite well
understood nowadays: There are the algorithmic versions 
of the above mentioned papers 
\cite{BierstoneMillman}, \cite{EV1}, \cite{ENV}. One even has a 
computer implementation by Bodn\'ar and Schicho \cite{Schicho}. 
\end{rem}

Let $
\Pi:{\frak U}
 \to
S$ be as in Thm.~\ref{thmmot1}.
Let $\bar{s}_0\in S(\bar{k})$ be a geometric base point of 
$S$ which is defined over $k.$ One can find a 
$G$-equivariant embedding
${\frak U}_{\bar{s}_0}\to \CX_{\bar{s}_0},$ 
where $G$ is as in the last section, 
and where $\CX_{\bar{s}_0}$ is a projective variety over $\bar{s}_0$
which is defined over $k.$   
By the above result,
there exists a morphism which is defined over $k$
$$\sigma_{\CX_{\bar{s}_0}}:\tilde{\CX}_{\bar{s}_0}\to \CX_{\bar{s}_0}$$  such that
$\tilde{\CX}_{\bar{s}_0}$ is smooth projective  and such that 
$$D_{\bar{s}_0}:=\tilde{\CX}_{\bar{s}_0}\setminus \sigma_{\CX_{\bar{s}_0}}^{-1}({\frak U}_{\bar{s}_0})$$ is a normal
crossings divisor.
By Thm.~\ref{thmbiermil} (iii), 
we can assume that the action of $G$ on ${\frak U}_{\bar{s}_0}$ carries over
to an action on $\tilde{\CX}_{\bar{s}_0}.$ 
Thus the projector $\PPP_{\Delta_n}$ extends to a
projector 
${\PPP}_{\Delta_n}\in \End(\tilde{\CX}_{\bar{s}_0})\otimes E$ 
(which is denoted 
by the 
same symbol).  

\begin{cor}\label{corisogalmod} Let $\ka$ be a number field and let 
$ \V^{n}\in \LS_{E_\lambda}^\et(S)$
be as in Thm.~\ref{thmpureweight}. Let 
${\frak U},\, \CX_{\bar{s}_0},$ and $\tilde{\CX}_{\bar{s}_0}$  be as above.
Then there exists  an isomorphism 
of $G_k$-modules
$$  (\V^{n})_{{\bar{s}_0}}\simeq \im\left(\PPP_{\Delta_n}(H^{n-1}(\tilde{\CX}_{{\bar{s}_0}},E_\lambda))\To 
\PPP_{\Delta_n}(H^{n-1}({\frak U}_{{\bar{s}_0}},E_\lambda))\,\right)\,.$$
\end{cor}

\proof 
Since $\tilde{\CX}_{{\bar{s}_0}}$ is smooth projective and 
${\frak U}_{{\bar{s}_0}}$ can be seen as a dense open subset of $\tilde{\CX}_{{\bar{s}_0}},$ 
it follows from 
 \cite{Katz96}, 9.4.3,  that 
  $$\im\left(H^{n-1}(\tilde{\CX}_{{\bar{s}_0}},E_\lambda)\to 
H^{n-1}({\frak U}_{{\bar{s}_0}},E_\lambda)\,\right)\simeq W^{n-1}(H^{n-1}({\frak U}_{{\bar{s}_0}},E_\lambda))\,.$$
By smooth base change,
\begin{multline}\nonumber 
\im\left(\PPP_{\Delta_n}(H^{n-1}(\tilde{\CX}_{{\bar{s}_0}},E_\lambda))\to 
\PPP_{\Delta_n}(H^{n-1}({\frak U}_{{\bar{s}_0}},E_\lambda))\right)\simeq \\
W^{n-1}\left(\PPP_{\Delta_n}(
H^{n-1}({\frak U}_{{\bar{s}_0}},E_\lambda))\right)\,.\end{multline}
But
$$W^{n-1}\left(\PPP_{\Delta_n}(
H^{n-1}({\frak U}_{{\bar{s}_0}},E_\lambda))\right)
\simeq W^{n-1}((\V_{n})_{{\bar{s}_0}})\nonumber 
\simeq (\V^n)_{{\bar{s}_0}}\;,$$
where the first isomorphism follows from base change and 
Thm.~\ref{thmmot1} and 
the second isomorphism follows from Thm.~\ref{thmpureweight}. 
 \Endproof

\subsection{Consequences of the geometric interpretation 
of the middle convolution.}\label{secconse}

In this section we derive some consequences of the results of the last 
section. Our main result 
is a description of the occurring 
determinants, see 
Thm.~\ref{thmdet}. The concept of a compatible system of Galois representations 
will play a crucial role:

\begin{defn} {\rm Let $k,E$ be  number fields.
A {\em strictly compatible system  of ($n$-dimensional, 
$\lambda$-adic)  
$E$-rational Galois representations
 of ${ G}_k$} consists of a
 collection $\{ \rho_\lambda: { G}_k\to \GL_n(E_\lambda)\}_{
\lambda \in \PP^f(E)}$
of Galois representations and  a finite set $Z \subseteq \PP^f(k)$ of primes of $k$ (the {\em exceptional set})
 such that the following holds:
 \begin{itemize}

\item For any prime $\lambda \in \PP^f(E)$ and any
prime $\pi \in \PP^f(k)\setminus 
 Z$ whose characteristic is different from 
 $\chara(\lambda),$
the Galois representation $\rho_\lambda$
is unramified at $\pi$ and the coefficients of the 
characteristic polynomial $${\rm det}(1-\rho_\lambda({\rm Frob}_\pi)\cdot t)$$
are contained 
in $E.$

\item  For any two primes  $\lambda_1,\lambda_2\in \PP^f(E)$ and for any 
prime
$\pi \in \PP^f(k)\setminus Z$ whose characteristic
is different from  $\char(\lambda_1)$  and from $\char(\lambda_2)$   
there is an equality of 
characteristic polynomials
$$  {\rm det}(1-\rho_{\lambda_1}({\rm Frob}_\pi)\cdot t) = 
 {\rm det}(1-\rho_{\lambda_2}({\rm Frob}_\pi)\cdot t)\,.$$

\end{itemize}}
\end{defn}

\begin{defn}{\rm  Let $X$ be a nonsingular projective variety
 over an algebraically closed field.
 Let $\Pic(X)\simeq H^1(X,\OO_X^*)$ 
be the Picard group of $X.$ 
The {\em Neron-Severi group} of $X$ is defined to be the 
group  $\NS(X):=\Pic(X)/\Pic^\circ(X).$ 
Let $\Pic^n(X)\leq \Pic(X)$ denote the 
group of divisors which are numerically equivalent to zero. 
We set $\N(X):=\Pic(X)/\Pic^n(X).$
}\end{defn}

\begin{rem}\label{remns}{\rm  Let $X$ be a nonsingular projective surface over 
an algebraically closed subfield  of $\CC.$  Then the following 
holds: \begin{enumerate}
\item The {\em Neron-Severi group} $\NS(X)$ is a finitely generated 
group.  Moreover, the group $\N(X)$ is a free abelian group and 
$\Pic^n(X)/\Pic^\circ(X)$ is finite, see \cite{MilneEC}, Lemma V.3.26 and 
Lemma V.3.27. The rank of $\N(X)$  is called 
the {\em Picard number} of $X.$
\item The Kummer sequence 
\begin{equation}\nonumber \Pic(X)\stackrel{\ell}{\To} \Pic(X) \To H^2(X,\mu_\ell)\end{equation}
 induces 
the class map $$ cl_X:C^1(X)\To \Pic(X)\stackrel{\alpha}{\To} H^2(X,\QQ_\ell(1))\, ,$$ 
see \cite{MilneEC}, Rem. VI.9.6.
The map $\alpha$ factors as $$\Pic(X)\To \N(X)\stackrel{\beta}{\To}H^2(X,\QQ_\ell(1))\, ,$$
where the map $\beta$ is injective, see \cite{MilneEC}, Section V.3. 
\item The cup product of $H^2(X,\QQ_\ell(1))$ restricted to $\N(X)$
is induced by the intersection product of algebraic cycles, see 
\cite{MilneEC}, Prop. VI.9.5. 
Moreover,
the intersection product on 
$\N(X)$ is non-degenerate (this follows from the definition  of numerical 
equivalence and 
\cite{MilneEC}, 
Lemma V.3.27). \end{enumerate}}
\end{rem}

\begin{prop}\label{propnsbeidsei}
Let $k$ be a number field, let $X_k$  be a smooth and 
geometrically irreducible 
surface over $k,$ and let 
 $X=X_{k}\otimes \bar{k}.$ 
Then the following holds:
\begin{enumerate}
\item 
The 
group $\NS(X)_\ell= \N(X)\otimes_\ZZ \QQ_\ell$ is
a $G_k$-submodule of  $H^2(X,\QQ_\ell(1)).$
\item   Let $E$ be a number field. For $\lambda\in \PP^f(E),$ let 
$$\NS(X)_\lambda:=\NS(X)\otimes_\ZZ E_\lambda=\N(X)\otimes_\ZZ E_\lambda$$
be the $\lambda$-adic Neron-Severi group. 
Let $G\leq \Aut_k(X)$ be a finite group of automorphisms of $X_k$ 
and let  $\PPP\in E[G]$ be an idempotent element.
Then $\PPP(\NS(X)_\lambda)$ is a $G_k$-submodule of $\NS(X)_\lambda$ and 
the system of Galois representations 
$$  \left(\nabla_\lambda: G_k\to \GL(\PPP(\NS(X)_\lambda)\,\right)_{\lambda\in \PP^f(E)}$$
is a strictly compatible system of $E$-rational Galois representations.  
\end{enumerate}
\end{prop}

\proof  
The group $\Aut(X)$ acts in the usual way on $H^2(X,\QQ_\ell(1))$ 
and via transport of structures
on $\Pic(X)$ which induces 
the structure of an $\Aut(X)$-module on $\NS(X)$ and thus on $\NS(X)\otimes \QQ_\ell.$ 
 Since the class map is compatible 
with the action of $G_k\leq \Aut(X)$ on $C^1(X)$ and on $H^2(X,\QQ_\ell(1)),$ see \cite{MilneEC}, Prop. VI.9.2,
the group $\NS(X)\otimes \QQ_\ell$ is a $G_k$-submodule of $H^2(X,\QQ_\ell(1)).$

By Rem.~\ref{remns} (ii), there is  a sequence 
of $G_k$-modules
$$ \Pic(X)\To \N(X)\To \NS(X)\otimes_\ZZ E \To \NS(X)_\lambda \,.$$ 
There exist finitely many divisor classes $d_1,\ldots,d_s$
which generate $\N(X)$ and which  are permuted by $G_k$ and by $G.$ 
 Let $f\in G_k$ be any element. Then 
$f$ commutes with $\PPP,$
 and the characteristic polynomial of $f\PPP$ acting on 
$\NS(X)\otimes_\ZZ E$ is an element in the polynomial ring $E[t].$ It follows that 
the characteristic polynomial of $f\PPP$ acting on $\NS(X)_\lambda$ is an element in the polynomial ring $E[t]$ which is independent 
of $\lambda.$ Assume that $X$ has 
good reduction at a finite prime $p$ of $k$ and that $\chara(\lambda)\not=
\chara(p).$ 
By the above arguments, the characteristic polynomial of $\Frob_p\PPP$ acting on 
$\NS(X)_\lambda$ is an element in the polynomial ring $E[t]$ which is independent 
of $\lambda.$ It follows that the characteristic polynomial of $\Frob_p$ on 
$\PPP(\NS(X)_\lambda)$ is an element of $E[t]$ which is 
independent 
of $\lambda.$
Since $\nabla_\lambda $ is unramified outside any $p$ with 
$\chara(\lambda)=
\chara(p)$ and outside 
the finitely many 
places of bad reduction of $X_k,$ the system $(\nabla_\lambda)_{\lambda\in \PP^f(E)}$
is a strictly compatible system of $E$-rational Galois representations.  
\Endproof

\begin{prop}\label{corisogalmod2} Let $k,E$ be number fields.
For $\lambda \in \PP^f(E),$ let
$$ \V=\V^3=(({\cal F}_1\ast {\cal F}_2)\otimes\GGG_1)\ast{\cal F}_3 \in \LS_{E_\lambda}^\et(S_k)\;$$
be as in Thm.~\ref{thmpureweight} with $n=3.$ 
 Assume that 
the associated 
tuple 
$T_{\V}$ generates an infinite and absolutely irreducible subgroup 
of $\GL(\V_{{\bar{s}_0}}),$ where $s_0 \in S(k).$
Let 
$$\rho_\lambda^{{\bar{s}_0}}:G_k\to \GL(\V_{{\bar{s}_0}})$$ be the 
Galois representation on the stalk $\V_{{\bar{s}_0}}.$ Then  
 the system of Galois representations 
$(\rho_\lambda^{{\bar{s}_0}})_{\lambda\in \PP^f(E)}$ 
is a strictly compatible system 
of $E$-rational Galois representations. \end{prop}

\proof  It suffices to prove that 
the Tate twisted system  $ (\rho_\lambda^{{\bar{s}_0}}(1))_{\lambda\in \PP^f(E)}$  
is a strictly compatible system 
of $E$-rational Galois representations.

By Cor.~\ref{corisogalmod}, 
$$  \V_{{\bar{s}_0}}(1)\simeq \im\left(\PPP_{\Delta_n}(H^2(\tilde{\CX}_{{\bar{s}_0}},E_\lambda(1)))\To 
\PPP_{\Delta_n}(H^2({\frak U}_{{\bar{s}_0}},E_\lambda(1)))\,\right)\,,$$
where $\tilde{\CX}_{{\bar{s}_0}}$ is 
a smooth projective surface and $\tilde{\CX}_{{\bar{s}_0}}\setminus {\frak U}_{{\bar{s}_0}}$ is 
the union of smooth divisors which intersect transversally.
By Rem.~\ref{remns}, the restriction 
of the Poincar\'e pairing on $H^2(\tilde{\CX}_{{\bar{s}_0}},E_\lambda(1))$
to the $\lambda$-adic 
Neron-Severi group 
$$\NS(\tilde{\CX}_{{\bar{s}_0}})_\lambda:=
\NS(\tilde{\CX}_{{\bar{s}_0}})\otimes E_\lambda \leq 
H^2(\tilde{\CX}_{{\bar{s}_0}},E_\lambda(1))$$  is non-degenerate. Also,
the Poincar\'e pairing on $H^2(\tilde{\CX}_{{\bar{s}_0}},E_\lambda(1))$
is  compatible with the $\pi_1^\et(S)$-action (this follows from Prop.~\ref{propweilzwei} and 
 Rem.~\ref{rempoinabs} (ii)).
 Thus the $\pi_1^\et(S,{\bar{s}_0})$-module $H^2(\tilde{\CX}_{{\bar{s}_0}},E_\lambda(1))$ decomposes into a direct sum $\NS(\tilde{\CX}_{{\bar{s}_0}})_\lambda
\oplus \NS(\tilde{\CX}_{{\bar{s}_0}})_\lambda^\perp,$ where 
$\NS(\tilde{\CX}_{{\bar{s}_0}})_\lambda^\perp$ denotes the
 orthogonal complement of $\NS(\tilde{\CX}_{{\bar{s}_0}})_\lambda.$

By the excision sequence, the kernel of the map 
$$H^2(\tilde{\CX}_{{\bar{s}_0}},E_\lambda(1))\To 
 H^2({\frak U}_{{\bar{s}_0}},E_\lambda(1))$$
is contained in the Neron-Severi group $\NS(\tilde{\CX}_{{\bar{s}_0}})_\lambda.$
(This follows from an  application of 
the long exact sequence of \cite{MilneEC}, Rem. VI.5.4(b), together with the fact that 
the class map $cl_{\tilde{\CX}_{{\bar{s}_0}}}$ can be defined using the 
Gysin map, see \cite{MilneEC}, Section VI.9.)
Since the decomposition $\NS(\tilde{\CX}_{{\bar{s}_0}})_\lambda
\oplus \NS(\tilde{\CX}_{{\bar{s}_0}})_\lambda^\perp$ is preserved 
by $G$ (this follows from 
\cite{MilneEC}, Prop. VI.9.2), it follows that  the  $\pi_1^\et(S,{\bar{s}_0})$-module  
$\PPP_{\Delta_n}(H^2(\tilde{\CX}_{{\bar{s}_0}},E_\lambda(1)))$ 
can be written as a direct sum 
$$ M^{{\bar{s}_0}}_{1,\lambda}\oplus M^{{\bar{s}_0}}_{2,\lambda}:=
\left(\PPP_{\Delta_n}({ 
\NS(\tilde{\CX}_{{\bar{s}_0}})_\lambda})\right) \oplus 
 \left(\PPP_{\Delta_n}({ 
\NS(\tilde{\CX}_{{\bar{s}_0}})_\lambda^\perp})\right)\,.$$

It follows from  the 
assumptions on $T_\V$ that 
$ M^{{\bar{s}_0}}_{1,\lambda}$ coincides with the kernel of the map 
$$\PPP_{\Delta_n}(H^2(\tilde{\CX}_{{\bar{s}_0}},E_\lambda))\To 
\PPP_{\Delta_n}(H^2({\frak U}_{{\bar{s}_0}},E_\lambda))\,.$$ Consequently, 
one obtains
an isomorphism of $\pi_1^\et(S,{\bar{s}_0})$-modules
$\V_{{\bar{s}_0}}\simeq M^{{\bar{s}_0}}_{2,\lambda}.$  
Let $\rho_\lambda:G_k\to \GL(M^{{\bar{s}_0}}_{2,\lambda})$ be the induced Galois
representation. It remains to show that the system 
$(\rho_\lambda)_{\lambda \in \PP^f(E)}$ is 
a strictly compatible system 
of $E$-rational Galois representations. This follows from the 
following arguments:  The projector $\PPP_{\Delta_n}$ can be written 
as a finite sum 
$$ \PPP_{\Delta_n}=\sum a_i g_i\, ,$$ where $g_i$ is an element of the group
$${G}=F_1\otimes F_2\otimes  F_3\otimes G_1
\leq \End(\tilde{\CX}_{{\bar{s}_0}})\otimes E$$
and $a_i\in E$ (here, $F_1\otimes F_3 \otimes  F_3\otimes G_1$ denotes 
the central product of the groups $F_1,F_2, F_3,G_1$ as in 
Section~\ref{secast}).  
As in \cite{KisinWortmann}, Example
to Def. 1.2, one can use the Lefschetz fixed point formula 
for the endomorphisms $\Frob_p\circ g_i,$ 
together with the Weil-conjectures,
in order to show that  the characteristic polynomial of 
$\Frob_p\circ \PPP_{\Delta_n}$ acting on 
$H^2(\tilde{\CX}_{{\bar{s}_0}},E_\lambda)$
 is $E$-rational (i.e., contained in $E[t]$)
and independent of $\lambda$ (whenever the expression 
$\Frob_p\circ \PPP_{\Delta_n}$ makes sense, i.e., $\tilde{\CX}_{{\bar{s}_0}}$
has good reduction at $p$ and $\chara(\lambda)\not=p$).
The last statement says that 
 the characteristic polynomial of $\Frob_p$ acting on 
$ M^{{\bar{s}_0}}_{1,\lambda}\oplus M^{{\bar{s}_0}}_{2,\lambda}$ is $E$-rational 
and independent of $\lambda.$
 
By Prop.~\ref{propnsbeidsei} (ii),  the characteristic polynomial
 of $\Frob_p$ on 
 $M^{{\bar{s}_0}}_{1,\lambda}$ is $E$-rational and 
independent of $\lambda.$
This and the independence of $\Frob_p$ on 
$ M^{{\bar{s}_0}}_{1,\lambda}\oplus M^{{\bar{s}_0}}_{2,\lambda}$ imply 
 that the  characteristic polynomial
 of $\Frob_p$ on 
 $M^{{\bar{s}_0}}_{2,\lambda}$
 is $E$-rational and 
independent of $\lambda.$ By the properties 
of the \'etale cohomology, 
$\rho_\lambda$ is unramified outside the (finitely many) places 
of bad reduction of $\tilde{\CX}_{{\bar{s}_0}}$ and outside the 
places having characteristic $\ell.$ Thus 
$(\rho_\lambda)_\lambda$ 
is indeed a strictly compatible system of $E$-rational 
Galois representations. 
\Endproof



The following Proposition will be essential in the later 
applications:

\begin{prop}\label{propeins} Let $E$ be a number field and let 
$k$ be a totally real number field. Let 
$(\psi_\lambda:G_k\to E_\lambda)_{\lambda \in \PP^f(E)}$ 
be a strictly compatible system 
of one-dimensional,  $\lambda$-adic, $E$-rational
Galois representations and 
let $\chi_\ell:G_k\to \QQ_\ell^\times$ 
be the $\ell$-adic cyclotomic character, where 
$\ell=\char(\lambda).$
Then there exists a finite character 
$\eps:G_k \to E^\times$ and an integer $m\in \ZZ$ such that 
$$\psi_\lambda = \eps \otimes
 \chi_\ell^m\,.$$
\end{prop}

\proof By \cite{Schappacher}, Prop.~1.4 in Chap. 1, 
any such compatible system
arises from an algebraic 
 Hecke character (the Prop.~1.4 in loc.~cit. is a consequence
of Henniart's result on the algebraicity 
of one-dimensional compatible systems, see \cite{Henniart81} and 
\cite{Khare03}). Any Hecke character of a totally real 
field with values in $E$ is equal to a power of the norm character 
times a finite 
order character $e$ (which has values in the group $\langle \zeta_d\rangle$
of roots 
of unity contained in $E$), see \cite{Schappacher}, Chap. 0.3. 
The compatible system associated to the 
 norm character is the system $(\chi_\ell)$ 
of $\ell$-adic cyclotomic characters.
Let $(\eps_\lambda)$ be the compatible system $(\eps_\lambda)$ associated 
to $e.$ By construction,  
 $$
\eps_\lambda(\Frob_p)=e(\,(p)\,)$$ whenever this expression is well defined 
(see \cite{Schappacher}).
Thus $\eps_\lambda$ takes values in the finite group
$\langle \zeta_d\rangle$ and is  thus (by Cebotarev density) a finite character
$\eps$ which is independent of $\lambda.$
\Endproof

The determinant of 
a compatible system of $E$-rational Galois representations 
is again  a compatible system of $E$-rational Galois representations.
Thus, by Rem.~\ref{remaha}, Prop.~\ref{corisogalmod2}, 
and 
Prop. ~\ref{propeins}, one finds:

\begin{thm}\label{thmdet}  Let $E_\lambda$ 
be the completion of a number field $E$ 
at a finite prime $\lambda$ of $E$ with $\chara(\lambda)=\ell.$ Let 
$k$ be a totally real number field. 
Let 
${\cal F}_1,\,{\cal F}_2,\,{\cal F}_3$ and $\GGG$ be 
irreducible and non-trivial \'etale 
local systems on punctured affine lines over $k$ having finite monodromy. Let 
$$ \V=(({\cal F}_1\ast {\cal F}_2)\otimes\GGG)\ast{\cal F}_3
\in \LS_{E_\lambda}^\et(S)\,.$$  Assume that 
the associated 
tuple 
$T_{\V}$ generates an infinite and absolutely irreducible subgroup 
of $\GL(\V_{{\bar{s}_0}}),$ where $s_0\in S(k).$
Let $\rho_\V:\pi_1^{\rm geo}(S,{\bar{s}_0})\rtimes G_k\to \GL(\V_{{\bar{s}_0}})$ be the 
Galois representation associated to $\V.$
Then there exists a finite character 
$\eps:G_k \to E^\times$ and an integer $m\in \ZZ$ such that 
$$\det(\rho_\V) = \det(\rho_\V|_{\pi_1^{\rm geo}(S,{\bar{s}_0})})\otimes \eps \otimes \chi_\ell^m\,.$$
\end{thm}

\section{Applications to the inverse Galois
problem}\label{secappl}

\subsection{Subgroups of linear groups.}\label{seclars}

In this section we collect some facts on linear groups which are needed
in later applications.\\

Let $V$ be a vector space. 
A subgroup $G\leq \GL(V)$ is called {\em primitive} if there
is no nontrivial direct sum decomposition $V=\bigoplus V_i$ which is
preserved by the elements of $G.$
See \cite{dr00}, Section 6, for a proof of the following
proposition which is useful in assuring the primitivity of some
groups occurring below:

\begin{prop}\label{ree2}
Let ${\bf T}=(T_1,\ldots,T_r)\in \GL_n(\FF_q)^r,\,q=p^s,$ be an absolutely irreducible tuple (i.e., $T_1,\ldots,T_r$ generate an absolutely 
irreducible subgroup of 
$\,\GL_n(\FF_q)$)
such
that $T_1\cdots T_r$ is a scalar and let 
$m=\sum_{i=1}^r{\rm rk}(T_i-1).$
Let $x\in {\Bbb N}$ be the maximal length of a Jordan block
occurring in the Jordan decompositions of $T_1,\ldots,T_r$ which is
not divisible by $p$ and let $V_1\oplus \cdots \oplus V_l$ be a
$\langle {\bf T} \rangle$-invariant decomposition of ${\Bbb
F}_q^{\;n}.$ Let \( \phi : \langle {\bf T} \rangle \mapsto S_l\) be
the induced map. Then $ \phi (T_i) =1 $ for $ {\rm rk} (T_i-1) <
{\rm dim}(V_1)$ and

$$ {\rm dim}(V_1)\geq
{\rm max}\{ x, n-m/2+1/2 (a+b) \}\,, $$ where
$$a= \sum_{
T_i\; {\rm semisimple}, \;
 {\rm rk}(T_i-1)< {\rm dim}(V_1)} {\rm rk}(T_i-1)$$
and
$$b=
\sum_{T_i\; \rm unipotent}\;\left( \sum_{\rm Jordan\; blocks} \mbox{\rm
length of Jordan blocks of $T_i$ not divisible by $p$}\right). $$
\end{prop}

 We call an element $X \in \GL(V)$ a
{\it perspectivity} (resp. {\it biperspectivity}) if $X$ is
semisimple with ${\rm rk}(X-1) =1$ (resp. ${\rm rk}(X-1) =2) $ or
if $X$ is unipotent with ${\rm rk}(X-1) =1 $ (resp. ${\rm rk}(X-1)
=2$ and $(X-1)^2=0).$ A semisimple perspectivity is called a {\it
homology}, a unipotent perspectivity is called a {\em transvection}.\\

The following result of Wagner (\cite{Wagner78}) will  also be used:

\begin{prop}\label{propWagnerachtundsieb}
Let $G$ be a primitive irreducible subgroup of $\GL_n(\FF_q),$ where 
 $n>2.$ If $G$ contains a homology of order $m>2,$ then
one of the following  holds:

\noindent (a)  $\SL_n(\FF_{q_0}) \leq G \leq Z(G) \cdot \GL_n(\FF_{q_0}),\;
\FF_{q_0} \leq \FF _q $
              and $m \mid (q_0-1).$

\noindent (b)  $\SU_n(\FF_{q_0}) \leq G \leq Z(G)\cdot \GU_n(\FF_{q_0}),\;
\FF_{q_0^2} \leq \FF_q $
              and $m \mid (q_0+1).$

\noindent (c)  $\SU_3(2) \leq G \leq Z(G)\cdot \GU_3(2),\; m=3,\;
n=3,\; p \not = 2$
              and $3 \mid (q-1).$

\noindent (d)  $\SU_4(2) \leq G \leq Z(G)\cdot \GU_4(2),\;  m=3,\;
n=4,\; p \not =
 2$
              and $3 \mid (q-1).$
\end{prop}

\subsection{Galois realizations of cyclic and dihedral groups.}

In this section we recall some known Galois realizations. They
will form the building blocks for new Galois realizations, using
the
middle convolution. \\

\begin{prop} \label{propreal1} \begin{enumerate}
\item Let $\uo:=\{x_1,x_2\},\, x_1\not= x_2,$ 
where  $x_1,x_2\in \PP^1(\QQ).$
 Then there exists an \'etale Galois cover
$f: {\frak F} \to \PP^1_\QQ\setminus \uo$ with Galois group  isomorphic to
$\ZZ/2\ZZ.$

\item Let $m\in \NN_{>2}.$  
Then there exist infinitely many elements
 $\vo \in \OO_{\varphi(m)}({\QQ})$ such that there is an
\'etale Galois cover $f:{\frak F}\to \AA^1_\QQ \setminus \vo$
 with Galois group isomorphic to
 $\ZZ/m\ZZ=\langle \rho \rangle$ whose associated tuple is
$$ \g_f=(\rho^{m_1},\ldots,\rho^{m_{\varphi(m)}},1)$$
(see Section~\ref{seccovv} for the notion of $\g_f$).
Here  $m_i \in \ZZ/m\ZZ^\times$  and $\varphi$ denotes Euler's $\varphi$-function.

\item  Let $m\in \NN_{>2},$ let  $\ZZ/m\ZZ=\langle \rho \rangle,$
and let $\ZZ/2\ZZ=\langle \sigma \rangle.$ Let
$D_m=\ZZ/m\ZZ \rtimes \ZZ/2\ZZ $ denote the
dihedral group of order $2m,$ where 
 $\rho^\sigma =\rho^{-1}$.
  Then there exists an element 
 $${\bf z}=\{x_1,\ldots,x_{2+\varphi(m)}\}\in \OO_{2+\varphi(m)}({\QQ})$$ and
an
\'etale Galois cover $f:{\frak F}\to \AA^1_\QQ \setminus {\bf z}$
with Galois group isomorphic to $D_m$ such that
$$ \g_f=(g_1,g_2,\rho^{m_1},\ldots,\rho^{m_{\varphi(m)}},1)\, ,$$
where $g_1,g_2$ are not contained in $\langle \rho
\rangle$ and where $m_i\in \ZZ/m\ZZ^\times.$ Moreover, $x_1$ and $x_2$ can be assumed to
be $\QQ$-rational points.

\item  We use the notation of {\rm (iii)}. Let $m$ be 
odd. 
For any $r\geq 2+\varphi(m)$
and (if $r> 2+\varphi(m)$)  
for any choice of elements $x_{3+\varphi(m)},\ldots, x_r\in
\AA^1(\QQ)$ such that 
$${\bf z}':={\bf z} \cup \{x_{3+\varphi(m)},\ldots, x_r\} \in \OO_r(\QQ)\, ,$$
there exists  an \'etale Galois cover $f:{\frak F}\to
\AA^1_\QQ \setminus {\bf z}'$ with Galois
group isomorphic to $D_{2m}=\langle \delta \rangle D_m$ (where 
$\delta$ is the central involution)  such that
$$ \g_f=(g_1,g_2,\rho^{m_1},\ldots,\rho^{m_{\varphi(m)}},\delta,\ldots,\delta)
\in D_{2m}^{r+1}. $$ 
\end{enumerate}
\end{prop}

\proof See \cite{Voelklein}, Chap. 7, or \cite{MalleMatzat}, Chap. I.5.1,
 for (i) and 
(ii). For claim (iii), note that 
the dihedral group $D_m$ is a factor group of 
the wreath product $H=\langle \rho\rangle \wr \langle \sigma \rangle,$
see \cite{MalleMatzat}, Prop. IV.2.3. Let 
 $\vo=\{y_1,\ldots,y_{\varphi(m)}\}$ be as in (ii) and let 
$$ U:=\{(x,y)\in \CC^2\mid y\not= -{x}{y_i}-y_i^2,\, i=1,\ldots, \varphi(m),
\,\, 
 x^2\not=4y\}\,.$$
The construction given in loc.~cit., proof of Prop. IV.2.1 (using 
  (ii)), shows that there 
exists an \'etale  Galois cover $\tilde{f}:X\to U$ 
with Galois group isomorphic to $H$  which is defined over 
$\QQ.$
 By factoring out
the kernel of the surjection $H\to D_m,$ one sees that 
there exists 
a Galois cover ${f}':X'\to U,$ with Galois group isomorphic to  $D_m.$
Claim (iii) now follows from restricting the cover $f'$ to a suitable
 (punctured)
line in $U$ which is defined over $\QQ.$
 Claim (iv) follows from (i) and
(iii) by taking  fibre products of covers. \Endproof

\subsection{Special linear groups as Galois groups.}

\begin{thm}\label{thmrealierung1} Let $\ell$ be a prime, let 
$q=\ell^s\,(s\in \NN),$ and let $n\in \NN.$ Assume that 
$$q\equiv 5\mod 8\quad \mbox{and that}\quad 
n> 6+2\varphi(m)\,,\quad \mbox{ where}\quad  
m:=(q-1)/4 \, .$$
Let 
$$E:=\QQ(\zeta_{m}+\zeta_{m}^{-1},\zeta_4)\, ,$$
where $\zeta_i\,\,(i\in \NN)$ denotes a primitive $i$-th root of unity. Let 
$\lambda$ be a prime of 
$E$ with $\char(\lambda)=\ell$ and let $O_\lambda$ be the valuation ring of 
the completion $E_\lambda.$ 
Then  the 
special linear  group $\SL_{2n+1}(O_\lambda)$ occurs regularly as
Galois group over $\QQ(t).$ \end{thm}

Before giving the proof of the theorem let us mention the
following corollary which immediately follows from reduction modulo
$\lambda$ :

\begin{cor} \label{corsldreia}
The special linear group $\SL_{2n+1}(\FF_q)$  
occurs regularly as Galois
group over $\QQ(t)$ if $$q\equiv 5\mod 8\quad \mbox{and}\quad 
n> 6+2\varphi((q-1)/4)\,.$$ \Endproof
\end{cor}

{\bf Proof of the theorem:} First we consider the case where $m\geq 3$  and 
$n=2r-4,$ where  $r\geq  2+\varphi({m}):$ Let 
$$f_1:{\frak F}_1\to \AA^1_\QQ\setminus
\uo_1,\quad \uo_1=\{x_1,\ldots,x_r\}\, ,$$ be a $D_{{2m}}$-cover
as in Prop.~\ref{propreal1}
(iv), where we assume that the points $x_1,\, x_2$ are
$\QQ$-rational. Let $\chi_1: D_{{2m}} \hookrightarrow
\GL_2(E_\lambda)$ be an orthogonal embedding of $D_{{2m}}.$ Thus
$$\FFF_1:=\L_{(f_1,\chi_1)}\in \LS^\et_{E_\lambda}(\AA^1_\QQ\setminus
\uo_1)$$
 is an \'etale local system of rank two (compare to  Section~\ref{seccovv}
for the notation of $\L_{(f_1,\chi_1)}$). The
associated tuple is
$$T_{\FFF_1}=(A_1,\ldots,A_{r+1})\in \GL_2(E_\lambda)^{r+1}\, ,$$
where $A_1,A_2$ are reflections, and $A_3,\cdots,A_{r+1}$ are
diagonal matrices  with eigenvalues
$$(\zeta_{{m}}^{m_1},\zeta_{{m}}^{-m_1}),\ldots,
(\zeta_{{m}}^{m_{\varphi({m})}},\zeta_{{m}}^{-m_{\varphi({m})}}),(-1,-1),
\ldots,(-1,-1)\, ,$$ where 
$\zeta_{{m}}^{m_1=1},\ldots,\zeta_{{m}}^{m_{\varphi({m})}}$ are 
the primitive powers of $\zeta_{{m}}$ (compare to our convention in Section
~\ref{secetaledef} for the notion of an associated tuple).

Let $f_2:{\frak F}_2 \to \AA^1_\QQ\setminus \uo_2,\, \uo_2:=\{0\},$ be a double cover
as
in Prop.~\ref{propreal1} (i), let $\chi_2:\ZZ/2\ZZ\hookrightarrow
E_\lambda^\times ,$  and let $$\FFF_2:=\L_{(f_2,\chi_2)}\in
\LS^\et_{E_\lambda}(\AA^1_\QQ\setminus \uo_2)\,.$$ Thus $\FFF_2$ is a
Kummer sheaf with $T_{\FFF_2}=(-1,-1).$

The middle convolution $\FFF_1\ast\FFF_2$ is an \'etale local
system on $\AA^1_\QQ\setminus \uo_1$ (since $\uo_1\ast
\uo_2=\uo_1$). By Prop.~\ref{dimension}, the rank of
$\FFF_1\ast\FFF_2$ is $2r-4$ and
$$T_{\F_1\ast\F_2}=(B_1,\ldots,B_{r+1})\in \GL_{2r-4}^{r+1}\, ,$$
where  (by Lemma~\ref{lemmonodromy1} and Rem.~\ref{remanaletal}) the matrices  $B_1,B_2$ are
transvections, the matrices $B_3,\ldots,B_{2+\varphi({m})}$ are
biperspectivities with non-trivial eigenvalues
$$(-\zeta_{{m}}^{m_1},-\zeta_{{m}}^{-m_1}),\ldots,
(-\zeta_{{m}}^{m_{\varphi({m})}},-\zeta_{{m}}^{m_{\varphi({m})}})\quad {\rm (respectively),}$$ the
matrices $B_{3+\varphi({m})},\ldots,B_r$ are unipotent
biperspectivities,  and where the last
matrix is equal to $-1$ by  Lemma  
\ref{lemmonodromy2} 
(see Section~\ref{seclars}, for 
the notion of a transvection and a (bi)perspecti\-vi\-ty).

Let $$\GGG_1:=\L_{(g,\xi)}|_{\AA^1_\QQ\setminus \uo_1}\in
\LS^\et_{E_\lambda}(\AA^1_\QQ\setminus \uo_1)$$ be the restriction of
the local system $\L_{(g,\xi)},$ where $g:{\frak G}\to \AA^1_\QQ\setminus
\{x_1\}$ is as in Prop.~\ref{propreal1} (i) and  $\xi : \ZZ/2\ZZ
\hookrightarrow E_\lambda^\times$.  Then the associated tuple of
$(\FFF_1\ast \FFF_2)\otimes\GGG_1$ is
$$ T_{(\FFF_1\ast \FFF_2)\otimes\GGG_1}=(-B_1,B_2,\ldots,B_r,-B_{r+1}=1)\,.$$

Let $f_3:{\frak F}_3 \to \AA^1_\QQ\setminus \uo_3$ be a $\ZZ/4\ZZ$-cover
as in Prop.
\ref{propreal1} (ii),  where we have chosen $\uo_3$
such that $\uo_1\ast \uo_3$ is generic in the sense of 
Section~\ref{sectconvdef}. Let
 $\chi_3:\ZZ/4\ZZ \hookrightarrow
E_\lambda^\times$  be an embedding and set $\FFF_3:=\L_{(f_3,\chi_3)}\in
\LS_{E_\lambda}^\et(\AA^1_\QQ\setminus \uo_3).$ Thus
$$T_{\FFF_3}=(i,-i,1),\quad {\rm where }\quad i:= \zeta_4\,.$$

Then the middle convolution
 $$\V:=((\FFF_1\ast \FFF_2)\otimes\GGG_1)\ast \FFF_3$$ 
is an \'etale local system
 on $S=\AA^1_\QQ\setminus \uo_1\ast \uo_3.$ The rank of
$\V$ is $4r-7=2n+1$ by Prop.
\ref{dimension}.  Let $T_{\V}=(C_1,\ldots,C_{2r+1}).$  By Lemma~\ref{lemmonodromy1},
$$C_1\sim J(-\zeta_4,2)\oplus_{k=3}^{2r-4}J(-\zeta_4,1)\oplus_{2r-4}^{4r-7}J(1,1)$$ (where we use the notation introduced in 
 Section~\ref{secdefacon}),
$C_2$ is a homology of order four, the
elements $C_3,\ldots,C_{2+\varphi({m})}$ are semisimple
biperspectivities with non-trivial eigenvalues
$$(-i\zeta_{{m}}^{m_1},-i\zeta_{{m}}^{-m_1}),\ldots,
(-i\zeta_{{m}}^{m_{\varphi({m})}},-i\zeta_{{m}}^{m_{\varphi({m})}})\, ,$$ the
elements $C_{3+\varphi({m})},\ldots,C_m$ are biperspectivities
with non-trivial eigenvalues $(i,i),$ and the matrix $C_{i+r}$ is
the Galois conjugate of the matrix $C_i\, (i=1,\ldots,r).$
Especially, one sees that $\langle C_1,\ldots,C_{2r}\rangle \leq
\SL_{4r-7}(E_\lambda)\times \langle \zeta_4 \rangle.$


Choose a  $\QQ$-rational 
point $s_0$ of $ S=\AA^1_\QQ\setminus \uo_1\ast
\uo_3.$ 
 By Thm.~\ref{thmdet} and 
Remark~\ref{rempoinabs} (ii) ($\V$ is punctually pure of weight $2$ by 
Prop.~\ref{propweilzwei}),  
$$\det(\rho_{{\V}})=\det(\rho_{{\V}}|_{\pi_1^{\rm geo}(S,{\bar{s}_0})})
\otimes \chi_\ell^{-(4r-7)}\otimes \epsilon\, ,$$
where $\epsilon:G_\QQ\to \langle
\zeta_4\rangle\subseteq E^\times $ and where $G_\QQ$ is embedded in 
$$\pi_1^\et(S,{\bar{s}_0})=\pi_1^{\rm geo}(S,{\bar{s}_0})\rtimes G_\QQ$$
by the choice of $s_0.$
Thus 
$$\im(\det(\rho_{{\V}(1)}))\leq \langle\zeta_4\rangle \, ,$$
where ${\V}(1)$ stands for the
Tate twist of $\V$ (see Def.~\ref{deftwist}). 
Let  $$\rho:=\rho_{{\V}}(1)\otimes \delta\;: \;\pi_1^{\rm geo}(S,{\bar{s}_0})\rtimes G_\QQ \To 
\GL_{4r-7}(E_\lambda)\, ,$$ where 
$\delta=\det(\rho_{{\V}(1)})^{-1}.$ It follows that 
\begin{equation}\label{eqgeomm}\im({\rho})\leq \SL_{4r-7}(E_\lambda).
\end{equation}
Let 
 ${H}=\im(\rho)$ and $\bar{H}:=\im(\bar{\rho}),$ where 
$\bar{\rho}$ is the residual representation of $\rho.$ 
Let further ${H}^\geom:=\im(\rho^\geom)$ and 
$\bar{H}^\geom:=\im(\bar{\rho}^\geom).$
By Thm.~\ref{thmirrd} and Rem.~\ref{remanaletal}, $\bar{H}^\geom$ is absolutely 
irreducible.
 Using Prop.~\ref{ree2} it is easy to see
that $\bar{H}^\geom$ is a primitive subgroup of
$\GL_{4r-7}(\FF_q).$ Thus, by the existence of the homology
$C_2,$  the theorem of Wagner (Prop.
\ref{propWagnerachtundsieb}), and 
the assumptions on $m,$
 one obtains $\SL_{4r-7}(\FF_q)\leq \bar{H}^\geom$
(where 
$\FF_q$ denotes the residue field).
By \eqref{eqgeomm}, 
$$\bar{H}^\geom=\SL_{4r-7}(\FF_q)\,.$$
Since $E_\lambda$ is unramified over 
$\QQ_\ell,$ the residual map 
$\SL_{4r-7}(O_\lambda)\to \SL_{4r-7}(\FF_q)$   is 
Frattini (see \cite{Wei}, Cor. A). Thus 
$${H}^\geom=\SL_{4r-7}(O_\lambda)=H\, ,$$ where the last equality follows 
from \eqref{eqgeomm}. 
The fundamental group
$\pi_1^\et(S,{\bar{s}_0})$ is a factor of
$G_{{\QQ}(t)}$ and $\pi_1^{\rm geo}(S,{\bar{s}_0})$ coincides with the image of $G_{\bar{\QQ}(t)}$ in
$\pi_1^\et(S,{\bar{s}_0}).$ 
Thus
the group
$\SL_{4r-7}(O_\lambda)$ occurs regularly as Galois group over
$\QQ(t).$ This proves the claim in the case $m\geq 3$ and $n=2r-4.$\\

The proof for $m\geq 3$  and $n=2r-3$ uses the following sheaves:
Let
$r\geq 4+\varphi({m}).$ Let $$\FFF_1:=\L_{(f_1,\chi_1)}\in \LS^\et_{E_\lambda}(\AA^1_\QQ\setminus
\uo_1)$$ be defined as above, where we assume 
that the points $x_{\varphi(m)+3}$ and 
$x_{\varphi(m)+4}$ coincide with $i,-i$ and 
that $x_1,x_2$ and $x_r$ are rational points. 
 By a suitable tensor operation,
one obtains a local system 
$\FFF_1'\in \LS^\et_{E_\lambda}(\AA^1_\QQ\setminus
\uo_1)$  whose
associated tuple is
$$T_{\FFF_1'}=(A_1',\ldots,A_{r+1}')\in \GL_2(E_\lambda)^{r+1}\, ,$$
where $A_1',A_2'$ are reflections, and $A_3',\cdots,A_{r+1}'$ are
diagonal matrices  with eigenvalues
\begin{multline}(\zeta_{{m}}^{m_1=1},\zeta_{{m}}^{-m_1}),\ldots,
(\zeta_{{m}}^{m_{\varphi({m})}},\zeta_{{m}}^{-m_{\varphi({m})}}),\\
\quad (i,i),\,(-i,-i),\,(-1,-1),\ldots,(-1,-1),(1,1))\,.\end{multline}
Let $\FFF_2=\L_{(f_2,\chi_2)}$ be the Kummer 
sheaf as above and let $\FFF_3':=\FFF_2.$ 
Let 
$\GGG'=\L_{(g',\xi')},$ where $g': {\frak G}'\to \AA^1_\QQ\setminus\{x_1,x_r\}$ 
is the double cover ramified at $x_1$ and $x_r$ and 
$\xi'$ is the embedding of the Galois group 
into $E_\lambda.$
Now continue as above, using the sheaf 
$$\V=((\FFF_1'\ast \FFF_2)\otimes \GGG'|_{\AA^1_\QQ\setminus \uo_1})\ast 
\FFF_3'$$ which has rank  $4r-5.$

The case where $m=1$ follows from  the same arguments as above
using the dihedral group $D_3$ instead of $D_m.$
\Endproof


\section{Galois realizations of $\PGL_2(\FF_\ell)$ 
and a family of $\Kd$-surfaces}\label{secspecil}

In this section, we construct regular Galois extensions 
$L_\ell/\QQ(t)$ with $$\Gal(L_\ell/\QQ(t))\simeq \PGL_2(\FF_\ell),\quad \ell
\not\equiv \pm 1,2 \mod 8\,,$$ and 
determine the behaviour of the Galois groups under the specialization $t\mapsto 1.$

\subsection{Computation of specializations of convolutions.}

Let $L/\QQ(t)$ be a finite Galois extension with 
Galois group
$G.$  By specializing 
$t$ to $s_0\in \QQ,$ one obtains a Galois extension $L'/\QQ$ 
with Galois group isomorphic to a subgroup $G_{s_0}$ of $G.$
 
It is the aim of this subsection 
to give a procedure for computing the Galois groups of 
specializations of Galois extensions of ${\QQ(t)}$ 
which are constructed 
via the convolution. \\

Let  $$ \V=(\cdots ((({\cal F}_1\ast {\cal F}_2)\otimes\GGG_1)\ast{\cal F}_3)\otimes 
\cdots \ast {\cal F}_n)\otimes \GGG_{n-1}$$ be the iterated convolution
of irreducible local systems with finite monodromy as in 
Thm.~\ref{thmpureweight}. Assume that the 
underlying field $k$ is the field of rational numbers. Let 
\begin{multline}\nonumber 
\Pi=\pr^{[1,n]}_{[n,n]}\otimes \pi_n: {\frak U}:=
{\frak F}_1\otimes\cdots \otimes {\frak F}_{n} \otimes 
{\frak G}_1\otimes \cdots\otimes  {\frak G}_{n-1}\\
 \To S=\AA^1_\QQ\setminus \uo_1\ast \cdots \ast \uo_n\end{multline}
  be the associated affine morphism over $S$ (see Section~\ref{secast}).

 Now assume that
$L\subseteq \overline{\QQ(t)}$ is the fixed field of the kernel of 
the residual representation  
$$\bar{\rho}_{\V}\,:\,\pi_1^\et(S,\bar{s}_0)\,\To\, \GL(\overline{\V}_{{\bar{s}_0}})\,.$$
 In order to determine the specialization $G_{\bar{s}_0}$ 
of $G=\Gal(L/\QQ(t)),$
 one can proceed as 
follows:

\begin{itemize}
\item As in Section~\ref{secgeommi}, compute an equivariant 
normal crossings compactification $\tilde{\CX}_{\bar{s}_0}$ of 
the fibre ${\frak U}_{\bar{s}_0}.$
\item Compute the cohomology 
$H^{n-1}(\tilde{\CX}_{\bar{s}_0},E_\lambda)$ and  
 determine $\V_{\bar{s}_0}$ as a subfactor of
$H^{n-1}(\tilde{\CX}_{\bar{s}_0},E_\lambda)$ using the projector 
$\PPP_{\Delta_n}$ (compare to Cor.~\ref{corisogalmod}).
\item Using this  description of $\V_{\bar{s}_0}$
as a subfactor of  $H^{n-1}(\tilde{\CX}_{\bar{s}_0},E_\lambda)$ 
 and using the number of $\FF_p$-rational points
of $\tilde{\CX}_{\bar{s}_0}\otimes \bar{\FF}_p,$
compute the trace of 
$$\Frob_p\in 
G_\QQ\leq \pi_1^\et(S,\bar{s}_0)=\pi_1^\geom(S,\bar{s}_0)\rtimes G_\QQ$$
 on $\V_{\bar{s}_0}$ for 
sufficiently many primes $p$ where $\tilde{\CX}_{\bar{s}_0}$ has good reduction.
\item Determine the Hodge numbers of (the de Rham version of)
$\V_{\bar{s}_0}.$ Using these Hodge numbers and the 
results of Section~\ref{tameinertia} below on crystalline representations, 
compute the action 
of the tame inertia group 
$$I^t_\ell\leq G_\QQ\leq \pi_1^\et(S,\bar{s}_0)$$ on $\V_{\bar{s}_0}.$
\item Then use  this information and group theory in order
to determine the image of 
$ \bar{\rho}_{\V}|_{G_\QQ}.$ One then has
$$G_{\bar{s}_0}=\im(\bar{\rho}_{\V}|_{G_\QQ})\,.$$
\end{itemize}

\begin{rem}{\rm  \begin{enumerate}
\item The above procedure  equally works in the case of  Tate twists.

\item The above  procedure is illustrated in the next sections
in the case of a Tate twist of a rank three orthogonal 
local system which arises from the convolution of three rank 
one systems. \end{enumerate}}
\end{rem}

\subsection{The underlying local systems and Galois extensions.}\label{secunderl}

In this section, we derive regular Galois extensions of $\QQ(z)$ 
whose Galois groups are  projective linear groups. This is done via the 
calculus of convolution of local systems. The specializations of these 
Galois extensions are considered in Section~\ref{secspeck}.\\

Let 
$$f: {\frak F}\To \AA^1_\QQ\setminus \uo,\quad {\rm where}\quad 
\uo=\{\pm 1\}\,,$$ be a Galois cover with Galois group 
 $G=\ZZ/2\ZZ=\langle \sigma \rangle$ and associated tuple
$\g_f:=(\sigma,\sigma,1),$ see Section~\ref{seccovv} and Prop. 
\ref{propreal1}. 
 Let $\chi:G \hookrightarrow \QQ_\ell $ and let 
$$\LL:=\LL_{(f,\chi)}\in \LS^\et_{\QQ_\ell}(\AA^1_\QQ\setminus \{\pm 1\})$$
(see Section~\ref{seccovv} for the notation $\LL_{(f,\chi)}$).
 
Let  $$\LL_{-1}=\LL_{(f',\chi)}\in \LS^{\et}_{\QQ_\ell}
(\AA^1_\QQ\setminus \{0\})$$ be the 
Kummer sheaf with $T_{\LL_{-1}}=(-1,-1),$ see Section~\ref{sectconvdef}
(i.e., $f':{\frak F}'\to \AA^1_\QQ\setminus\uo',\, \uo'= \{0\},$ is the double cover 
with associated tuple $(\sigma,\sigma)$).

Let 
$$ \V:=(\LL\ast \LL)\ast \LL_{-1}\in  \LS^\et_{\QQ_\ell}(S)\, ,$$
where $$S:=\AA^1_{\QQ,z}\setminus \uo\ast\uo\ast\uo'=
\AA^1_{\QQ,z}\setminus \{0,\pm 2\}\,.$$

Using the Convolution Program  \cite{DettweilerConv}
one verifies that
$\LL\ast \LL$  is a $2$-dimensional symplectic 
\'etale local system of $\QQ_\ell$-modules 
on $\AA^1_\QQ\setminus \{\pm 1,0\},$
whose associated tuple  $T_{\LL\ast \LL}$ is 
$$\left(\; \left(\begin{array}{rrr} -3& -8\\
2&5\end{array}\right),\;\;\left(\begin{array}{rrr} 1& -4\\
0&1\end{array}\right),\;\;\left(\begin{array}{ccc} 1& 0\\
2&1\end{array}\right),\;\;\left(\begin{array}{ccc} -3& -4\\
4&5\end{array}\right)\, \right)\,.$$
Thus,   for $\ell\not=2,$ the residual 
representation $\bar{\rho}_{\L\ast\L}$ is  absolutely
irreducible.\\

Thus
$\V$  is an irreducible (by Rem.~\ref{remlab} and 
\cite{dr00}, Cor. 3.6, or by Thm.~\ref{thmirrd})  $3$-dimensional
\'etale local system of $\QQ_\ell$-modules 
on $\AA^1_\QQ\setminus \{\pm 2,0\}.$ Moreover,
the associated tuple is $T_{\V}=(M_1,\ldots,M_4),$ where  
$$M_1= \left(\begin{array}{rrr} -1& -4&4\\
0&1&0\\
0&0&1\end{array}\right),\quad M_2:=\left(\begin{array}{rrr} 1& 0&0\\
-2&-1&2\\
0&0&1\end{array}\right)\, ,$$
$$M_3:=\left(\begin{array}{ccc} 1& 0&0\\
0&1&0\\
4&4&-1\end{array}\right),\quad M_4:=\left(\begin{array}{ccc} -1& -4&4\\
2&7&-6\\
4&12&-9\end{array}\right)\,.$$

Set 
$$\rho_\ell:=\rho_{\V}\,:\,\pi_1^\et(S,\bar{1})= \pi_1^\geom(S,\bar{1})\rtimes G_\QQ\To \GL_3(\QQ_\ell)\,.$$
By Propositions~\ref{proppoin} and \ref{fform}, one obtains
 a bilinear pairing 
$$ \V\otimes \V\longrightarrow \QQ_\ell(-2)$$ such that
 the image of  the 
restriction of $\rho_\ell$
to $\pi_1^\geom(S)$  is contained in the general orthogonal 
group $\O_3(\QQ_\ell).$ 
Taking Tate twists, one obtains
an orthogonal pairing 
\begin{equation} \label{eqtwist}\V(1)\otimes \V(1)\longrightarrow \QQ_\ell\,.\end{equation}
Thus  the whole 
image of $\rho_\ell(1)=\rho_\ell\otimes \chi_\ell$ is contained in the 
 orthogonal group $\O_3(\QQ_\ell),$ 
compare to Rem.~\ref{rempoinabs}.
Let $\bar{\V}$ and $\bar{\V}(1)$
denote the local system of 
$\FF_\ell$-modules corresponding to the residual representation 
$\bar{\rho}_\ell$ and its Tate twist 
$\bar{\rho}_\ell(1)=\bar{\rho}_\ell\otimes \bar{\chi}_\ell,$ respectively.

\begin{prop}\label{proporthog} Let  $\ell\not\equiv \pm 1\mod 8$ be an odd prime.
Then 
$$\langle \bar{M}_1,\ldots,\bar{M}_4\rangle=\O_3(\FF_\ell)\, ,$$
where $$(\bar{M}_1,\ldots,\bar{M}_4)=T_{\bar{\V}}=T_{\bar{\V}(1)}$$
is the reduction modulo $\ell$ of the tuple $({M}_1,\ldots,{M}_4)$  above.
\end{prop}

\proof The monodromy tuple $(\bar{M}_1,\ldots,\bar{M}_4)$ consists of three 
reflections $\bar{M}_1,\bar{M}_2$ and $\bar{M}_3$ which generate an irreducible 
subgroup and whose product is the negative 
of a unipotent element of unipotent rank two.
Thus the commutator subgroup ${\rm SO}_3(\FF_\ell)'$
of $\SO_3(\FF_\ell)$ 
is contained in $\langle \bar{M}_1,\ldots,\bar{M}_4\rangle.$ 
Moreover, it can be checked that the matrix $M_1M_2$ has eigenvalues 
$$(\alpha_1,\,\alpha_2,\,\alpha_3)= (3+ 2\sqrt{2},\,3-2\sqrt{2},\,1)\,.$$ 
Thus, by reducing the coefficients modulo $\ell,$ one sees that 
$\bar{M}_1\bar{M}_2$ is an element in the torus $T$ 
of $\SO_3(\FF_\ell)$ of order 
$\ell+1$ if  $2$ is not a square modulo $\ell$ (which is 
equivalent to $\ell\not\equiv \pm 1\mod 8$). 
The element $3+2\sqrt{2}$ is the square of 
$1+\sqrt{2}.$ Moreover, 
$$(1+\sqrt{2})(1-\sqrt{2})=-1\,.$$ Thus, by an argument using the 
norm, the element $\bar{M}_1\bar{M}_2$
is not the square of an element of $T$. Since the spinor norm
$\SO_3(\FF_\ell)\to \{\pm 1\}$ 
 restricted to $T$ is surjective and has value $-1$ at 
$\bar{M}\in T$ if and only if $\bar{M}$ is a square in $T,$ it follows that 
$\SO_3(\FF_\ell)$ is contained in $\langle \bar{M}_1,\ldots,\bar{M}_4\rangle$
(the commutator subgroup ${\rm SO}_3(\FF_\ell)'$ is the kernel of the 
spinor norm). Thus $\langle \bar{M}_1,\ldots,\bar{M}_4\rangle=\O_3(\FF_\ell).$
\Endproof

The last result yields a new proof of the 
following classical result (compare to \cite{Matzat87}, Folgerung 2, p. 181):
 
\begin{cor}\label{cordeltal} For any $\ell\not\equiv \pm 1,2\mod 8,$ 
there 
exists a regular Galois extension $L_\ell/\QQ(z)$ 
with $\Gal(L_\ell/\QQ(z))\simeq \PGL_2(\FF_\ell).$
\end{cor}

\proof By Prop.~\ref{proporthog} and Equation \eqref{eqtwist},
for any $\ell\not\equiv \pm 1,2\mod 8,$ there is  a surjective
Galois representation $\bar{\rho}_\ell(1):\pi_1^\et(S,\bar{1})\to \O_3(\FF_\ell).$ 
Since we are in dimension $3,$ 
$$\im\left(\,\bar{\rho}_\ell(1)\otimes \det(\bar{\rho}_\ell(1))\,\right)= \SO_3(\FF_\ell)\,.$$ Let $\mu_\ell:G_{\QQ(z)}\to \SO_3(\FF_\ell)$ be the 
composition of the canonical surjection $G_{\QQ(z)}\to \pi_1^\et(S,\bar{1})$
with $\bar{\rho}_\ell(1)\otimes \det(\bar{\rho}_\ell(1)).$
It follows that  the  
subfield $L_\ell$ 
of $\overline{\QQ(z)}$ fixed by 
the kernel of $\mu_\ell$
is a  Galois extension with Galois group isomorphic 
to $\SO_3(\FF_\ell)\simeq \PGL_2(\FF_\ell)$ which is regular
since the image of $\pi_1^{\rm geo}(S)$ coincides with 
$\SO_3(\FF_\ell).$ \Endproof

\subsection{The underlying $\Kd$-surface.}\label{seckdrei}

By Section~\ref{secast}, the local system $\V$ is a subfactor 
of $R^2\Pi_*(\QQ_l),$ where 
$ \Pi: {\frak U}\To S$ is a smooth affine morphism of relative 
dimension two. It is the aim of this section to determine $\Pi$ and 
to compute a normal crossings compactification of the fibre 
${\frak U}_{{1}}$ of ${\frak U}$ over the point $1\in S(\QQ).$ \\

By adapting the notation of  Section 
\ref{secmotiv} (using $x,y,z$ instead of $x_1,x_2,x_3$
and using the local systems ${\frak F},{\frak F},{\frak F}'$ of Section 
\ref{secunderl} instead of 
$\bF_1,\bF_2,\bF_3$),
let  $$\OO_{[1,3]}:=\{(x,y,z)\mid\, x, y-x \not= \pm 1;\,\,
y,z \not=0,\pm 2;\,\, z-y \not=0\}\, ,$$ let 
$$\pr=\pr^{[1,3]}_{[3,3]}:\OO_{[1,3]}\To S$$ be the third projection,
 and let 
$$ {\frak U}:=\{(x,y,w,z)\in \AA^1_\QQ\times \OO_{[1,3]}
\mid\, w^2=(x^2-1)((y-x)^2-1)
(z-y)\,\}\,.$$
 The map $$\pi_3:{\frak U}\To \OO_{[1,3]}\,,\quad  (x,y,w,z)\Mapsto 
(x,y,z)$$ is an \'etale Galois 
cover with Galois group $G=\langle \sigma \rangle \simeq \ZZ/2\ZZ.$
Thus ${\frak U}$ is nothing else than 
the variety ${\frak F}\otimes {\frak F}\otimes {\frak F}'$ in the sense of Section 
\ref{secast} and $\pr \circ\pi_3:{\frak U}\to S$ coincides with the 
map $\Pi.$\\

Let ${\frak U}_{1}$ denote the fibre of ${\frak U}$ over $1\in S(\QQ).$ 
Clearly, the fibre ${\frak U}_{1}$ embeds into a (singular) double cover 
$\tau:\CX_1 \to \PP^2_{[x,y,v]}$ with ramification locus 
$$ B:=B_1\cup \ldots \cup B_6\,,$$ 
where
$$\begin{array}{l}
B_1=\{[x,y,v]\in \PP^2 \mid \, v=0\}\,,\\
B_2=\{[x,y,v]\in \PP^2 \mid \, x=-v\}\,,\\
B_3=\{[x,y,v]\in \PP^2 \mid \, x=v\}\,,\\
B_4=\{[x,y,v]\in \PP^2 \mid \, y=v\}\,,\\
B_5=\{[x,y,v]\in \PP^2 \mid \, y=x+v\}\,,\\
B_6=\{[x,y,v]\in \PP^2 \mid \, y=x-v\}\,.\\
\end{array}$$
The singular locus of $B$ is  the union 
of the $11$ points
$$\begin{array}{llll} L_1=[1,0,0]\,,& L_2=[0,1,0]\,,&L_3=[1,1,0]\,,& 
L_4=[1,2,1]\,,\\
L_5=[-1,0,1]\,,&L_6=[1,0,1]\,,&L_7=[-2,-1,1]\,,&L_8=[-1,1,1]\,,\\
L_9=[1,1,1]\,,& L_{10}=[1,2,1]\,,& L_{11}=[1,0,1]\,.
\end{array}$$

The following resolution process of the singularities 
of $\CX_1$ is the canonical resolution 
of $\CX_1$ in the sense of \cite{BPV}, Chap. III.7:\\

Let $\gamma:{\PP^2}' \to \PP^2$ be the 
 blow up of  $\PP^2$
 at $L_1,L_2,\ldots,L_{11}.$  Let 
 $\tau':\CX'\to {\PP^2}'$ be the pullback of the cover $\tau$ along 
$\gamma.$  
 Then $\tau'$ is ramified over the exceptional divisors
$E_2,E_3$ over $L_2$ and $L_3$ (resp.) and over the strict transform $B'$ of 
$B.$ The union $B'\cup E_2\cup E_3$  is 
a normal crossings divisor $N,$ whose singular locus is
where the strict transform 
 of the lines $B_1,B_2,B_3$ meets $E_2$ and where
the strict transform of the lines $B_1,B_5,B_6$ meet $E_3.$ Blowing 
up  ${\PP^2}'$ at these intersections, one obtains 
a map $\tau'': \tilde{\PP}^2\to {\PP^2}'.$
Let $\tilde{B}$ be the strict transform 
of $\tilde{B}:=B'\cup E_2\cup E_3$ and let  
 $\tilde{\tau}:\tilde{\CX}_{1}\to \tilde{\PP}^2$ be the pullback 
of $\tau'$ along $\tau''.$  By construction,
$\tilde{\tau}$ is a double cover which has no singularities since the 
ramification locus 
$\tilde{B}$ of the covering $\tilde{\CX}_{1}\to \tilde{\PP}^2$ 
has no singularities.

\begin{prop}\label{proppicard}  
\begin{enumerate}
\item Let $\sigma: \tilde{\CX}_1 \to {\CX}_1$ be the natural map.
 The  variety $\tilde{\CX}_1$ is a $G$-equivariant normal crossings 
desingularization of  $\CX_1$ and $\tilde{\CX}_1 \setminus \sigma^{-1}({\frak U})$ is a 
normal crossings divisor. 
\item 
The variety $\tilde{\CX}_{\bar{1}}$ is a 
 $\,{\rm K}3$-surface of Picard number $\geq 19$ (see 
Section~\ref{secconse} for the definition of the Picard number).
\end{enumerate}
\end{prop}

\proof The claim (i) immediately follows from the above 
construction of $\sigma.$ 

For (ii), note that all the occurring 
singularities of $B$ are simple in the sense of 
 \cite{BPV}. Thus, by inserting the degree 
$6$ of $B$ into the formulas  on p. 183 of loc.~cit., one sees that
 the minimal resolution of ${\CX}_{\bar{1}}$ is a $\Kd$-surface. By 
construction, the variety  $\tilde{\CX}_{\bar{1}}$ is the canonical 
(hence minimal, since all singularities are 
simple) resolution of the double cover ${\CX}_{\bar{1}} \to \PP^2,$
compare to \cite{BPV}, Chap. III.7. Hence $\tilde{\CX}_{\bar{1}}$ is a ${\rm K}3$-surface. 

To prove the claim on the Picard number of $\tilde{\CX}_{\bar{1}},$
 one considers 
the following divisors on $\tilde{\CX}_{\bar{1}}:$ The canonical resolution
process yields $9$ copies of $\PP^1$ $P_1,\ldots,P_9$
over the double points 
of $B$ and $2\cdot 4$ copies of $\PP^1$ $P_{10},\ldots,P_{13},$
resp. $P_{14},\ldots,P_{17},$
over the two triple points 
of $B.$ Let $P_{18},\,P_{19}$ be the two copies of $\PP^1$ over 
the line $y=0.$ The intersection matrix $I$ of $P_1,\ldots,P_{19}$ 
is the following matrix
(up to a suitable ordering of $P_1,\ldots,P_9$ and of 
$P_{10},\ldots,P_{17}$): 
\begin{center}
\begin{tiny}$$\left(\begin{array}{ccccccccccccccccccccc}

-2&0&0&0&0&0&0&0&0&0&0&0&0&0&0&0&0&1&1\\
0&-2&0&0&0&0&0&0&0&0&0&0&0&0&0&0&0&1&1\\
0&0&-2&0&0&0&0&0&0&0&0&0&0&0&0&0&0&1&1\\
0&0&0&-2&0&0&0&0&0&0&0&0&0&0&0&0&0&0&0\\
0&0&0&0&-2&0&0&0&0&0&0&0&0&0&0&0&0&0&0\\
0&0&0&0&0&-2&0&0&0&0&0&0&0&0&0&0&0&0&0\\
0&0&0&0&0&0&-2&0&0&0&0&0&0&0&0&0&0&0&0\\
0&0&0&0&0&0&0&-2&0&0&0&0&0&0&0&0&0&0&0\\
0&0&0&0&0&0&0&0&-2&0&0&0&0&0&0&0&0&0&0\\
0&0&0&0&0&0&0&0&0& -2&1&1&1&0&0&0&0&0&0\\
0&0&0&0&0&0&0&0&0& 1&-2&0&0&0&0&0&0&0&0\\
0&0&0&0&0&0&0&0&0& 1&0&-2&0&0&0&0&0&0&0\\
0&0&0&0&0&0&0&0&0& 1&0&0&-2&0&0&0&0&0&0\\
0&0&0&0&0&0&0&0&0& 0&0&0&0&-2&1&1&1&0&0\\
0&0&0&0&0&0&0&0&0& 0&0&0&0&1&-2&0&0&0&0\\
0&0&0&0&0&0&0&0&0& 0&0&0&0&1&0&-2&0&0&0\\
0&0&0&0&0&0&0&0&0& 0&0&0&0&1&0&0&-2&0&0\\
1&1&1&0&0&0&0&0&0& 0&0&0&0&0&0&0&0&-2&x\\
1&1&1&0&0&0&0&0&0& 0&0&0&0&0&0&0&0&x&-2
\end{array}\right)$$\end{tiny}\end{center}
\vspace{.5cm}
(here,  $x$ denotes the intersection number 
of $P_{18}$ and $P_{19}$). The last claim can 
be seen using the following arguments: 
By \cite{BPV}, Chap. III.7 (Table 1),  the divisors $P_1,\ldots,P_{17}$
are curves with self intersection number $-2.$ The intersection 
behaviour of $P_{10},\ldots,P_{13}$ and of 
$P_{14},\ldots,P_{17}$ can also 
be read off from loc.~cit., Table 1. Up to an ordering, we can assume 
that the intersection of $P_{18},$ resp. $P_{19},$ with $P_1,\,P_2,\,P_3$
is equal to one and zero for $P_i,\,i=4,\ldots,17.$ 
The pencil of lines through $[1:0:0]$ 
gives rise to an elliptic fibration
$\epsilon: \tilde{\CX}_{\bar{1}}\to
 \PP^1_{y,v}.$
 Zariski's 
lemma (see \cite{BPV}, Lemma III.8.1)
applied to $\epsilon$ shows that the self-intersection 
of $P_{18}$ and $P_{19}$ is $< 0.$   Thus, by 
loc.~cit., Prop. VIII.3.6, $P_{i}^2=-2,\,i=18,19.$ 
Zariski's lemma also implies that 
the intersection number $x$ of $P_{18}$ and $P_{19}$ is 
$0,1$ or $2$ (compare to loc.~cit., Chap. V.7). The determinant
of the intersection matrix $I$ is 
$$ 8192x^2 + 24576x + 16384\, ,$$ thus $I$ is invertible.
 Consequently,
the divisors $P_1,\ldots,P_{19}$ span a submodule 
of the Neron-Severi group of rank $19.$ 
\Endproof

 Consider the Neron-Severi group 
$$\NS(\tilde{\CX}_{\bar{1}})_\ell:=
\NS(\tilde{\CX}_{\bar{1}})\otimes \QQ_\ell\leq 
H^2(\tilde{\CX}_{\bar{1}},\QQ_\ell(1))\,.$$ The restriction of the
Poincar\'e pairing on $H^2(\tilde{\CX}_{\bar{1}},\QQ_\ell(1))$ to 
$\NS(\tilde{\CX}_{\bar{1}})_\ell$ coincides with the intersection pairing 
and is known to be non-degenerate, see Rem.~\ref{remns}. 
Since the Poincar\'e pairing on $H^2(\tilde{\CX}_{\bar{1}},\QQ_\ell(1))$
is compatible with the action of $G_\QQ,$  the orthogonal complement
$\NS(\tilde{\CX}_{\bar{1}})^\perp_\ell$ (of $\NS(\tilde{\CX}_{\bar{1}})_\ell$ in 
$H^2(\tilde{\CX}_{\bar{1}},\QQ_\ell(1))$)  can be seen as a $G_\QQ$-module in 
a natural way. 

Since $\tilde{\CX}_{\bar{1}}$ is a ${\rm K}3$-surface, 
the dimension of $H^2(\tilde{\CX}_{\bar{1}},\QQ_\ell(1))$ is equal to 
$22.$ It follows that the dimension of 
$\NS(\tilde{\CX}_{\bar{1}})^\perp_\ell$ is smaller or equal to $3.$

\begin{thm}\label{thmhodgetyp} 
There exists an isomorphism of $G_\QQ$-modules:
$$ \NS( \tilde{\CX}_{\bar{1}})^\perp_\ell\simeq \V_{\bar{1}}(1)\quad \Longleftrightarrow \quad \NS( \tilde{\CX}_{\bar{1}})^\perp_\ell(-1)\simeq \V_{\bar{1}}\,.$$
\end{thm}

\proof  Let $$\PPP=
\PPP_{\chi\otimes \chi\otimes \chi'}=1/2(1-\sigma)\in 
\End(\tilde{\CX}_{\bar{1}})\otimes \QQ\, ,$$ where
$\sigma$ is as in the last section. Then $\PPP$ annihilates a subspace of 
$\NS( \tilde{\CX}_{\bar{1}})_\ell$ (and of 
$\NS( \tilde{\CX}_{\bar{1}})_\ell(-1)$) of dimension $18$  (spanned  
by the divisors $P_1,\ldots,P_{17},P_{18}+P_{19}$ which are as in the 
proof of Prop.  \ref{proppicard}). 
By Thm.~\ref{thmmot1} and by Thm.~\ref{thmpureweight}, 
$\V_{\bar{1}}(1)$ is fixed by the projector 
 $\PPP.$ 
 Since  $\V_{\bar{1}}\leq H^2(\CX_{\bar{1}},\QQ_\ell)$
is (as $\pi_1^\geom(S,\bar{1})$-module) irreducible of 
rank three,  it has to coincide 
with the complement of the Neron-Severi group (compare to the proof of 
Cor.~\ref{corisogalmod2}). 
\Endproof

\begin{rem}{\rm 
\begin{enumerate}
\item The proof of the last theorem and 
 Prop.~\ref{proppicard} imply that 
the Picard number of $\tilde{\CX}_{\bar{1}}$ is equal to $19.$
\item Using exactly the same arguments as above, one can show that 
$\tilde{\CX}_{\bar{1}}$ is contained in a family $\tilde{\CX}/S$
of ${\rm K}3$-surfaces of Picard number $19.$ 
\end{enumerate} }
\end{rem}

\subsection{Computation of Frobenius elements.}

\begin{prop}\label{propfrobb} Let  $n\in \NN_{>0},$ let $p>3$ be a  prime, and 
let $q=p^n.$
Let 
$$N(q):=\#\{(w,x,y)\in \FF_{q}^{\;3}
\mid w^2=(x^2-1)((y-x)^2-1)(y-1)\}$$
and let 
$$\rho_\ell=\rho_\V: \pi_1^\et(S,\bar{1})\To \GL(\V_{\bar{1}})\simeq 
 \GL_3(\QQ_\ell)\,,$$
where $\V$ is as in Section~\ref{secunderl}. Let 
$\Frob_p$ be a Frobenius element in $G_\QQ$ which is viewed 
as an element in $$\pi_1^\et(S,\bar{1})
=\pi_1^\geom(S,\bar{1})\rtimes G_\QQ\,,$$ and let $\Frob_q=\Frob_p^n.$
 If $\ell \not= p,$ then the Galois representation 
$\rho_\ell|_{G_\QQ}$ is unramified at $p$
and 
the  trace of $\rho_{\ell}(\Frob_{q})$ is given as follows:
$$ \tr( \rho_{\ell}(\Frob_{q}))=N(q)+q-q^{2}-
\left(1+\left( \frac{-1}{q}\right)\;\right)\cdot q\,,$$
where $\left( \frac{\cdot}{\cdot }\right)$ denotes the Legendre symbol.
\end{prop}

\proof  Let $B$ be the ramification divisor of the (singular) 
double cover  $\CX_1\to \PP^2.$
Since $p>3,$ the reduction modulo $p$ of 
$B$ consists of $6$ reduced irreducible components which intersect 
at $11$ different points. It 
follows that the 
 desingularization process given in Section~\ref{seckdrei} is
equally  valid in characteristic $p$ and that 
$\tilde{\CX}_1$  has good reduction at $p$
(where $\tilde{\CX}_1$ is denotes the desingularization of
${\CX}_1$ as in Section~\ref{seckdrei}). 
Consequently, the Galois representation $\rho_\ell|_{G_\QQ}$ is unramified 
at $p.$
The desingularization process also yields 
 divisors $$\bar{P}_i\subseteq \tilde{\CX}_1\otimes \bar{\FF}_{{q}},\quad i=1,\ldots,19\,,$$ which coincide with 
 the reduction of the divisors $P_i$ which occur in 
the proof of Prop.~\ref{proppicard}.

 Let 
$\tilde{N}(q)$ be the 
number of $\FF_{q}$-rational points of $\tilde{\CX}_1\otimes \FF_{{q}}.$
The base change theorem and the Lefschetz trace formula imply that 
$$\tilde{N}({q})=1+q^{2}+\tr\left(\Frob_{q}|_{\NS(\tilde{\CX}_{\bar{1}})_\ell(-1)}\right)
+\tr( \rho_{\ell}\left(\Frob_{q})\right)$$
(compare to Thm.~\ref{thmhodgetyp}).
The divisors $\bar{P}_1,\ldots,\bar{P}_{17}$ are defined over $\FF_{{q}}$ and 
so
each one contributes a summand ${q}$ to the trace  
of $\Frob_{q}$ on 
$$\NS(\tilde{\CX}_{\bar{1}})_\ell(-1)\,\leq\, H^2(\tilde{\CX}_{\bar{1}},\QQ_\ell)\,=\,
\NS(\tilde{\CX}_{\bar{1}})_\ell(-1)\oplus \V_{\bar{1}}.$$ The divisors $\bar{P}_{18},\,
\bar{P}_{19}$
are defined over $\FF_{{q}}$ if and only if the Legendre 
symbol
$$ \left( \frac{-1}{{q}}\right)$$ is equal to one. Thus they 
contribute $$ \left(1+\left( \frac{-1}{{q}}\right)\right){q}$$ to the 
trace of $\Frob_{q}$ on 
${\NS(\tilde{\CX}_{\bar{1}})_\ell}(-1)$ and one 
finds
$$ \tr( \rho_{\ell}(\Frob_{q}))=\tilde{N}({q})-1-q^{2}-17{q}-
\left(1+\left( \frac{-1}{{q}}\right)\right){q}\,.$$

Each double (resp. each triple) point of the divisor $B\otimes \bar{\FF}_{{q}}$
 is defined over 
$\FF_{{q}}.$ 
So the contribution of the exceptional lines 
of the blow ups of the double points 
to the number of $\FF_{{q}}$-rational points of $\tilde{\CX}_1\otimes 
\bar{\FF}_{{q}}$
 sums up to $9({q}+1).$
 The blow ups of the two triple points 
yield each $4$ copies of $\PP^1$ which are defined over 
$\FF_{{q}}$ and where three of the copies of $\PP^1$ 
intersect the remaining one 
at different  rational points. Thus the two triple 
points contribute $8({q}+1)-6$ to the  
 number of $\FF_{{q}}$-rational points of 
$\tilde{\CX}_1\otimes 
\bar{\FF}_{{q}}.$ The line at  infinity away from the singularities
is in the ramification locus 
and contributes ${q}-2$ points. Consequently,
\begin{eqnarray}\nonumber 
 \tilde{N}({q})&=&N({q}) -8+({q}-2)+17({q}+1)-6\nonumber \\
&=&N({q})+1+18{q}\,.\nonumber \end{eqnarray}
Adding everything yields the claim. 
\Endproof

\begin{cor}\label{corfrobat}
 The eigenvalues 
of the Frobenius elements $$\rho_\ell(1)(\Frob_p),\quad {\rm for}\quad 
\ell\not=p \quad {\rm and}\quad  3<p\leq 29\, ,$$
are 
$$\alpha_p,\,\alpha_p^{-1},\,\left(\frac{3}{p}\right)\,,$$
where 
$\alpha_p$ is as follows:
$$ \alpha_5=\frac{1}{ 5}(1+\sqrt{-24})\,,\quad 
\alpha_7=\frac{1}{ 7}(5+\sqrt{-24})\,,\quad 
\alpha_{11}=\frac{1}{ 11}(-7+\sqrt{-72})\,,$$
$$\alpha_{13}=\frac{1}{ 13}(-11+\sqrt{-48})\,,\quad \alpha_{17}={1}\,,\quad 
\alpha_{ 19}=\frac{1}{19}(1+\sqrt{-360})\,,$$
$$\alpha_{23}=\frac{1}{ 23}(-7+\sqrt{-480})\,,\quad 
\alpha_{29}=\frac{1}{ 29}(-4+\sqrt{-216})\,.$$
\end{cor}

\proof  Using Magma \cite{Magma}, one checks 
that
$$\begin{array}{|c|c|c|c|c|c|c|c|c|}\hline
N(5)&N(7)&N(11)&N(13)&N(17)&N(19)&N(23)&N(29)\\
\hline
27&45&107&173&323&325&515&891\\
\hline\end{array}\,\,.
$$
By Prop.~\ref{propfrobb}, the traces $t_p:=\tr\left(\rho_\ell(\Frob_p)\right),\, 3<
p\leq 29,$
are 
$$\begin{array}{|c|c|c|c|c|c|c|c|}
\hline
t_5&t_7&t_{11}&t_{13}&t_{17}&t_{19}&t_{23}&t_{29}\\
\hline
-3&3&-3&-9&17&-17&9&21\\
\hline \end{array}\,\,.
$$
Similarly, the traces $t_{p^2}:=\tr\left(\rho_\ell(\Frob_p^2)\right),\, 3<
p\leq 29,$
are 
$$\begin{array}{|c|c|c|c|c|c|c|c|}\hline 
t_{5^2}&t_{7^2}&t_{11^2}&t_{13^2}&t_{17^2}&t_{19^2}&t_{23^2}&t_{29^2}\\
\hline
-21&51&75&315&867&-357&-333&1659\\
\hline \end{array}\,\,.
$$
Since, by Poincar\'e duality (see Prop.~\ref{proppoin}), the 
image of $\rho_\ell(1)$ is contained in the orthogonal group,
the eigenvalues of $\rho_\ell(1)(\Frob_p)$ are of the form
$(\alpha_p,\pm 1,\alpha_p^{-1}).$
The 
claim now follows from solving the trace equation 
$$ \alpha_p +\alpha_p^{-1} \pm 1=\frac{t_p}{p}\,.$$ Here,
the sign $\pm 1$ coincides with 
the Legendre symbol $\left(\frac{3}{p}\right)$ which can be 
verified by checking the correctness of the equation 
$$ \alpha_p^2+\alpha_p^{-2}+1=\frac{t_{p^2}}{p^2}\,.$$
\Endproof

\subsection{Crystalline 
representations
 and some consequences of $p$-adic Hodge theory.}\label{tameinertia}

In this section, we collect some results of 
S. Wortmann \cite{Wor} on crystalline representations
and give a useful result on the image of the tame inertia 
under $\rho_\ell(1).$\\

Using the  notation of Section~\ref{secgeneralnotation}.
Let $k$ be a number field,  let $\nu \in \PP^f(k),\; \chara (\nu)=\ell,$ and
let 
$L=k_\nu.$
Let $I=I_\nu,$ let  $R=R_\nu,$ and let $I^t=I_\nu^t.$ 
Let $m \in \NN > 0$ be prime to $\ell$.
This implies that the group of  $m$-th roots of unity 
$$\mu_m:=\{ \xi \in L^{nr} \mid \xi^m=1 \}$$
is isomorphic to its image in  $\bar{\FF}_\nu,$ where 
$L^{nr}$ is the maximal unramified extension of $L.$ 
For any choice of a uniformizer $\pi$ of $L^{nr}$, the extension
$L_m:=L^{\nr}(\pi^{1/m})$ is a totally ramified extension of degree $m$.
Take
$ \sigma \in \Gal(L_m/L^{\nr})$ and set
\[ \Psi_m: \Gal(L_m/L^{\nr}) \to \mu_m,\; \Psi_m(\sigma) \pi^{1/m}=
    \sigma(\pi^{1/m})\,.\]
A character $I^t \to \FF_{\ell^n}^\times$ is called a 
{\it fundamental character of level
$n$} if it is equivalent to 
 $\Psi_{\ell^n-1}^{\ell^k}$ for some $1 \leq k \leq n.$\\

Let $V$ be a finite dimensional vector space over $L,$ let  
$\rho : G_k \to \GL(V)$  be a Galois representation, and let
$\overline{\rho}: G_k \to \GL(\overline{V})$ be the 
residual representation of $\rho.$
Denote by $D_{\cris}$ the functor that associates to $V$ the filtered 
Dieudonn{\'e} -module $D_\cris(V),$ compare to 3.1 in \cite{Wor}.
Let $\{ d_1,\ldots, d_s \}$ be the set of indices where the filtration
of $D_\cris(V)$ jumps and denote by $\alpha_i$ the multiplicity of
$d_i,$ i.e., the dimension of the associated quotient.
Then there is the following result, see  Prop. 3 in \cite{Wor}:

\begin{prop}\label{sigrid}
 Assume that $L$ is absolutely unramified.
 If $V$ is crystalline
 and if the length of the filtration of $D_{\cris}(V)$ is $<\ell,$
 then the following holds:
\begin{enumerate}

 \item The semi-simplification of the $I$-module $\overline{V}$ is well-defined
 and the action of $I$ factors through the tame quotient $I^t.$

\item For a simple  subquotient $W$ of the $I$-module $\overline{V}$ of dimension
      $h(W)$ there is an isomorphism
 $\End_{\FF_\ell}(W)$ $ \cong \FF_{\ell^{h(W)}}.$
      Fixing an isomorphism gives $W$ the structure of a one-dimensional 
       $\FF_{\ell^{h(W)}}$- vector space on which $I^t$ acts via multiplication with
      $$\Psi_{\ell^{h(W)}-1}^{i_0+ \cdots + i_{h(W)} \ell^{h(W)-1}}\, ,$$
      where the indices $-i_j$ run through $\{ d_1,\dots,d_s\}$ such that each
      component of $(d_1,\dots,d_s)$ (counted with multiplicities) appears as some
     index $-i_j$ for some subquotient.
\end{enumerate}
\end{prop}

\begin{prop}\label{Die} Let $\ell>3$ be a prime 
of $\QQ$ and let 
$$\rho_\ell: \pi_1^\et(S,\bar{1})\To \GL(\V_{\bar{1}})\simeq 
 \GL_3(\QQ_\ell)\,,$$
be as in  Section~\ref{secunderl}.
 Let $I_\ell^t\leq G_\QQ$ denote 
the tame inertia subgroup at $\ell,$ where we view 
$G_\QQ$ as a subgroup of $\pi_1^\et(S,\bar{1})$ via the splitting
$\pi_1^\et(S,\bar{1})=\pi_1^\geom(S,\bar{1})\rtimes G_\QQ.$
Then there are  the following possibilities
for the action of the tame inertia $\bar{\rho}_\ell(1)|_{I_\ell^t}:$
$$ (\Psi_{\ell-1}^{-1},\Psi_{\ell-1}^0,\Psi_{\ell-1}^1)\quad  or \quad 
 (\Psi_{\ell-1}^0, \Psi_{\ell^2-1}^{1-l})\,.$$

\end{prop}

\proof Since $\tilde{\CX}_{\bar{1}}$ is a ${\rm K}3$-surface with Picard 
number $19,$
it follows from \cite{Wor}, Prop. 45,  that the conditions 
of Prop.~\ref{sigrid} are satisfied for $k=\QQ$ and for 
$$V=\V_{\bar{1}}=\NS(\tilde{\CX}_{\bar{1}})_\ell(-1)^\perp$$ with 
$(d_1,\,d_2,\,d_3)=(0,1,2)$
and $(\alpha_1,\,\alpha_2,\,\alpha_3)=(1,1,1).$
The claim now follows from listing all possibilities 
 and from the observation
that twisting with $\chi_\ell$ raises the occurring exponents 
by one (also using that 
the image of $\bar{\rho}_\ell(1)$ is contained in the orthogonal group). 
\Endproof

\subsection{Computation of the specialization $z\mapsto 1$.} \label{secspeck}

\begin{thm}\label{thmspecki} Let $\ell\not\equiv \pm 1\mod 8$ be an odd prime.
Let   $L_\ell/\QQ(z)$ be the regular Galois extension with 
 $$\Gal(L_\ell/\QQ(z))\,=\,\SO_3(\FF_\ell)\,\simeq \,\PGL_2(\FF_\ell)$$ as in 
Cor.~\ref{cordeltal}. Let $L'_\ell/\QQ$ denote the 
specialization of $L_\ell/\QQ(z)$ under $z\mapsto 1.$ 
Then the following holds:
\begin{enumerate}
\item If, in addition, $\ell\geq 11$ and
$\left(\frac{-3}{\ell}\right)=-1\,,$ then the Galois 
group $\Gal(L_\ell/\QQ(z))$ 
is preserved under specialization $z\mapsto 1,$ i.e.,
$$\Gal(L'_\ell/\QQ)\,\simeq\, \PGL_2(\FF_\ell)$$
(where $\left( \frac{\cdot}{\cdot }\right)$ denotes the Legendre symbol).
\item
For almost all $\ell\not\equiv \pm 1,2\mod 8,$ the Galois 
group $\Gal(L_\ell/\QQ(z))$ is preserved under specialization $z\mapsto 1.$ 
\end{enumerate}

\end{thm}

\proof Let $H_\ell$ be the 
image of 
$G_\QQ$  under 
the Galois representation $$\delta_\ell:=\bar{\rho}_\ell(1)\otimes \det(\bar{\rho}_\ell(1)):\pi_1(S,\bar{1})\To \SO_3(\FF_\ell)$$ which occurs in the proof 
of Cor.~\ref{cordeltal} (where 
$G_\QQ$ is viewed as a subgroup 
of $\pi_1^\et(S,\bar{1})= \pi_1^\geom(S,\bar{1})\rtimes G_\QQ$).
 We have to show 
that $H_\ell=\SO_3(\FF_\ell)\simeq \PGL_2(\FF_\ell).$
By replacing $\delta_\ell$ by its semisimplification if necessary,
we may assume that the following possibilities
can occur for $H_\ell:$ 
$$ \ZZ/n\ZZ,\,D_n,\,A_4,\,S_4,\, A_5,\, [\SO_3(\FF_\ell),\SO_3(\FF_\ell)]\simeq 
\PSL_2(\FF_\ell),\,\SO_3(\FF_\ell)\simeq 
\PGL_2(\FF_\ell)\, ,$$
where $D_n$ denotes the dihedral group of order $2n$ and $n$ divides
$\ell-1$ or $\ell+1.$ If $\ell>3,$  Prop.~\ref{Die} implies 
that  $H_\ell$ contains a maximal torus $T_\ell,$
where $T_\ell$ denotes 
the image of the tame inertia. 

Assume first that $\ell=11.$ 
 Then the Frobenius element 
$$F_{5}:=\delta_\ell(\Frob_{5})=-\bar{\rho}_\ell(1)(\Frob_{5})$$
has order $5$ and the Frobenius element  
$$F_{13}:=\delta_\ell(\Frob_{13})=\bar{\rho}_\ell(1)(\Frob_{13})$$
has order $4.$ Since the elements $F_5$ and $F_{13}$ are contained in 
 maximal tori of different orders $10$ and $12,$ they cannot commute.
It follows that $H_\ell$ is neither a cyclic nor a dihedral group.  
The presence of the maximal torus $T_\ell$ excludes the groups
$A_4,\,S_4,\, A_5,$ because they do not contain elements of 
order $\geq 6.$  Thus $H_\ell$ contains the commutator 
subgroup of $\SO_3(\FF_\ell).$ 
Then the presence 
of a maximal torus implies that $H_\ell=\SO_3(\FF_\ell)$ (using 
the spinor norm as in the proof of Prop.~\ref{proporthog}).

By compatibility, the same arguments generalize
 to the case where $\ell\geq 11$ and
$\left(\frac{-3}{\ell}\right)=-1\,.$ This proves claim (i).

To prove claim (ii), one notes that by (i), the compatible system
$(\rho_\ell|_{G_\QQ})_{\ell}$
is not potentially induced by  Hecke characters. This means that 
the restriction $(\rho_\ell^{\rm s}|_{G_k})_\ell,$ where $k$ is any number field,
is not induced by the sum of Hecke characters. 
By \cite{Wor}, Prop. 45, the $\ell$-adic orthogonal 
complements of the Neron-Severi group 
of a ${\rm K}3$-surface of Picard number $19$ form a rank three 
orthogonal compatible system for which 
 Thm. B. of loc.~cit. holds. This results states that
if the system 
is not potentially induced by Hecke characters,
then its image contains  the commutator 
subgroup of $\SO_3(\FF_\ell)$ for almost all $\ell.$ Thus $H_\ell$ contains 
the commutator subgroup of $\SO_3(\FF_\ell)$ for almost all $\ell.$
Again, the presence
of the maximal torus $T_\ell$ implies that $H_\ell=\SO_3(\FF_\ell).$ 
\Endproof

\begin{rem} \begin{enumerate}
\item Using other  Frobenius elements (not only $F_5$ and
 $F_{13}$) it is easy to generalize Thm.~\ref{thmspecki} (i) 
to a larger set of primes. Also, the same arguments can 
be used to investigate 
the specializations of the (non-regular) Galois extensions
which appear at the primes $\ell =\pm 1 \mod 8$ via 
$\rho_\ell(1)\otimes \det(\rho_\ell(1)).$ 
\item Using more or less the same arguments, one can determine 
the specializations of $\Gal(L_\ell/\QQ(z))$ at other 
rational numbers $s_0\in \QQ\setminus \{0,\pm 2\}.$
\item Claim (ii) of Thm.~\ref{thmspecki} can be made effective once
a semi-stable alteration of ${\frak X}_1\otimes K,$ 
where $K$ is an explicitly given number field, is known (compare to  
\cite{deJong}). The  determination 
 of such an alteration seems to be an interesting
computational challenge. 
\end{enumerate}
\end{rem}

\pagebreak
                  \addcontentsline{toc}{section}{{References}}

                    \bibliographystyle{plain}
                    \bibliography{p}


\end{document}